\documentclass[11pt]{report}

\usepackage{epsf}
\usepackage{verbatim}
\usepackage{makeidx}
\usepackage{graphicx}
\usepackage{psfig}
\usepackage{amssymb}
\usepackage{amsmath}
\usepackage{latexsym}
\usepackage{amsbsy}
\usepackage{epsfig}

\usepackage{chicagob}

\usepackage{endnoteb}

\usepackage{theorem}

\usepackage{smallcaptions}

\newtheorem{theorem}{Theorem}[chapter]
\newtheorem{definition}[theorem]{Definition}
\newtheorem{example}[theorem]{Example}
\newtheorem{proposition}[theorem]{Proposition}

\newcommand{\commentout}[1]{}

\newenvironment{ownquote}
               {
\list{}{\rightmargin\leftmargin}%
                \small \item[]}
               {\endlist}

\begin{document}
\bibliographystyle{chicago}

\newpage \thispagestyle{empty} \title{A Tutorial Introduction to the \\
  Minimum Description Length Principle} \author{
  Peter Gr\"unwald \\
  Centrum voor Wiskunde en Informatica \\
  Kruislaan 413,
  1098 SJ  Amsterdam\\
  The Netherlands\\
  {\tt pdg@cwi.nl } \\
  {\tt www.grunwald.nl} } \date{} \maketitle{} \abstract{This tutorial
  provides an overview of and introduction to Rissanen's Minimum
  Description Length (MDL) Principle. 
The first chapter provides a
  conceptual, entirely non-technical introduction to the subject.  It
  serves as a basis for the technical introduction given in the second
  chapter, in which all the ideas of the first chapter are made
  mathematically precise. This tutorial will appear as the first two
  chapters in the collection {\em Advances in Minimum Description
    Length: Theory and Applications\/} \cite{GrunwaldMP04}, to be
  published by the MIT Press.  }

\setcounter{page}{2}

%\vglue-20mm
\newpage
\thispagestyle{empty}

\begin{sloppypar}
\setcounter{tocdepth}{2}
\tableofcontents
\end{sloppypar}
%\newpage
%\thispagestyle{empty}
%
%\ 
\newpage
\thispagestyle{empty}

 \newpage
 \thispagestyle{empty}

\chapter[Introducing MDL]{Introducing the MDL Principle}
  \label{chap:survey}
\newlength{\ownboxwidth} \newlength{\ownleftboxwidth}
  \newlength{\ownrightboxwidth} \setlength{\ownboxwidth}{\textwidth}
  \addtolength{\ownboxwidth}{-0.7 cm} \setlength{\ownleftboxwidth}{0.1
    cm} \setlength{\ownrightboxwidth}{0.01 cm}
 
\newcommand{\prs}{\par}
\newcommand{\emptyline}{\vspace{\baselineskip}}
\newcommand{\halfemptyline}{\vspace{.25cm}}

\newcounter{repeatproba}
\newcounter{repeatprobb}
\newcounter{repeatprobc}
\newcounter{lemmastronglawcounter}
\newcounter{snavel}

\newcounter{myenumeratecounter}
\newenvironment{myenumerate}{\begin{list}
{\arabic{myenumeratecounter}.}
{
\usecounter{myenumeratecounter}
\setlength{\itemsep}{0.0 cm}
}
}
{\end{list}}

\newcommand{\ownbox}[2]{
\noindent \ \vspace{\baselineskip} \\
\noindent
\fbox{\hspace*{\ownleftboxwidth} \parbox{\ownboxwidth}{\
\\ \normalsize
{\bf #1
} \\ \noindent #2 }\hspace*{\ownrightboxwidth}
} \ \vspace{\baselineskip} \\}

\newcommand{\myargmax}[1]{\ensuremath{\underset{#1}{\ \arg \max \ }}}
\newcommand{\myargmin}[1]{\ensuremath{\underset{#1}{\ \arg \min \ }}}
\newcommand{\isbydefinition}{\ensuremath{{:=}}}
%pdg: \isbydefinition changed 15-8-03. Was '\equiv'
\newcommand{\stsum}{\sum_{i=1}^n}
\newcommand{\transpose}{\ensuremath{{\text{\sc T}}}}
\newcommand{\indicator}{\ensuremath{{\mathbf 1}}}

\newcommand{\distance}{\ensuremath{{d}}}

\newcommand{\freq}{\ensuremath \gamma}
\newcommand{\remainder}{\ensuremath R}
\newcommand{\constant}{\ensuremath{K}}
\newcommand{\const}{\ensuremath{c}}

% SETS AND THEIR ELEMENTS 

% SETS 1: general sets

\newcommand{\reals}{{\mathbb R}}
\newcommand{\booleans}{{\mathbb B}}
\newcommand{\integers}{{\mathbb Z}}
\newcommand{\naturals}{{\mathbb N}}
\newcommand{\simplexN}[1]{\ensuremath{\Delta}^{(#1)}}

% TOPOLOGY

\newcommand{\closure}{\ensuremath{\mathrm{cl \ }}}
\newcommand{\interior}{\ensuremath{\mathrm{int \ }}}
\newcommand{\ball}{\ensuremath{B}}

% SETS 2: Outcomes and Sample Spaces.

\newcommand{\justaset}{\ensuremath {\mathcal U}}
\newcommand{\justanevent}{\ensuremath {\mathcal E}}
\newcommand{\justasymbol}{\ensuremath {s}}
% \justasymbol denoted to use an element of sample space in contexts
% where x is already used (ie Does s occur in sample x_1, ..., x_n ?)

% we use l for dimensionality of sample space, k for dimensionality of
% parameter space, and m for dimensionality of general space, if we do not
% want to commit to an interpretation (eg beginning of chapter 2)

\newcommand{\data}{\ensuremath D}
\newcommand{\Dtrain}{\ensuremath D_{\text{train}}}
\newcommand{\Dtest}{\ensuremath D_{\text{test}}}
%needed for treatment of cross-validation

\newcommand{\dataalphabet}{\ensuremath \mathcal A}
\newcommand{\samplespace}{\ensuremath {{\mathcal X}}}
%peter 22-7-03: note Xspace and samplespace now the same symbol!
\newcommand{\xspace}{\ensuremath {{\mathcal X}}}
\newcommand{\yspace}{\ensuremath {{\mathcal Y}}}
\newcommand{\zspace}{\ensuremath {{\mathcal Z}}}
\newcommand{\outcomex}{\ensuremath x}
\newcommand{\outcomey}{\ensuremath y}
\newcommand{\bogus}{\ensuremath \Box}
\newcommand{\maxoutcome}{\ensuremath m}
\newcommand{\region}{\ensuremath{{\mathbf R}}}

%SETS AND ELEMENTS 3: estimators/hypotheses/parameters

%\newcommand{\mldim}[1]{\ensuremath{\hat{\theta}_{#1}}}
\newcommand{\mlN}[1]{\ensuremath{\hat{\theta}^{(#1)}}}
\newcommand{\ml}{\ensuremath \hat{\theta}}
\newcommand{\mapest}{\ensuremath \breve{\theta}}
\newcommand{\param}{\ensuremath {\theta}}
\newcommand{\paramN}[1]{\ensuremath {\theta}^{(#1)}}
\newcommand{\meanest}{\ensuremath \bar{\theta}}
\newcommand{\truepar}{\ensuremath \theta^*}
\newcommand{\trueprob}{\ensuremath \prob^*}

\newcommand{\twop}[1]{\ensuremath{{\ddot{#1}}}} 
\newcommand{\twopN}[2]{\ensuremath{{\ddot{#1}^{(#2)}}}}
\newcommand{\mcparam}[1]{\ensuremath{\theta_{[1|#1]}}}
\newcommand{\mcmlparam}[1]{\ensuremath{\hat{\theta}_{[1|#1]}}}
\newcommand{\mctrueparam}[1]{\ensuremath{{\theta}^*_{[1|#1]}}}
\newcommand{\mcmlparamb}[2]{\ensuremath{\hat{\theta}_{[#1|#2]}}}
\newcommand{\utm}{\ensuremath{\text{\sc{ul}}}}

\newcommand{\phyp}{\ensuremath {H}}
\newcommand{\ghyp}{\ensuremath {\mathcal H}}
%phyp = point hypothesis
%ghyp = general (point or compound) hypothesis
\newcommand{\hfun}{\ensuremath {h}}

%SETS AND ELEMENTS 4: parameter SETS

%\newcommand{\parasetN}[1]{\ensuremath{\Theta_{#1}}}
\newcommand{\parasetN}[1]{\ensuremath{\Theta^{(#1)}}}
\newcommand{\paraset}{\ensuremath{\Theta}}
\newcommand{\twopparaset}[2]{\ensuremath{{\ddot{\Theta}^{(#1)}_{#2}}}}

% codes and codelengths

\newcommand{\clengthnr}{\ensuremath{l}}

\newcommand{\lunif}{\ensuremath{L_{U}}}
\newcommand{\punif}{\ensuremath{P_{U}}}
\newcommand{\Lint}{\ensuremath{L_{\naturals}}}
\newcommand{\punivint}{\ensuremath{\bar{P}_{\naturals}}}
\newcommand{\indexcode}{\ensuremath C_{\text{index}}}
\newcommand{\shannoncode}{\ensuremath C_{\text{Shannon}}}
\newcommand{\Lindex}{\ensuremath L_{\text{index}}}
\newcommand{\Lshannon}{\ensuremath L_{\text{Shannon}}}
\newcommand{\prefixlengths}{\ensuremath{\mathcal{L}}}
\newcommand{\code}{\ensuremath C}
\newcommand{\CN}[1]{\ensuremath C^{(#1)}}
\newcommand{\codelength}{\ensuremath{L}}
\newcommand{\codelengths}{\ensuremath{\mathcal{L}}}
\newcommand{\regret}{\ensuremath{\mathcal R}}
\newcommand{\maxregret}{\ensuremath{\mathcal R}_{\max}}

% functions, vectors, matrices

\newcommand{\justafunction}{\ensuremath{g}}
\newcommand{\justavector}{\ensuremath{v}}
\newcommand{\justamatrix}{\ensuremath{M}}
% used when we do not want to enforce the interpretation of a parameter space etc.

% probability things

\newcommand{\empprob}{\ensuremath{\mathbb P}}

\newcommand{\prob}{\ensuremath P}
\newcommand{\pd}{\ensuremath P}
\newcommand{\densityf}{\ensuremath{f}}

\newcommand{\Exp}{\ensuremath{\text{\rm E}}}
\newcommand{\cov}{\ensuremath{\text{\rm cov}}}
\newcommand{\var}{\ensuremath{\text{\rm var}}}
\newcommand{\noiseterm}{\ensuremath{Z}}

\newcommand{\expprior}{\ensuremath{q}}

%models 

\newcommand{\Bernoulli}{\ensuremath{\mathcal{B}}}
\newcommand{\BernoulliN}[1]{\ensuremath{{\mathcal{B}}^{(#1)}}}
\newcommand{\fakemarkov}{\ensuremath{\dddot{\Bernoulli}}}
\newcommand{\M}{\ensuremath {{\cal M}}}
\newcommand{\cH}{\ensuremath {{\cal H}}}
\newcommand{\MN}[1]{\ensuremath {{\cal M}^{(#1)}}}
\newcommand{\HN}[1]{\ensuremath {{\cal H}^{(#1)}}}
\newcommand{\cM}{\ensuremath {{\cal C}}}
  %these are 'candidate hypotheses'

% Fisher inf stuff

\newcommand{\observedmatrix}{\ensuremath \hat{I}}
%H already used for entropy
\newcommand{\fishermatrix}{\ensuremath I}
\newcommand{\curvelength}[1]{\ensuremath{\int_{#1} \fisherdens d\theta }}
\newcommand{\fisherdens}{\ensuremath{\sqrt{| \fishermatrix(\theta)|}}}
\newcommand{\fishervolume}{\ensuremath{V_{\text{d}}}}

% CHAPTER 3 

\newcommand{\level}{\mbox{\sc level}}
\newcommand{\descendantset}{\mbox{\sc d}}
\newcommand{\precision}{\ensuremath{d}}

\newcommand{\kl}{\mbox{$D$}}
\newcommand{\entropy}{{\ensuremath{\text{\rm H}}}}

% CHAPTER 4 

% CHAPTER 6, 7

\newcommand{\jp}{\ensuremath{w_{\text{Jef}}}}
\newcommand{\Lref}{\ensuremath{\mathcal{L}_{\text{ref}}}}
\newcommand{\pmodel}{\ensuremath{\theta}}
\newcommand{\prior}{\ensuremath{w}}
\newcommand{\discreteprior}{\ensuremath{W}}
\newcommand{\decision}{\ensuremath{\delta}}
\newcommand{\decisionspace}{\ensuremath{\Delta}}
\newcommand{\loss}{\ensuremath{\text{\sc loss}}}
\newcommand{\emptysequence}{\ensuremath{\text{\sc doen}}}
\newcommand{\realprob}{\ensuremath{\text{\sc doen}}}
\newcommand{\popt}{\ensuremath{\text{\sc doen}}}
\newcommand{\opthypothesis}{\ensuremath{\text{\sc doen}}}
%MOET ANDERS!

\newcommand{\puniv}{\ensuremath{\bar{P}}}
\newcommand{\punivtp}{\ensuremath{\bar{P}_{\text{2-p}}}}
\newcommand{\punivplugin}{\ensuremath{\bar{P}_{\text{plug-in}}}}
\newcommand{\punivcrude}{\ensuremath{\bar{P}_{\text{crude}}}}
\newcommand{\punivref}{\ensuremath{\bar{P}_{\text{refined}}}}
\newcommand{\punivseq}{\ensuremath{\bar{P}_{\text{\rm seq}}}}
\newcommand{\punivseqn}[1]{\ensuremath{\bar{P}_{\text{seq},#1}}}
\newcommand{\punivopt}{\ensuremath{\bar{P}_{\text{\rm opt}}}}
\newcommand{\optprior}{\ensuremath{{\prior}_{\text{\rm opt}}}}
\newcommand{\punivrnml}{\ensuremath{\bar{P}_{\text{\rm rnml}}}}
\newcommand{\punivnml}{\ensuremath{\bar{P}_{\text{\rm nml}}}}
\newcommand{\punivmeta}{\ensuremath{\bar{P}_{\text{\rm meta}}}}
\newcommand{\punivpreq}{\ensuremath{\bar{P}_{\text{\rm preq}}}}
\newcommand{\punivnmln}[1]{\ensuremath{\bar{P}_{\text{\rm nml},#1}}}
\newcommand{\punivbayes}{\ensuremath{\bar{P}_{\text{\rm Bayes}}}}
\newcommand{\luniv}{\ensuremath{\bar{L}}}
\newcommand{\lunivtp}{\ensuremath{\bar{L}_{\text{\rm 2-p}}}}
\newcommand{\lunivnml}{\ensuremath{\bar{L}_{\text{\rm nml}}}}
\newcommand{\lunivbayes}{\ensuremath{\bar{L}_{\text{\rm Bayes}}}}

\newcommand{\volume}{\ensuremath{V_{\text{p}}}}
\newcommand{\Mdisc}{\ensuremath{\ddot{\mathcal M}}}
\newcommand{\probset}{\ensuremath{{\mathcal P}}}
\newcommand{\complexity}{\ensuremath{\text{\bf COMP}}}
\newcommand{\rcomplexity}{\ensuremath{\text{\bf RCOMP}}}
\newcommand{\typset}{\ensuremath{{\mathcal A}}}

%----------------------------------------------------------
% POST chapter 7

%MML

\newcommand{\wftwopartlength}{\ensuremath{L_{\ensuremath{\text{mml-wf}}}}}
\newcommand{\mmlplength}{\ensuremath{L_{\ensuremath{\text{mml-p}}}}}

%MaxEnt
\newcommand{\constraint}{\ensuremath{\mathcal C}}
\newcommand{\constraintdataset}{\ensuremath {{\mathbf C}^n}}
\newcommand{\entropified}[1]{\ensuremath{{\langle {#1} \rangle}}}
\newcommand{\function}{\ensuremath \phi}
\newcommand{\range}{\ensuremath \mathbf U}
\newcommand{\expval}{\ensuremath t}
\newcommand{\mep}{\ensuremath{\prob_{me}}}
\newcommand{\mept}{\ensuremath{\prob_{me,\phi}(\cdot | t)}}
\newcommand{\freqdataset}{\ensuremath {{\mathbf G}^n}}
\newcommand{\freqdatasetn}[1]{\ensuremath {{\mathbf G}^{#1}}}
\newcommand{\constraintprobset}{\ensuremath {{\mathbf M}}}
\newcommand{\constraintdatasetn}[1]{\ensuremath {{\mathbf C}^{#1}}}
\newcommand{\empiricalconstraint}{\ensuremath {{\cal C}_{ e}}}
\newcommand{\average}[2]{\ensuremath \overline{#1}^{#2}}

% Unclear wether needed

%\newcommand{\maxnumclass}{\ensuremath K}
%\newcommand{\numcomps}{\ensuremath m}
%\newcommand{\numdimensions}{\ensuremath k}
%\newcommand{\traindata}{\ensuremath{D}}
\newcommand{\lagrange}{\ensuremath \beta}

% priors

%\newcommand{\jprior}{\ensuremath{\pi}}

%old

%\index{symbols}{zzz@$\tilde{c}$ where $c \in \{
%\pmodel,\lagrange,\hypothesis,\sigma \}$: model that minimizes expected error}
%\index{symbols}{zzz@$\hat{c}$ where $c \in \{
%\pmodel,\lagrange,\hypothesis,\sigma \}$: model that minimizes
%empirical error/maximizes likelihood}
%\index{symbols}{zzz@${c}_{mdl}$ where $c \in \{
%\pmodel,\lagrange,\hypothesis,\sigma \}$: model that minimizes
%two-part codelength}
%\index{symbols}{zzz@$\ddot{c}$ where $c \in \{
%\pmodel,\lagrange,\hypothesis \}$: (unspecified) estimator}
%\index{symbols}{zzz@$\breve{c}$ where $c \in \{
%\pmodel,\lagrange,\hypothesis \}$: Bayesian MAP estimator}
%\index{symbols}{zzz@$\bar{c}$ where $c \in \{
%\pmodel,\lagrange,\hypothesis \}$: mean of Bayesian posterior}
%\index{symbols}{zzz@$c^*$ where $c \in \{
%\pmodel,\lagrange,\hypothesis \}$: `true' model generating the data}

%\newcommand{\convergesatrate}{\ensuremath{\lesssim}}

% use next line ONLY if you have three or more authors; otherwise remove it.
%\clearpage

% now follows the abstract (to be replaced by your own abstract).

%you can replace the following sections with the body of your .tex file

\section{Introduction and Overview}
\label{sec:introduction}
How does one decide among competing explanations of data given limited
observations? This is the problem of {\em model selection}. It stands
out as one of the most important problems of inductive and statistical
inference.  The Minimum Description Length (MDL) Principle is a
relatively recent method for inductive inference that provides a
generic solution to the model selection problem. MDL is based on the
following insight: 
any regularity in the data can be used to {\it compress\/} the data,
i.e. to describe it using fewer symbols than the number of symbols
needed to describe the data literally. The more regularities there
are, the more the data can be compressed. Equating `learning'
with `finding regularity', we can therefore say that the more we are
able to compress the data, the more we have {\it learned\/} about the
data. Formalizing this idea leads to a general theory of inductive
inference with several attractive properties:
\begin{description}
\item[1. Occam's Razor] MDL chooses  a
  model that trades-off goodness-of-fit on the observed data with
  `complexity' or `richness' of the model. As such, MDL embodies a
  form of Occam's Razor, a principle that is both intuitively
  appealing and informally applied throughout all the sciences.
\item[2. No overfitting, {\em automatically\/}]
MDL procedures {\em automatically\/} and {\em inherently\/} protect
  against overfitting and can be used to estimate both the parameters
  and the structure (e.g., number of parameters) of a model. 
In contrast, to avoid overfitting when estimating the structure of a
model, 
traditional methods such as maximum likelihood must be
  {\em modified\/} and extended 
with additional, typically {\em ad hoc\/} principles. 
\item[3. Bayesian interpretation] MDL is closely related to Bayesian
  inference, but avoids some of the interpretation difficulties of the
  Bayesian approach\endnote{See Section~\ref{sec:bayes},
    Example~\ref{ex:bayescraze}.}, especially in the realistic case
  when it is known a priori to the modeler that none of the models
  under consideration is true. In fact:
\item[4. No need for `underlying truth'] In contrast to other statistical
  methods, MDL procedures have a clear interpretation independent of
  whether or not there exists some underlying `true' model.
\item[5. Predictive interpretation] Because data compression is formally
  equivalent to a form of probabilistic prediction, MDL methods
  can be interpreted as searching for a model with good predictive
  performance on {\em unseen\/} data.
\end{description}
In this chapter, we introduce the MDL Principle in an entirely
non-technical way, concentrating on its most important applications,
model selection and avoiding overfitting. In Section~\ref{sec:basic}
we discuss the relation between learning and data compression.
Section~\ref{sec:modselintro} introduces model selection and outlines
a first, `crude' version of MDL that can be applied to model
selection. Section~\ref{sec:cruderefined} indicates how these `crude'
ideas need to be refined to tackle small sample sizes and differences
in model complexity between models with the same number of parameters.
Section~\ref{sec:philosophy} discusses the philosophy underlying MDL,
and considers its relation to Occam's Razor. Section~\ref{sec:history}
briefly discusses the history of MDL. All this is summarized in
Section~\ref{sec:summary}.
\section{The Fundamental Idea: \\ Learning as Data Compression}
\label{sec:basic}
We are interested in developing a method for {\em learning\/} the laws
and regularities in data. The following example will illustrate what
we mean by this and give a first idea of how it can be related to
descriptions of data.
\paragraph{Regularity $\ldots$}
\label{ex:bernoulli}
Consider the following three
sequences. We assume that each sequence is 10000 bits long, and we
just list the beginning and the end of each sequence.
\begin{eqnarray}
\label{seqa}
& {\tt
00010001000100010001} \ \ \ldots \ \ {\tt 0001000100010001000100010001\ } & \\
\label{seqb}
& {\tt
01110100110100100110} \ \ \ldots \ \ {\tt 1010111010111011000101100010\ } & \\
\label{seqc}
& {\tt
00011000001010100000} \ \ \ldots \ \ {\tt 0010001000010000001000110000\ }&
\end{eqnarray}
The first of these three sequences is a 2500-fold repetition of {\tt
  0001}. Intuitively, the sequence looks regular; there seems to be a
simple `law' underlying it; it might make sense to conjecture that
future data will also be subject to this law, and to predict that
future data will behave according to this law.  The second sequence
has been generated by tosses of a fair coin. It is intuitively
speaking as `random as possible', and in this sense there is no
regularity underlying it. Indeed, we cannot seem to find such a
regularity either when we look at the data.  The third sequence
contains approximately four times as many 0s as 1s. It looks less regular,
more random than the first; but it looks less random than the second.
There is still some discernible regularity in these data, but of a
statistical rather than of a deterministic kind. Again, noticing that
such a regularity is there and predicting that future data will behave
according to the same regularity seems sensible.
\paragraph{...and Compression}
We claimed that any
regularity detected in the data can be used to {\it compress\/} the data,
i.e. to describe it in a short manner. 
%Such a description should
%always completely determine the data it describes\endnote{But see Section~\ref{sec:kolmogorov.}} - hence given a
%description or {\it encoding\/} $D'$ of a particular sequence of data
%$D$, we should always be able to fully reconstruct $D$ on the basis of
%$D'$. 
Descriptions are always relative to some {\it description
method\/} which maps descriptions $D'$ in a unique manner to data sets
$D$. A particularly versatile description method is a general-purpose
computer language like {\sc C} or {\sc Pascal}. A description of 
$D$ is then any computer program that prints $D$ and then halts. Let
us see whether our claim works for the three sequences above. Using
a language similar to Pascal, we can write a program $$ {\tt for\ i\
=\ 1\ to\ 2500;\ print\ `0001`;\ next;\ halt } $$ which prints
sequence (1) but is clearly a lot shorter. Thus, 
sequence (1) is indeed highly compressible. On
the other hand, we show in Section~\ref{sec:prob}, that if one
generates a sequence like (2) by tosses of a fair coin, then with
extremely high probability, the shortest program that prints (2) and
then halts will look something like this: $$ {\tt print\
`011101001101000010101010........1010111010111011000101100010`;\ halt }
$$ This program's size is about equal to the length of the
sequence. Clearly, it does nothing more than repeat the sequence.

The third sequence lies in between the first two: generalizing $n =
10000$ to arbitrary length $n$, we show in Section~\ref{sec:prob} that the
first sequence can be compressed to $O(\log n)$ bits; with
overwhelming probability, the second sequence cannot be compressed at all;
and the third sequence can be compressed to some length $\alpha n$,
with $0 < \alpha < 1$.
\begin{example}{\bf [compressing various regular sequences]}
  \rm The regularities underlying sequences (1) and (3) were of a very
  particular kind. To illustrate that {\em any\/} type of regularity
  in a sequence may be exploited to compress that sequence, we give a
  few more examples:
\begin{description}
\item[The Number $\pi$] Evidently, there exists a computer program for generating
  the first $n$ digits of $\pi$ -- such a program could be based, for
  example, on an infinite series expansion of $\pi$. 
This computer program has constant
  size, except for the specification of $n$ which takes no more than 
$O(\log n)$ bits. Thus, when $n$ is very large, the size of the
program generating the first $n$ digits of $\pi$ will be very small
compared to $n$: the $\pi$-digit sequence is deterministic, and therefore
extremely regular.
\item[Physics Data] Consider a two-column table where the first column
  contains numbers representing various heights
  from which  an object was dropped. The second column
  contains the corresponding times it
  took for the object to reach the ground. Assume both heights and
  times are recorded to some finite precision. In
  Section~\ref{sec:modselintro} we illustrate that such a table can be
  substantially compressed by first describing the coefficients of the
  second-degree polynomial $H$ that expresses Newton's law; 
then describing the heights; and then describing the deviation of the
time points from the numbers predicted by $H$.
\item[Natural Language] Most sequences of words are not valid sentences according
  to the English language. This fact can be exploited to
  substantially compress English text, as long as it is syntactically
  mostly correct: by first describing a
  grammar for English, and then describing an English text $D$ 
with the help
  of that grammar \cite{Grunwald96d}, $D$ can be described using
  much less bits than are needed without the assumption that word order
  is constrained.
\end{description}
\end{example}
\subsection{Kolmogorov Complexity and Ideal MDL} 
\label{sec:kolmo}
To formalize our ideas, we need to decide on a description method,
that is, a formal language in which to express properties of the data.
The most general choice is a general-purpose\endnote{By this we mean
  that a universal Turing Machine can be implemented in it
  \cite{LiV97}.} computer language such as {\sc C} or {\sc {Pascal}}.
This choice leads to the definition of the {\it Kolmogorov
  Complexity\/} \cite{LiV97} of a sequence as the length of the
shortest program that prints the sequence and then halts. The lower
the Kolmogorov complexity of a sequence, the {\it more regular} it is.
This notion seems to be highly dependent on the particular computer
language used. However, it turns out that for every two
general-purpose programming languages $A$ and $B$ and every data
sequence $D$, the length of the shortest program for $D$ written in
language $A$ and the length of the shortest program for $D$ written in
language $B$ differ by no more than a constant $c$, which does not
depend on the length of $D$. This so-called {\em invariance theorem\/}
says that, {\em as long as the sequence $D$ is long enough}, it is not
essential which computer language one chooses, as long as it is
general-purpose. Kolmogorov complexity was introduced, and the
invariance theorem was proved, independently by \citeN{Kolmogorov65},
\citeN{Chaitin69} and \citeN{Solomonoff64}. Solomonoff's paper, called
{\em A Theory of Inductive Inference}, contained the idea that the
ultimate model for a sequence of data may be identified with the
shortest program that prints the data. Solomonoff's ideas were later
extended by several authors, leading to an `idealized' version of MDL
\cite{Solomonoff78,LiV97,GacsTV01}. This idealized MDL is very general
in scope, but not practically applicable, for the following two
reasons:
\begin{description}
\item[1. uncomputability] It can be shown that there exists no computer program that,  for every set of data $D$, when given $D$ as
input, returns the  shortest program that prints $D$ \cite{LiV97}.
\item[2. arbitrariness/dependence on syntax] In practice we are
  confronted with small data samples for which the invariance theorem
  does not say much. Then the hypothesis chosen by idealized MDL may
  depend on arbitrary details of the syntax of the programming
  language under consideration.
\end{description}
\subsection{Practical MDL} 
Like most authors in the field, we concentrate here
on non-idealized,
practical versions of  MDL that explicitly deal with the two problems
mentioned above. The basic idea is to scale down Solomonoff's
approach so that it does become applicable. This is achieved by using description methods
that are less expressive than general-purpose computer languages. Such
description methods $C$ should be restrictive enough so that for any data
sequence $D$, we can always compute the length of the shortest
description of $D$ that is attainable using method $C$; but they
should be 
general enough to allow us to compress many of the intuitively
`regular' sequences. The price we pay is that, using the `practical'
MDL Principle, there will always be some regular sequences which we
will not be able to compress. But we already know that there can be
{\it no \/} method for inductive inference at all which will always
give us all the regularity there is --- simply because there can be no
automated method which for any sequence $D$ finds the shortest
computer program that prints $D$ and then halts.  Moreover, it will
often be possible to guide a suitable choice of $C$ by a priori
knowledge we have about our problem domain. For example, below we
consider a description method $C$ that is based on the class of all
polynomials, such that with the help of $C$ we can compress all data
sets which can meaningfully be seen as points on some polynomial.
\section{MDL and Model Selection} 
\label{sec:modselintro}
\begin{figure}
\centerline{\epsfxsize=6cm \epsfbox{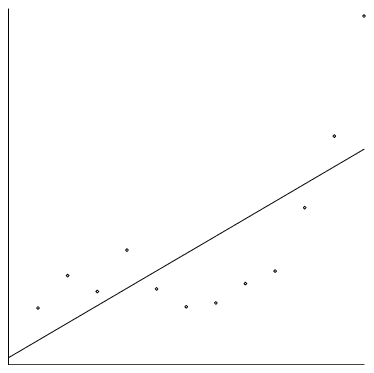} \hspace*{-2 cm}
\epsfxsize=6cm \epsfbox{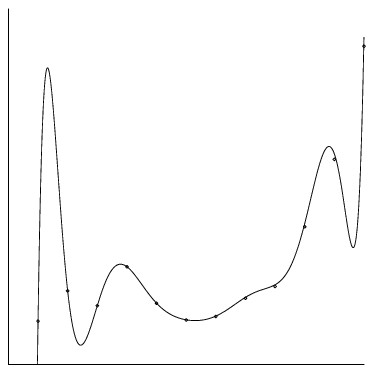} \epsfxsize=6cm
\vspace*{-1 cm} \hspace*{ -2 cm} \epsfbox{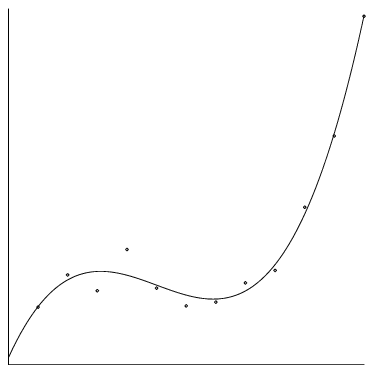}}
\caption{\label{fig:polynomial} A simple, a complex and a
trade-off (3rd degree) polynomial.}
\end{figure}

Let us recapitulate our main insights so far:
\ownbox{MDL: The Basic Idea}{The goal of statistical inference may be
  cast as trying to find regularity in the data.  `Regularity' may be
  identified with `ability to compress'. MDL combines these two
  insights by {\em viewing learning as data compression\/}: it tells
  us that, for a given set of hypotheses ${\cal H}$ and data set $D$,
  we should try to find the hypothesis or combination of hypotheses in
  ${\cal H}$ that compresses $D$ most.}
\noindent
This idea can be applied to all sorts of
inductive inference problems, but it turns out to be most fruitful in (and its
development has mostly concentrated on) problems of {\em model
  selection\/} and, more generally, dealing with {\em overfitting}.
Here is a standard example (we explain the difference between `model'
and `hypothesis' after the example).
\begin{example}{\bf [Model Selection and Overfitting]}
\label{ex:modsel}
\rm  Consider the points in Figure~\ref{fig:polynomial}. 
We would like to learn how the $y$-values depend on the
  $x$-values. To this end, we may want to 
fit a polynomial to the points. Straightforward linear regression will give us
  the leftmost polynomial - a straight line that seems overly simple:
  it does not capture the regularities in the data well. Since for any
  set of $n$ points there exists a polynomial of the $(n-1)$-st degree
  that goes exactly through all these points, simply looking for the
  polynomial with the least error will give us a polynomial like the
  one in the second picture. This polynomial seems overly complex: it
  reflects the random fluctuations in the data rather than the general
  pattern underlying it.  Instead of picking the overly simple or the
  overly complex polynomial, it seems more reasonable to prefer a
relatively  simple polynomial with small but nonzero error, as in the rightmost
  picture. This intuition is confirmed by numerous experiments on
  real-world data from a broad variety of sources
  \cite{Rissanen89,Vapnik98,Ripley96}: if one naively fits a
  high-degree polynomial to a small sample (set of data points), then one
  obtains a very good fit to the data. Yet if one {\em tests\/}
  the inferred polynomial on a second set of data coming from the same
  source,  it typically fits this test data very badly in the
  sense that there is a large distance between the polynomial and the
  new data points.  We say
  that the polynomial {\em overfits\/} the data. 
%The problem of how to
%  avoid  overfitting is  one of the hardest  in
%  inductive inference.
%
%The goal of {\em model selection\/} is to infer a `good' model for the
%available data. The model should neither be overly simple (in this
%case it `underfits' and
%thus fails to capture some meaningful patterns in the data) nor overly
%complex (in which case it `overfits' and fits the noise in the data
%rather than the underlying regularities).
  Indeed, all model selection methods that are used in practice either
  implicitly or explicitly choose a trade-off between goodness-of-fit
  and complexity of the models involved.  In practice, such trade-offs
  lead to much better predictions of test data than one would get by
  adopting the `simplest' (one degree) or most
  `complex\endnote{Strictly speaking, in our context it is not very
    accurate to speak of `simple' or `complex' polynomials; instead we
    should call the {\em set\/} of first degree polynomials `simple',
    and the {\em set\/} of 100-th degree polynomials `complex'.}'
  ($n-1$-degree) polynomial. MDL provides
  one particular means of achieving such a trade-off.
\end{example}
It will be useful to make a precise distinction between `model' and
`hypothesis': 
\ownbox{Models vs. Hypotheses}{
We use the phrase {\em point hypothesis\/}
to refer to a {\em single\/} probability distribution or function. An
example is the polynomial $5x^2 + 4x + 3$. A point hypothesis is also
known as a `simple hypothesis' in the statistical literature.

We use the word {\em model\/} to refer to a family (set) of probability distributions or functions
with the same functional form. An example is the set of all
second-degree polynomials. A model is also known as a `composite
hypothesis' in the statistical literature.

We use {\em hypothesis\/} as a generic term, referring to both point
hypotheses and models.
}
In our terminology, 
the problem described in Example~\ref{ex:modsel} is a `hypothesis
selection problem' if we are interested in selecting both the degree
of a polynomial and the corresponding parameters; it is a `model
selection problem' if we are mainly interested in selecting the degree. 

\ 
\\
\noindent
To apply MDL to polynomial or other types of hypothesis and model
selection, we have to make precise the somewhat vague insight
`learning may be viewed as data compression'.  This can be done in
various ways.  In this section, we concentrate on the earliest and
simplest implementation of the idea. This is the so-called {\em
  two-part code\/} version of MDL:
\begin{figure}[h]
  \ownbox{Crude\endnote{The terminology `crude MDL' is not standard.
      It is introduced here for pedagogical reasons, to make clear the
      importance of having a single, unified principle for designing
      codes.  It should be noted that Rissanen's and Barron's early
      theoretical papers on MDL already contain such principles,
      albeit in a slightly different form than in their recent papers.
      Early practical applications \cite{QuinlanR89,Grunwald96d} often
      do use {\em ad hoc\/} two-part codes which really are `crude' in
      the sense defined here.}, 
Two-part
  Version of MDL Principle (Informally Stated)}{Let $\HN{1}, \HN{2},
  \ldots$ be a list of candidate models (e.g., $\HN{k}$ is the set of
  $k$-th degree polynomials), each containing a set of
  point hypotheses (e.g., individual polynomials). The best point hypothesis
  $H \in \HN{1} \cup \HN{2} \cup \ldots $ to explain the data $D$
  is the one which minimizes the sum $L(\phyp) + L(D|\phyp)$, where
\label{box:hypsel}
\begin{itemize}
\item $L(\phyp)$ is the length, in bits, of the description of the
  hypothesis; and
\item $L(D|\phyp)$ is the length, 
in bits, of the description of the data when encoded
with the help of the hypothesis.
\end{itemize}
The best {\em model\/} to explain $D$ is the smallest model
containing the selected $\phyp$.} \vspace*{-0.5 cm}
%\caption{\label{fig:twopart} The two-part MDL Principle: first, crude
%  implementation of the MDL ideas.}
\end{figure}

\noindent
%This `crude\endnote{See the previous note.}' version of MDL may perform
%suboptimally in practice and is treated here mainly to illustrate the
%basic ideas in a simple manner. Refined versions are discussed in the
%next section.
%We make the meaning of $L(\phyp)$
%and $L(D| \phyp)$ precise in Section~\ref{sec:crude}.
%Intuitively, a hypothesis with a large description
%length (large $L(\phyp)$) is what we would intuitively call a
%`complex' one\endnote{See previous note.}.
\begin{example}{ \bf [Polynomials, cont.]}
\rm
In our previous example, the candidate hypotheses were polynomials. We can
describe a polynomial by describing its coefficients in a certain
precision (number of bits per parameter). Thus, the higher the degree
of a polynomial or the precision, the more\endnote{See the previous note.}  bits we
need to describe it and the more `complex' it becomes. A description
of the data `with the help of' a hypothesis means that the better the
hypothesis fits the data, the shorter the description will be. A
hypothesis that fits the data well gives us a lot of {\em
  information\/} about the data.  Such information can always be used
to compress the data (Section~\ref{sec:prob}).  Intuitively, this is
because we only have to code the {\em errors\/} the hypothesis makes
on the data rather than the full data. In our polynomial example, the
better a polynomial $\phyp$ fits $D$, the fewer bits we need
to encode the discrepancies between the actual $y$-values $y_i$ and
the predicted $y$-values $H(x_i)$. 
%In general, the better the fit, the
%fewer bits we need to describe the data with the help of the
%hypothesis (i.e. the smaller $L(D|\phyp)$).  
We can typically find a
very complex point hypothesis (large $L(\phyp)$) with a very good fit (small
$L(D|\phyp)$).  We can also typically find a very simple point hypothesis
(small $L(\phyp)$) with a rather bad fit (large $L(D|\phyp)$). The sum
of the two description lengths will be minimized at a
hypothesis that is quite (but not too) `simple', with a good (but not
perfect) fit.
% As desired, the MDL Principle gives us a trade-off
%between hypothesis complexity and goodness-of-fit on the data
%sequence.  
\end{example}
\section{Crude and Refined MDL}
\label{sec:cruderefined}
Crude MDL picks the $H$ minimizing the sum $L(H) + L(D| H)$. To make
this procedure well-defined, we need to agree on
precise definitions for the codes (description methods) giving rise to lengths $L(D|H)$ and
$L(H)$. We now discuss these codes in more detail. We will see that
the definition of $L(H)$ is problematic, indicating that we somehow
need to `refine' our crude MDL Principle.

\paragraph{Definition of $L(D|H)$}
Consider a two-part code as described above, and assume for the time
being that all $H$ under consideration define probability
distributions. If $H$ is a polynomial, we can turn it into a
distribution by making the additional assumption that the $Y$-values
are given by $Y = H(X) + Z$, where $Z$ is a
normally distributed noise term.

For each $H$ we need to define a code with length $L(\cdot \mid H)$
such that $L(D|H)$ can be interpreted as `the codelength of $D$ when
encoded with the help of $H$'. 
It turns out that for probabilistic hypotheses, there
is only one reasonable choice for this code. It is the so-called {\em Shannon-Fano
  code}, satisfying, for all data sequences $D$, $L(D|H) = - \log
P(D|H)$, where $P(D|H)$ is the probability mass or density of $D$
according to $H$ -- such a code always exists, Section~\ref{sec:prob}.
\paragraph{Definition of $L(H)$: A Problem for Crude MDL}
It is more problematic to find a good code for hypotheses $H$. Some
authors have simply used `intuitively reasonable' codes in the past,
but this is not satisfactory: since the description length $L(H)$ of
any fixed point hypothesis $H$ can be very large under one code, but quite
short under another, our procedure is in danger of becoming arbitrary.
Instead, we need some additional principle for designing a code for
${\cal H}$.  In the first publications on MDL
\cite{Rissanen78,Rissanen83}, it was advocated to choose some sort of
{\em minimax code\/} for ${\cal H}$, minimizing, in some precisely
defined sense, the shortest worst-case total description length
$L(H) + L(D|H)$, where the worst-case is over all possible data
sequences. Thus, the MDL Principle is employed at a `meta-level' to
choose a code for $H$. However, this code requires a cumbersome
discretization of the model space ${\cal H}$, which is not always
feasible in practice. Alternatively, Barron \citeyear{Barron85}
encoded $H$ by the shortest computer 
program that, when input $D$, computes
$P(D|H)$. While it can be shown that this leads to similar
codelengths, it is computationally problematic. Later, Rissanen
\citeyear{Rissanen84} realized that these problems could be side-stepped
by using a {\em one-part\/} rather than a {\em two-part code}. This
development culminated in 1996 in a completely precise prescription of
MDL for many, but certainly not all practical situations
\cite{Rissanen96}. We call this modern version of MDL {\em refined
  MDL\/}:
\paragraph{Refined MDL}
In refined MDL, we associate a code for encoding $D$ {\em
  not with a single $H \in {\cal H}$}, but with the full  model
${\cal H}$. Thus, given model ${\cal H}$, we encode data not in two parts but we
design a single {\em one-part code\/} with lengths $\luniv(D| {\cal H})$.
This code is designed such that {\em whenever there is a member of
  (parameter in) ${\cal H}$ that fits the data well, in the sense  that $L(D
  \mid H)$ is small,  then the codelength $\luniv(D | {\cal H})$ will
  also be small}. Codes with this property are called {\em
  universal codes\/} in the information-theoretic literature \cite{BarronRY98}. 
Among all such
universal codes, we pick the one that is {\em minimax optimal\/}
in a sense made precise in Section~\ref{sec:universal}.  For example, the set
${\cal H}^{(3)}$ of third-degree polynomials is associated 
with a code with lengths
$\luniv(\cdot \mid {\cal H}^{(3)})$ 
such that, the better the data $D$ are fit by
the best-fitting third-degree polynomial, the shorter the codelength
$\luniv(D\mid {\cal H})$. $\luniv(D \mid {\cal H})$ is  called the  {\em
  stochastic complexity\/} of the data given the model.

\paragraph{Parametric Complexity}
The second fundamental concept of refined MDL 
is the {\em parametric complexity\/} of
a parametric model ${\cal H}$ which we denote by $\complexity({\cal H})$.
This is a measure of the `richness' of model ${\cal H}$, indicating its
ability to fit random data. This complexity is related to the
degrees-of-freedom in ${\cal H}$, but also to the geometrical structure of
${\cal H}$; see Example~\ref{ex:modselb}.
To see how it relates to stochastic complexity, let, for given data
$D$, $\hat{H}$ denote the distribution in ${\cal H}$ which maximizes
the probability, and hence minimizes the codelength $L(D\mid \hat{H})$
of $D$. It turns out that 
%\begin{equationn}
$$
\text{stochastic complexity of $D$ given ${\cal H}$} = L(D \mid \hat{H}) +
\complexity({\cal H}).
$$
Refined MDL model selection between two parametric models (such as
the models of first and second degree polynomials) now proceeds by
selecting the model such that the stochastic complexity of the given
data $D$ is smallest.  Although we used a one-part code to encode
data, refined MDL model selection still involves a trade-off between
two terms: a
goodness-of-fit term $L(D \mid \hat{H})$ and a complexity term
$\complexity({\cal H})$. 
However, because
we do not explicitly encode hypotheses $H$ any more, there is no
arbitrariness any more. The resulting procedure can be interpreted
in several different ways, some of which provide us with rationales for MDL beyond
the pure coding interpretation (Sections~\ref{sec:compression}--\ref{sec:prequential}):
\begin{description}
\item[1. Counting/differential geometric interpretation] The parametric
  complexity of a model is the logarithm of the number of {\em
    essentially different}, {\em distinguishable\/} point hypotheses within
  the model.
\item[2. Two-part code interpretation] For large
  samples, the stochastic complexity can be interpreted as a two-part
  codelength of the data after all, where 
  hypotheses $H$ are encoded with a special code that works by 
first discretizing the model space ${\cal H}$ into a set of
  `maximally distinguishable hypotheses', and then assigning equal
  codelength to each of these. 
\item[3. Bayesian interpretation] In many cases, 
refined MDL model selection coincides with Bayes factor model
selection based on a {\em  non-informative prior\/} such as {\em
  Jeffreys' prior\/} \cite{BernardoS94}. 
\item[4. Prequential interpretation] Refined MDL model selection can be
  interpreted as selecting the model with the best predictive
  performance when sequentially predicting {\em unseen\/} test data, in the sense described in Section~\ref{sec:prequential}. This
  makes it an instance of Dawid's \citeyear{Dawid84} {\em prequential\/} model validation
  and also relates it to {\em cross-validation} methods.
\end{description}
Refined MDL allows us to compare models of different functional form. 
It
even accounts
for the phenomenon that 
different models with the same number of parameters may not be equally
`complex':
\begin{example}
\label{ex:modselb} 
\rm
Consider two models from psychophysics
describing the relationship between physical dimensions (e.g., light
intensity) and their psychological counterparts (e.g. brightness)
\cite{MyungBP00}: $y = ax^b + \noiseterm$ (Stevens' model) and $y = a
\ln(x+b) + \noiseterm$ (Fechner's model) where $\noiseterm$ is a
normally distributed noise term.  Both models have two free
parameters; nevertheless, it turns out that in a sense, Stevens' model
is more {\em flexible\/} or {\em complex\/} than Fechner's. Roughly speaking,
this means there are a lot more data patterns that can be
{\em explained\/} by Stevens' model than can be explained by Fechner's
model. \citeN{MyungBP00} generated many samples of size 4 from
Fechner's model, using some fixed parameter values. They then 
fitted both models to each sample. In 67\% of the trials,
Stevens' model fitted the data better than Fechner's, even though the
latter generated the data. Indeed, in refined MDL, the `complexity' associated
with Stevens' model is much larger than the complexity associated with
Fechner's, and if both models fit the data equally well, MDL will
prefer Fechner's model.
\end{example}
\index{overfitting} 
Summarizing, refined MDL removes the arbitrary aspect of crude,
two-part code MDL and associates parametric models with an inherent
`complexity' that does not depend on any particular description method
for hypotheses. We should, however, warn the reader that we only 
discussed
a special, simple situation in which we compared a finite
number of parametric models that satisfy certain regularity
conditions. Whenever the models do not satisfy these conditions, or if
we compare an infinite number of models, then the refined  ideas have to be
extended.  We then obtain a `general' refined MDL Principle, which
employs a combination of one-part and two-part codes. 
\commentout{Summarizing: \ownbox{Crude and Refined MDL}{ In
    practice, MDL is always applied to {\em probabilistic models\/}
    ${\cal H}$, such that each $H \in {\cal H}$ defines a probability distribution
    on samples.

In crude MDL, we use the Shannon-Fano code with $L(D|H) = - \log
P(D|H)$ as `the codelength of the data $D$ when encoded with the help
of point hypothesis $H \in {\cal H}$'. In refined MDL, we additionally use codes
with lengths $\luniv(D|{\cal H})$. This codelength is called the {\em
  stochastic complexity of $D$ given ${\cal H}$} and is to be interpreted as
`the codelength of data $D$ when encoded with the help of model ${\cal
  H}$.
Broadly speaking, $\luniv(D|H) = L(D\mid \hat{H}) + C$ where $\hat{H}$
is the distribution within ${\cal H}$ that best fits $D$, and $C$ is the
{\em parametric complexity\/} of ${\cal H}$, indicating how well
${\cal H}$ can
fit random noise.
}
}
\section{The MDL Philosophy}
\label{sec:philosophy}
The first central MDL idea is that every regularity in data may be
used to compress that data; the second central idea is that learning
can be equated with finding regularities in data.  Whereas the first
part is relatively straightforward, the second part of the idea implies that {\em methods for learning from data must have a clear
  interpretation independent of whether any of the models under
  consideration is `true' or not}. Quoting 
J. Rissanen \citeyear{Rissanen89}, the main
originator of MDL:
\begin{ownquote}
``We never want to make the false assumption that the
  observed data actually were generated by a
  distribution of some kind, say Gaussian, and then go on to analyze
  the consequences and make further deductions.  Our deductions may be
  entertaining but quite irrelevant to the task at hand, namely, to
  learn useful properties from the data.''\begin{flushright}{\it Jorma
    Rissanen [1989]}\end{flushright}
\end{ownquote}
Based on such ideas, Rissanen has developed a radical philosophy of
learning and statistical inference that is considerably different from
the ideas underlying mainstream statistics, both frequentist and
Bayesian. We now describe this philosophy in more detail:
\paragraph{1. Regularity as Compression}
According to Rissanen, the goal of inductive inference should be to
`squeeze out as much regularity as possible' from the given data. The
main task for statistical inference is to distill the meaningful
information present in the data, i.e. to separate structure
(interpreted as the regularity, the `meaningful information') from
noise (interpreted as the `accidental information').  For the three
sequences of Example~\ref{ex:bernoulli}, this would amount to the
following: the first sequence would be considered as entirely regular
and `noiseless'. The second sequence would be considered as entirely
random - all information in the sequence is accidental, there is no
structure present. In the third sequence, the structural part would
(roughly) be
the pattern that 4 times as many 0s than 1s occur; given this
regularity, the description of exactly which of all sequences with
four times as many 0s than 1s occurs, is the accidental information.

\paragraph{2. Models as Languages} Rissanen interprets models (sets of hypotheses) as nothing more than languages for
describing useful properties of the data -- a model ${\cal H}$ is {\em
  identified\/} with its corresponding universal code $\luniv(\cdot
\mid {\cal H})$. Different individual hypotheses within the models express
different regularities in the data, and may simply be regarded as {\em
  statistics}, that is, summaries of certain regularities in the data.
{\em These regularities are present and meaningful independently of whether
some $H^* \in {\cal H}$ is the `true state of nature' or not}. Suppose that
the model
${\cal H}$ under consideration is probabilistic. In
traditional theories, one typically assumes that some $P^* \in {\cal H}$
generates the data, and then `noise' is defined as a random quantity relative to this
$P^*$. In the MDL view `noise' is defined relative to the model ${\cal
  H}$
as the residual number of bits needed to encode the data once the
model ${\cal H}$ is given. Thus, noise is {\em not\/} a random
variable: it is a function only 
of the chosen model and
the {\em actually observed data}. Indeed, there is
no place for a `true distribution' or a `true state of nature' in this
view -- there are only models and data. To bring out the difference to
the ordinary statistical viewpoint, consider the phrase `these
experimental data are quite noisy'. According to a traditional
interpretation, such a statement means that the data were generated by
a distribution with high variance. According to the MDL philosophy,
such a phrase means only that the data are not compressible with the
currently hypothesized model -- as a matter of principle, it can {\em
  never\/} be ruled out that there exists
a different model under which the
data are very compressible (not noisy) after all!
\paragraph{3. We Have Only  the Data}
Many (but not all\endnote{For example, cross-validation cannot 
  easily be interpreted in such terms of `a method hunting for the
  true distribution'.}) other methods of inductive inference are based on
the idea that there exists some `true state of nature', typically a
distribution assumed to lie in some model ${\cal H}$. The methods
are then designed as a means to identify or approximate this state of nature based on as
little data as possible. According to Rissanen\endnote{The present
  author's own views are somewhat milder in this respect, but this is
  not the place to discuss them.}, such methods are fundamentally
flawed. The main reason is that the methods are designed
under the assumption  that the true state of nature is
in the assumed model ${\cal H}$, which is often not the case. Therefore,
{\em such methods only admit a clear
interpretation under assumptions that are typically violated in
practice}. Many cherished statistical methods are designed in this way - we
  mention hypothesis testing, minimum-variance unbiased estimation,
  several non-parametric methods, and even some forms of Bayesian inference -- see Example~\ref{ex:bayescraze}.
In contrast, MDL has a clear
interpretation which {\em depends only on the data}, and not on the
assumption of any underlying `state of nature'. 
\begin{ownquote}
\begin{example}{\bf [Models that are Wrong, yet Useful]}
\label{ex:sanity}
\rm Even though the models under consideration are often wrong, they
can nevertheless be very {\em useful}. Examples are the successful
`Naive Bayes' model for spam filtering, Hidden Markov Models for
speech recognition (is speech a stationary ergodic process? probably
not), and the use of linear models in econometrics and psychology.
Since these models are evidently wrong, it seems strange to base
inferences on them using methods that are designed under the
assumption that they contain the true distribution.  To be fair, we
should add that domains such as spam filtering and speech recognition
are not what the fathers of modern statistics had in mind when they
designed their procedures -- they were usually thinking about much
simpler domains, where the assumption that some distribution $P^* \in
{\cal H}$ is `true'
may not be so unreasonable.
\end{example}
\end{ownquote}\paragraph{4. MDL and Consistency}
\label{page:consistent}
Let ${\cal H}$ be a probabilistic model, such that each $P \in {\cal H}$ is a
probability distribution. Roughly, a statistical procedure is called
{\em consistent\/} relative to ${\cal H}$ if, for all
$P^* \in {\cal H}$, the following holds: suppose data are distributed according to $P^*$. Then
given enough data, the learning method will learn a good approximation
of $P^*$ with high probability. Many traditional statistical methods have been designed with
consistency in mind (Section~\ref{sec:statistics}).

The fact that in MDL, we do not assume a true distribution may suggest
that we do not care about statistical consistency.  But this is
not the case: we would still like our statistical method to be such
that in the {\em idealized\/} case where one of the distributions in
one of the models under consideration actually generates the data, our
method is able to identify this distribution, given enough data. If
even in the idealized special case where a `truth' exists within our
models, the method fails to learn it, then we certainly cannot trust
it to do something reasonable in the more general case, where there
may not be a `true distribution' underlying the data at all. So:
consistency {\em is\/} important in the MDL philosophy, but it is used
{\em as a sanity check (for a method that has been developed without
  making distributional assumptions) rather than as a design
  principle}.

In fact, mere consistency is not sufficient. We would like our method
to converge to the imagined true $P^*$ {\em fast}, based on as small a
sample as possible. Two-part code MDL with `clever' codes achieves
good rates of convergence in this sense (\citeN{BarronC91},
complemented by \cite{Zhang04}, show that in many situations,
the rates are {\em  minimax optimal\/}). The same seems to be
true for refined one-part code MDL \cite{BarronRY98}, 
although there is at least one surprising
exception where inference based on the NML and Bayesian universal
model behaves abnormally -- see \cite{CsiszarS00} for the details.

\ \\ \noindent
Summarizing this section, the MDL
philosophy is quite agnostic about whether any of the models under
consideration is `true', or whether something like a `true
distribution' even exists. Nevertheless, 
it has  been suggested  \cite{Webb96,Domingos99} that MDL embodies a
naive belief that `simple models' are `a priori more likely to be
true' than complex models. Below we explain why such claims are mistaken.
\section[MDL and Occam's Razor]{MDL and Occam's Razor} 
\label{sec:occam} 
\index{Occam's Razor} \index{Ockham, W. of} \index{model!true}
\index{arbitrariness} \index{simplicity} When two models fit the data
equally well, MDL will choose the one that is the `simplest' in the
sense that it allows for a shorter description of the data. As such,
it implements a precise form of Occam's Razor -- {\em even though as
  more and more data becomes available, the model selected by MDL may
  become more and more `complex'!\/} Occam's Razor is sometimes
criticized for being either (1) arbitrary or (2) false
\cite{Webb96,Domingos99}. Do these criticisms apply to MDL as well? 
\paragraph{`1. Occam's Razor (and MDL) is arbitrary'}
\label{page:arbitrary}
Because `description length' is a syntactic notion it may seem that
MDL selects an arbitrary model: different codes would have led to
different description lengths, and therefore, to different models.  By
changing the encoding method, we can make `complex' things `simple'
and vice versa.  This overlooks the fact we are not allowed to use
just any code we like! `Refined' MDL tells us to use a specific code, independent of any specific
parameterization of the model, leading to a notion of complexity that
can also be interpreted without any reference to `description
lengths' (see also 
Section~\ref{sec:perceived}).

\paragraph{`2. Occam's Razor is false'} 
\label{page:false}
It is often claimed that Occam's razor is false - we often try to
model real-world situations that are arbitrarily complex, so why
should we favor simple models? In the words of G. Webb\endnote{Quoted
  with permission from KDD Nuggets 96:28, 1996.}:  `What good are
simple models of a complex world?'

The short answer is: { even if the true data generating machinery
  is very complex, it may be a good strategy to prefer simple models
  for small sample sizes}. Thus, MDL (and the corresponding form of
Occam's razor) is a {\em strategy\/} for inferring models from data
  (``choose simple models at small sample sizes''),
not a statement about how the world works (``simple models are
more likely to be true'') -- indeed, a strategy cannot be true or
false, it is `clever' or `stupid'. And the strategy of preferring
simpler models is clever even if the data generating process is highly
complex, as illustrated by the following example:
\begin{example}
\label{ex:infinitelycomplex}{\bf [`Infinitely' Complex Sources]}
\rm Suppose that data are subject to the law
$Y = g(X) + Z$ where $g$ is some continuous function and $Z$ is some
noise term with mean $0$. If $g$ is not a polynomial, but $X$ only
takes values in a finite interval, say $[-1,1]$, we may still
approximate $g$ arbitrarily well by taking higher and higher degree
polynomials. For example, let $g(x) = \exp(x)$. Then, if we use MDL to
learn a polynomial for data $D = ((x_1,y_1), \ldots, (x_n,y_n))$, the
degree of the polynomial $\twop{f}^{(n)}$ selected by MDL at sample
size $n$ will increase with $n$, and with high probability,
$\twop{f}^{(n)}$ converges to $g(x) = \exp(x)$ in the sense that
$\max_{x \in [-1,1]} |\twop{f}^{(n)}(x) - g(x) | \rightarrow 0$.
%\commentout{
Of course, if we had better prior knowledge about the problem we could
have tried to learn $g$ using a model class ${\cal M}$ containing the
function $y = \exp(x)$. But in general, both our imagination and our
computational resources are limited, and we may be forced to use
imperfect models.
%}
\end{example}
If, based on a small sample, we choose the best-fitting
polynomial $\hat{f}$ within the set of {\em all\/} polynomials,
then, even though $\hat{f}$ will fit the data very well, it is likely
to be quite unrelated to the `true' $g$, and $\hat{f}$ may lead to
disastrous predictions of future data. The reason is that, for small samples, the set of
all polynomials is very large compared to the set of possible data
patterns that we might have observed. Therefore, any particular data
pattern can only give us very limited information about which
high-degree polynomial best approximates $g$. 
On the other hand, if we choose the
best-fitting $\hat{f}^\circ$ in some much smaller set such as the set
of second-degree polynomials, then it is highly probable that the
prediction quality (mean squared error) of $\hat{f}^\circ$ on future data
is about the same as its mean squared error on the data we observed:
the size (complexity) of the contemplated model is relatively small
compared to the set of possible data patterns that we might have
observed. Therefore, the particular pattern that we do observe gives
us a lot of information on what second-degree polynomial best
approximates $g$.

Thus, (a) $\hat{f}^\circ$ 
typically leads to better predictions of future
data than $\hat{f}$; and (b) unlike $\hat{f}$, $\hat{f}^\circ$ 
is {\em
  reliable\/} in that it gives a correct impression of how good it
will predict future data {\em even if the `true' $g$ is `infinitely' complex}.
This idea does not just appear in MDL, but is also the basis of
Vapnik's \citeyear{Vapnik98} Structural Risk Minimization approach and
many standard statistical methods for non-parametric inference. In
such approaches one acknowledges that the data generating machinery
can be infinitely complex (e.g., not describable by a finite degree
polynomial).  Nevertheless, it is still a good strategy to approximate
it by simple hypotheses (low-degree polynomials) as long as the sample
size is small. Summarizing: \ownbox{The Inherent Difference between
  Under- and Overfitting}{ If we choose an overly simple model for our
  data, then the best-fitting point hypothesis within the model is likely to
  be almost the best predictor, within the simple model, of future
  data coming from the same source.  If we overfit (choose a very
  complex model) and there is noise in our data, then, {\em even if
    the complex model contains the `true' point hypothesis}, the
  best-fitting point  hypothesis within the model is likely to lead to very
  bad predictions of future data coming from the same source.}  This
statement is very imprecise and is  meant more to convey the general idea
than to be completely true.  As will become clear in
Section~\ref{sec:perceived},
it becomes provably true if we use MDL's measure of model complexity; we
measure prediction quality by logarithmic loss; and we assume that one
of the distributions in ${\cal H}$ actually generates the data.
\section{History}
\label{sec:history}
The MDL Principle has mainly been developed by J. Rissanen in a series
of papers starting with \cite{Rissanen78}. It has its roots in the
theory of {\em Kolmogorov\/} or {\em algorithmic\/} complexity
\cite{LiV97}, developed in the 1960s by \citeN{Solomonoff64},
\citeN{Kolmogorov65} and Chaitin \citeyear{Chaitin66,Chaitin69}.
Among these authors, Solomonoff (a former student of the famous
philosopher of science, Rudolf Carnap) was explicitly interested in
inductive inference. The 1964 paper contains explicit suggestions on
how the underlying ideas could be made practical, thereby
foreshadowing some of the later work on two-part MDL. While Rissanen was not
aware of Solomonoff's work at the time, Kolmogorov's [1965] paper
did serve as an inspiration for Rissanen's \citeyear{Rissanen78} 
development of MDL. 

Another important inspiration for Rissanen was Akaike's \citeyear{Akaike73} AIC method for
model selection, essentially the first model
selection method based on information-theoretic ideas. Even though
Rissanen was inspired by AIC, both the actual method and the
underlying philosophy are substantially different from MDL.

MDL is much closer related to the {\em Minimum Message Length
  Principle}, developed by Wallace and his co-workers in a series of
papers starting with the ground-breaking \cite{WallaceB68}; other
milestones are \cite{WallaceB75} and \cite{WallaceF87}. Remarkably,
Wallace developed his ideas without being aware of the notion of Kolmogorov complexity.
Although Rissanen became aware of Wallace's work before the
publication of \cite{Rissanen78}, he developed his ideas mostly
independently, being influenced rather by Akaike and
Kolmogorov. Indeed, despite the
close resemblance of both methods in practice, the underlying
philosophy is quite different - see Section~\ref{sec:others}.

The first publications on MDL only mention two-part codes. Important
progress  was made by  \citeN{Rissanen84}, in which prequential codes
are employed for the first time and \cite{Rissanen87}, introducing the
Bayesian mixture codes into MDL. This led to
the development of the notion of stochastic complexity as the shortest
codelength of the data given a model \cite{Rissanen86a,Rissanen87}. However, the
connection to Shtarkov's {\em normalized maximum likelihood code \/}
was not made until 1996, and this prevented the full development of the
notion of `parametric complexity'. In the mean time, in his impressive 
Ph.D. thesis, \citeN{Barron85} showed how a specific version of the  
two-part code criterion has excellent frequentist statistical
consistency properties. This was extended by
\citeN{BarronC91} who achieved a breakthrough for two-part codes: they gave clear
prescriptions on how to design codes for hypotheses,  relating codes
with good minimax codelength properties to rates of convergence in
statistical consistency theorems. Some of the ideas of  \citeN{Rissanen87} and
\citeN{BarronC91} were, as it were, unified when \citeN{Rissanen96}
introduced a new definition of stochastic complexity based on the {\em
  normalized maximum likelihood code}
(Section~\ref{sec:universal}). The resulting theory was summarized for
the first time by \citeN{BarronRY98}, and is called `refined MDL' in
the present overview.
\section{Summary and Outlook}
\label{sec:summary}
We discussed how regularity is related to data compression, and how
MDL employs this connection by viewing learning in terms of data
compression. One can make this precise in several ways; in {\em
  idealized\/} MDL one looks for the shortest program that generates
the given data. This approach is not feasible in practice, and here we
concern ourselves with {\em practical\/} MDL. Practical
MDL comes in a crude version based on two-part codes and in a modern,
more refined version based on the concept of {\em universal coding}.
The basic ideas underlying all these approaches can be found in the
boxes spread throughout the text.

These methods are mostly applied to model selection but can also be
used for other problems of inductive inference. In contrast to most
existing statistical methodology, they can be given a clear
interpretation irrespective of whether or not there exists some `true'
distribution generating data -- inductive inference is seen as a
search for regular properties in (interesting statistics of) the data,
and there is no need to assume anything outside the model and the
data. In contrast to what is sometimes thought, there is {\em no\/}
implicit belief that `simpler models are more likely to be true' --
MDL does embody a preference for `simple' models, but this is best
seen as a strategy for inference that can be useful even if the environment is not
simple at all.

In the next chapter, we make precise both the crude and the refined
versions of practical MDL. For this, it is absolutely essential that
the reader familiarizes him- or herself with two basic notions of
coding and information theory: the relation between codelength functions and probability
distributions, and (for refined MDL), the idea of universal coding --
a large part of the chapter will be devoted to these.

   \begingroup
     \parindent 0pt
     \theendnotes
     \endgroup

%%% Local Variables: 
%%% mode: latex
%%% TeX-master: t
%%% End: 

\ \\ \ \newpage
\chapter[Tutorial on MDL]{Minimum Description Length Tutorial}
  \label{chap:grunwald}
%  \authorentry{Peter Gr\"unwald}
\setcounter{endnote}{0}

\newcommand{\smfrac}[2]{\mbox{$\frac{#1}{#2}$}}

% use next line ONLY if you have three or more authors; otherwise remove it.
%\clearpage

% now follows the abstract (to be replaced by your own abstract).

%you can replace the following sections with the body of your .tex file
\section{Plan of the Tutorial}
In Chapter~\ref{chap:survey} we introduced the MDL
Principle in an informal way. In this chapter we give an introduction
to MDL that is mathematically precise. Throughout the text,
we assume some basic familiarity with probability theory. While some
prior exposure to basic statistics is highly useful, it is not
required. The chapter can be read without any prior knowledge of
information theory.  The tutorial is organized according to the
following plan:
\begin{itemize}
\item The first two sections are of a preliminary nature: 
\begin{itemize}
\item Any understanding of MDL requires some minimal knowledge of
  information theory -- in particular the relationship between
  probability distributions and codes. This relationship is explained
  in Section~\ref{sec:prob}.
\item Relevant statistical
  notions such as `maximum likelihood estimation'  are
  reviewed in Section~\ref{sec:statistics}. There we also introduce the
  Markov chain model which will serve as an example model throughout the
  text.
\end{itemize}
\item Based on this preliminary material, in Section~\ref{sec:crude}
  we formalize a simple version of the MDL Principle,
  called the {\em crude two-part MDL Principle\/} in this text. We
  explain why, for successful practical applications, crude MDL
  needs to be refined.
\item Section~\ref{sec:universal} is once again preliminary: it
  discusses {\em universal coding}, the information-theoretic concept 
underlying refined
  versions of MDL.
\item Sections~\ref{sec:refined}--\ref{sec:beyond} 
define and discuss refined MDL. They 
are the key sections of the tutorial:
\begin{itemize}
\item Section~\ref{sec:refined} discusses basic refined MDL for
  comparing a finite number of simple statistical models and introduces
  the central concepts of {\em parametric\/} and {\em stochastic
    complexity}. It gives an {\em asymptotic expansion\/} of these
  quantities and interprets them from a compression, a geometric, a
  Bayesian and a predictive point of view.
\item Section~\ref{sec:unify} extends refined MDL to harder model
  selection problems, and in doing so reveals the general, unifying
  idea, which is summarized in Figure~\ref{fig:unified}.
\item Section~\ref{sec:beyond} briefly discusses how to extend MDL to
  applications beyond model section.
\end{itemize}
\item Having defined `refined MDL' in 
  Sections~\ref{sec:refined}--\ref{sec:beyond}, the next two sections
  place it in context:
\begin{itemize}
\item Section~\ref{sec:others} compares MDL to other approaches to
  inductive inference, most notably the related but different {\em
    Bayesian\/} approach.  
\item Section~\ref{sec:problems} discusses perceived as well as real
  problems with MDL. The perceived problems relate to MDL's relation
  to Occam's Razor, the real problems relate to the fact that
  applications of MDL sometimes perform suboptimally in practice.
\end{itemize}
\item Finally, Section~\ref{sec:conclusion} provides a conclusion.
\end{itemize}
\ownbox{Reader's Guide}{Throughout the text, paragraph headings
  reflect the most important concepts. Boxes summarize 
the most important findings. Together, paragraph headings and boxes provide an overview of
MDL theory. 

It is possible to read this chapter without having read
the non-technical overview of Chapter~\ref{chap:survey}. However, we
strongly recommend reading at least Sections~\ref{sec:modselintro} and
Section~\ref{sec:cruderefined} before embarking on the present chapter.}
\section{Information Theory I: Probabilities and Codelengths}
\label{sec:prob}
%To make our ideas precise, we need to know some basic information
%(coding theory). 
This first section is a mini-primer on information theory, focusing
on the relationship between probability distributions and codes. A
good understanding of this relationship is essential for a good
understanding of MDL.  After some preliminaries,  
Section~\ref{sec:prefix} introduces prefix codes, the
type of codes we work with in MDL. These are related to probability
distributions in two ways. In Section~\ref{sec:kraft} we discuss the
first relationship, which is related to the {\em Kraft inequality\/}:
for every probability mass function $P$, there exists a code with
lengths $- \log P$, and vice versa.  Section~\ref{sec:information}
discusses the second relationship, related to the {\em information
  inequality}, which says that if the data are distributed according
to $P$, then the code with lengths $- \log P$ achieves the minimum
expected codelength. Throughout the section we give examples relating
our findings to our discussion of regularity and compression in
Section~\ref{sec:basic} of Chapter~\ref{chap:survey}.
\paragraph{Preliminaries and Notational Conventions - Codes}
We use $\log$ to denote
logarithm to base 2. For real-valued $x$ we use $\lceil x \rceil $ to denote the {\em ceiling\/} of $x$, that is, $x$ rounded up to the nearest integer. 
We often abbreviate $x_1, \ldots, x_n$ to $x^n$.
Let $\xspace$ be a finite or countable set. A {\em code\/} for
$\xspace$ is defined as a 1-to-1 mapping from $\xspace$ to $\cup_{n
  \geq 1} \{0,1\}^n$.  $\cup_{n
  \geq 1} \{0,1\}^n$ is the set of binary strings (sequences of $0$s and
$1$s) of length $1$ or larger. For a given code $C$, we use $C(x)$ to
denote the encoding of $x$. Every code $C$ induces a function $L_C: 
\xspace \rightarrow \naturals$ called the {\em codelength function}.
\index{codelength function} 
$L_C(x)$ is the number of bits (symbols)
needed to encode $x$ using code $C$. 

Our definition of code implies that we only consider {\em lossless\/}
encoding in MDL\endnote{But see Section~\ref{sec:kolmogorov}.}: for any description $z$ it is always possible to retrieve
the unique $x$ that gave rise to it. More precisely, because the code
$C$ must be 1-to-1, 
there is at most one $x$ with $C(x) = z$. Then
$x = C^{-1}(z)$, where the inverse $C^{-1}$ of $C$ is sometimes called
a `decoding function' or `description method'.
\index{code}
\index{description method}
\paragraph{Preliminaries and Notational Conventions - Probability}
Let $P$ be a probability distribution defined on a finite or countable set
$\xspace$. 
We use $P(x)$ to denote the probability of $x$, and we 
denote the corresponding random variable by
$X$. If $P$ is a function on finite or
countable $\xspace$ such that $\sum_x P(x) < 1$, we call $P$ a {\em
  defective\/} distribution. A defective distribution may  be
thought of as a probability distribution that puts some of its mass on
an imagined outcome that in reality will never appear.

A {\em probabilistic source\/} $P$ is a sequence of probability
\index{probabilistic source} distributions $P^{(1)}, P^{(2)}, \ldots$
on $\xspace^1, \xspace^2, \ldots$ such that for all $n$,
$P^{(n)}$ and $P^{(n+1)}$ are {\em compatible\/}: $P^{(n)}$ is equal to 
the `marginal' distribution of $P^{(n+1)}$ restricted to $n$ outcomes. 
That is, for all $x^n \in {\cal X}^n$,
$P^{(n)}(x^n) =
\sum_{y \in \xspace} P^{(n+1)}(x^n, y)$. Whenever this cannot cause any confusion, we write $P(x^n)$ rather than $P^{(n)}(x^n)$.
A probabilistic source may be thought of as a probability
distribution on infinite sequences\endnote{Working directly with
  distributions on infinite sequences is more elegant, but it requires
  measure theory, which we want to avoid here.}. We say that the data
are {\em i.i.d.\/} (independently and identically distributed)
\index{i.i.d.}  under source $P$ if for each $n$, $x^n \in \xspace^n$,
$P(x^n) = \prod_{i=1}^n P(x_i)$.
\subsection{Prefix Codes}
\label{sec:prefix}
In MDL we only work with a subset of all possible codes, the so-called
{\em prefix codes}. A prefix code\endnote{Also known as {\em
    instantaneous codes\/} and called, perhaps more
  justifiably, `prefix-free' codes in \cite{LiV97}.} is a code such
that no codeword is a prefix of any other codeword. For example, let
$\xspace = \{a, b, c\}$. Then the code $C_1$ defined by $C_1(a) = 0$,
$C_1(b) = 10$, $C_1(c) = 11$ is prefix. The code $C_2$ with $C_2(a) =
0, C_2(b) = 10$ and $C_2(c) = 01$, while allowing for lossless decoding, is {\em
  not\/} a prefix code since $0$ is a prefix of $01$.
The prefix requirement is natural, and
nearly ubiquitous in the data compression literature. We now explain why.
\begin{example}
  \rm Suppose we plan to encode a sequence of symbols $(x_1, \ldots,
  x_n) \in \xspace^n$.  We already designed a
  code $C$ for the elements in $\xspace$.  The
  natural thing to do is to encode $(x_1,\ldots, x_n)$ by the concatenated string
  $C(x_1)C(x_2)\ldots C(x_n)$. In order for this method to succeed for
  all $n$, all $(x_1, \ldots, x_n) \in \xspace^n$, the
  resulting procedure must define a code, i.e. the function $C^{(n)}$ mapping
  $(x_1,\ldots, x_n)$ to $C(x_1)C(x_2)\ldots C(x_n)$ must be
  invertible.  If it were not, we would
  have to use some marker such as a comma to separate the codewords.
  We would then really be using a ternary rather than a binary
  alphabet. 

Since we always want to construct codes for sequences
  rather than single symbols, we only allow codes $C$ such that
the extension $C^{(n)}$ defines a code for all $n$. We say that such
codes 
have `uniquely decodable extensions'. 
It is easy to see that (a) every prefix code has uniquely decodable
extensions. Conversely, although this is not at all easy to see, it turns
out that  (b), for every code $C$ with uniquely
decodable extensions, there exists a prefix code $C_0$ such that for
all $n, x^n \in \xspace^n$, $L_{C^{(n)}}(x^n) = L_{C^{(n)}_0}(x^n)$ 
\cite{CoverT91}. Since in MDL we
are only interested in code-{\em lengths}, and never in actual codes,
we can restrict ourselves to prefix codes without loss of generality.

Thus, the restriction to prefix code may also be understood as a means
to send concatenated messages while avoiding the need to introduce
extra symbols into the alphabet.
\end{example}
Whenever
in the sequel we speak of `code', we really mean `prefix code'. We
call a prefix code $C$ for a set $\xspace$ {\em complete\/} if there exists
no other prefix code that compresses at least one $x$ more and no $x$
less then $C$, i.e. if there exists no code $C'$ such that for all
$x$, $L_{C'}(x) \leq L_C(x)$ with strict inequality for at least one
$x$.
\subsection[The Kraft Inequality - Codelengths \& Probabilities, I]{The Kraft Inequality - Codelengths and Probabilities, Part I}
\label{sec:kraft}
In this subsection we relate prefix codes to probability distributions. 
Essential for understanding the relation is the fact that no matter what code we use, {\em most sequences cannot be compressed}, as demonstrated by the following example:
\begin{example}{\bf [Compression and Small Subsets:
    Example~\ref{ex:bernoulli}, Continued.]}
\rm
\label{page:compression}
\label{ex:nocompression}
In Example~\ref{ex:bernoulli} we featured the following
three sequences:
\begin{eqnarray}
\label{seqla}
& {\tt
00010001000100010001} \ \ \ldots \ \ {\tt 0001000100010001000100010001\ } & \\
\label{seqlb}
& {\tt
01110100110100100110} \ \ \ldots \ \ {\tt 1010111010111011000101100010\ } & \\
\label{seqlc}
& {\tt
00011000001010100000} \ \ \ldots \ \ {\tt 0010001000010000001000110000\ }&
\end{eqnarray}
We showed 
that (a) the first sequence - an $n$-fold repetition of $0001$
- could be substantially compressed if we use as our code a
general-purpose programming language (assuming that valid 
programs must end with
a {\tt halt}-statement or a closing bracket, such codes satisfy the prefix
property). We also claimed that (b) the second
sequence, $n$ independent outcomes of fair coin tosses, cannot be
compressed, and that (c) the third sequence could be compressed to $\alpha
n$ bits, with $0 < \alpha < 1$. We are now in a position to prove
statement (b): strings which are `intuitively'
random cannot be substantially compressed.  
%By `it is virtually impossible' we mean `it happens with
%vanishing probability'. 
Let us take some arbitrary but fixed description method 
over the data alphabet consisting of the set of all binary sequences
of length $n$. Such a code maps binary strings to binary strings.
There are $2^n$
possible data sequences of length $n$. Only two of these
can be mapped to a description of length $1$ (since there are only two
binary strings of length 1: `0' and `1'). Similarly, only a subset of
at most $2^m$ sequences can have a description of length $m$. This
means that at most $\sum_{i=1}^m 2^m < 2^{m+1}$ data sequences can
have a description length $\leq m$. The fraction of data sequences of
length $n$ that can be compressed by more than $k$ bits is therefore
at most $2^{-k}$ and as such decreases exponentially in $k$. If data
are generated by $n$ tosses of a fair coin, then all $2^n$
possibilities for the data are equally probable, so the probability
that we can compress the data by more than $k$ bits is smaller than
$2^{-k}$. For example, the probability that we can compress the data
by more than 20 bits is smaller than one in a million.

% I checked the calculation above; it's OK.
We note that {\em after\/} the data (\ref{seqlb}) has been observed,
it is always possible to design a code which uses arbitrarily few bits
to encode this data - 
the actually observed sequence may be encoded as `1' for
example, and no other sequence is assigned a codeword. The point is
that with a code that has been designed {\em before\/} seeing any
data, it is virtually impossible to substantially compress randomly
generated data.  
\end{example}
The example demonstrates that achieving a short description length for
the data is equivalent to identifying the data as belonging to a tiny,
very {\em special\/} subset out of all a priori possible data
sequences.

\paragraph{A Most Important Observation} 
Let $\zspace$ be finite or countable. 
For concreteness, we may take
$\zspace = \{0,1\}^{n}$ for some large $n$, say $n = 10000$. 
From Example~\ref{ex:nocompression} we know that, no matter what
code we use to encode values in $\zspace$, `most' outcomes in
$\zspace$ will not be substantially compressible: at most two outcomes
can have description length $1 = - \log 1/2$; at most four outcomes
can have length $2 = - \log 1/4$, and so on. Now consider any
probability distribution on $\zspace$.  Since the probabilities $P(z)$
must sum up to one $(\sum_z P(z) = 1)$, `most' outcomes in $\zspace$
must have small probability in the following sense: at most 2 outcomes
can have probability $\geq 1/2$; at most 4 outcomes can have
probability $\geq 1/4$; at most $8$ can have $\geq 1/8$-th etc. This
suggests an analogy between codes and probability distributions: each
code induces a code length function that assigns a number to each $z$,
where most $z$'s are assigned large numbers. Similarly, each
distribution assigns a number to each $z$, where most $z$'s are
assigned small numbers.

It turns out that this correspondence can be made mathematically
precise by means of the {\em Kraft inequality\/} \cite{CoverT91}. We
neither precisely state nor prove this inequality; rather, in
Figure~\ref{fig:most} we state an
immediate and fundamental consequence:
{\em probability mass functions correspond to codelength
  functions}. The following example illustrates this and at the same
time introduces a type of code that will be frequently employed in the
sequel:
\begin{figure}\ownbox{
%The most important  observation of this tutorial: \\ 
Probability Mass Functions
  correspond to Codelength Functions}{
\label{box:most}
Let $\zspace$ be a finite or countable set and let $\prob$ be a probability
distribution on $\zspace$. 
%Let $\prob$ be a
%probability distribution on $\zspace^n$, the set of sequences of
%length $n$. 
Then there exists a prefix code $C$ for $\zspace$
such that for all $z \in \zspace$, $L_C(z) = \lceil - \log
\pd(z)\rceil $. $C$ is called the {\em
  code corresponding to $P$}.

Similarly, let $C'$ be a
prefix code for $\zspace$.
Then there exists a (possibly  defective) probability
distribution $\prob'$ such that for all $z \in \zspace$, 
$- \log \pd'(z) =
L_{C'}(z)$. $\prob'$ is called the {\em probability distribution
corresponding to $C'$}.
\begin{ownquote}
Moreover $C'$ is a {\em complete prefix code\/} iff $P$ is  proper
($\sum_{z} P(z) = 1$).
\end{ownquote}
Thus, large probability according to $P$ means 
small code length according to the code corresponding to $P$ and vice versa.

We are typically concerned with cases where $\zspace$ represents
sequences of $n$ outcomes; {\em that is, $\zspace = \xspace^n$ $(n \geq
1)$ where $\xspace$ is the sample space for one observation}.
}\vspace*{-0.5  cm}
\caption{\label{fig:most} The most important observation of this tutorial.}
\end{figure}
\label{page:probcode}
\begin{example}{\bf [Uniform Distribution Corresponds to Fixed-length Code]} 
\label{ex:uniform}
\rm
  Suppose $\zspace$ has $M$ elements. The uniform distribution $P_U$
  assigns probabilities $1/M$ to each element. We can arrive at a
  code corresponding to $P_U$ as follows. First, order and number the
  elements in $\zspace$ as $0, 1, \ldots, M-1$. Then, for each $z$
  with number $j$, set $C(z)$ to be equal to $j$ represented as a binary
  number with $\lceil \log M \rceil$ bits.  The resulting code has, for all $z \in \zspace$,
$L_C(z) = \lceil \log M \rceil = \lceil - \log P_U(z) \rceil.$
This is a code corresponding to $P_U$ (Figure~\ref{fig:most}).
In general, there exist several codes corresponding to $P_U$, one for
each ordering of $\zspace$. But all these codes share the same length
function $L_U(z) := \lceil - \log P_U(z) \rceil.$; therefore, $L_U(z)$ is
the unique codelength function corresponding to $P_U$. 

For example, if $M = 4$,
  $\zspace = \{a,b,c,d\}$, we can take $C(a) = 00, C(b) = 01, C(c) = 10,
  C(d) = 11$ and then $L_U(z) =2$ for all $z \in \zspace$. 
In general,
codes corresponding to uniform distributions assign fixed lengths to
each $z$ and are called {\em fixed-length\/} codes. To map a
non-uniform distribution to a corresponding code, we have to use a
more intricate construction \cite{CoverT91}. 
\end{example}
In practical applications, we almost always deal with probability
distributions $P$ and strings $x^n$ such that $P(x^n)$ decreases
exponentially in $n$; for example, this will typically be the case if
data are i.i.d., \index{i.i.d.}  such that $P(x^n) = \prod P(x_i)$.
Then $- \log P(x^n)$ increases linearly in $n$ and the effect of
rounding off $- \log P(x^n)$ becomes negligible.  Note that the code
corresponding to the product distribution of $P$ on $\xspace^n$ does
not have to be the $n$-fold extension of the code for the original
distribution $P$ on $\xspace$ -- if we were to require that, 
the effect of
rounding off would be on the order of $n$ . Instead, we {\em
  directly\/} design a code for the distribution on the larger space
$\zspace = \xspace^n$. In this way, the effect of rounding changes the
codelength by at most 1 bit, which is truly negligible.  For this and
other\endnote{For example, with non-integer codelengths the notion of
  `code' becomes invariant to the size of the alphabet in which we
  describe data.} reasons, we henceforth simply neglect the integer
requirement for codelengths. This simplification allows us to {\em
  identify\/} codelength functions and (defective) probability mass
functions, such that a short codelength corresponds to a high
probability and vice versa. Furthermore, as we will see, in MDL we are
not interested in the details of actual encodings $C(z)$; we only care
about the code lengths $L_C(z)$. It is so useful to think about these
as $\log$-probabilities, and so convenient to allow for non-integer
non-probabilities, that we will simply {\em redefine\/} prefix code
length functions as (defective) probability mass functions that can
have non-integer code lengths -- see Figure~\ref{fig:new}.  The
following example illustrates idealized codelength functions and at
the same time introduces a type of code that will be frequently used
in the sequel:
\begin{figure}
\ownbox{New Definition of Code Length Function}
{\label{box:idealizedcode} In MDL we are {\bf NEVER\/} concerned with
  actual encodings; we are only concerned with code {\em length\/}
  functions. The set of all codelength functions for finite or countable 
sample space
  $\zspace$ is defined as:
\begin{equation}
\label{eq:prefixlengthsb}
\prefixlengths_{\zspace} = \bigl\{ L: \zspace \rightarrow [0,\infty] \mid
\sum_{z \in \samplespace} 2^{-L(z)} \leq 1 \bigr\},
\end{equation}
or equivalently, $\prefixlengths_{\zspace}$ is the set of those
functions $L$ on $\zspace$ such that there exists a function $Q$ with
$\sum_z Q(z) \leq 1$ and for all $z$, $L(z) = - \log Q(z)$. ($Q(z) =
0$ corresponds to $L(z) = \infty$).

Again, $\zspace$ usually represents a
sample of $n$ outcomes: $\zspace = \xspace^n$ $(n \geq
1)$ where $\xspace$ is the sample space for one observation.}\vspace*{-0.5  cm}

\caption{Code lengths are probabilities.\label{fig:new}}
\end{figure}
\begin{example}{\bf `Almost' Uniform Code for the Positive Integers}
\label{ex:integers}  
\rm
Suppose we want to encode a number $k \in \{1,2, \ldots\}$.  In
  Example~\ref{ex:uniform}, we saw that in order to encode a number
  between $1$ and $M$, we need $\log M$ bits. What if we cannot
  determine the maximum $M$ in advance? We cannot just encode $k$
  using the uniform code for $\{1, \ldots, k\}$, since the resulting
  code would not be prefix. So in general, we will need more than $\log k$
  bits. Yet there exists a prefix-free code which performs `almost' as
  well as $\log k$. The simplest of such codes works as follows.  $k$
  is described by a codeword starting with $\lceil \log k \rceil$
  $0$s. This is followed by a $1$, and then $k$ is encoded using the
  uniform code for $\{1, \ldots, 2^{\lceil \log k \rceil} \}$. With
  this protocol, a decoder can first reconstruct $\lceil \log k
  \rceil$ by counting all $0$'s before the leftmost $1$ in the
  encoding. He then has an upper bound on $k$ and can use this
  knowledge to decode $k$ itself. This protocol uses less than $2
  \lceil \log k \rceil + 1$ bits.  Working with idealized, non-integer
  code-lengths we can simplify this to $2 \log k +1$ bits. To see this, consider the function $P(x) = 2^{-2 \log x - 1}$. An
  easy calculation gives
$$
\sum_{x \in 1, 2,\ldots} P(x) =
\sum_{x \in 1, 2,\ldots} 
2^{- 2 \log x -1} = \frac{1}{2} \sum_{x \in 1, 2,\ldots} x^{-2} 
< \frac{1}{2} + \frac{1}{2} \sum_{x=2,3, \ldots} \frac{1}{x(x-1)} = 
%\frac{1}{2} + \frac{1}{2} \sum_{x=2,3, \ldots} \bigl[ \frac{1}{x-1} + 
%\frac{1}{x} \bigr] = 
1,
$$
so that $P$ is a (defective) probability distribution. Thus, by our new
definition (Figure~\ref{fig:new}), there exists a prefix code with,
for all $k$, $L(k) = - \log P(k) = 2 \log k + 1$.  We call the
resulting code the `simple standard code for the integers'. In
Section~\ref{sec:universal} we will see that it is an instance of a
so-called `universal' code.
  
  The idea can be refined to lead to codes with lengths $\log k +
  O(\log \log k)$; the `best' possible refinement, with code lengths
  $L(k)$ increasing monotonically but as slowly as possible in $k$, is known as `the
  universal code for the integers' \cite{Rissanen83}. 
%The gain in
%  using these refinements over our simple code will be so small that
%  for the purposes of this tutorial, 
However, for our purposes in this tutorial, it is good enough to encode
  integers $k$ with $ 2\log k + 1$ bits.
\end{example}
\begin{ownquote}
\begin{example}{\bf [Example~\ref{ex:bernoulli} and~\ref{ex:nocompression}, 
Continued.]}
\rm We are now also in a position to prove the third and final claim of Examples~\ref{ex:bernoulli} and~\ref{ex:nocompression}.
Consider the three sequences (\ref{seqla}), (\ref{seqlb}) and
(\ref{seqlc}) on page~\pageref{seqla} again. 
It remains to investigate how much the third
sequence can be compressed. Assume for concreteness that, before
seeing the sequence, we are told that the sequence
contains a fraction of $1$s equal to $1/5 + \epsilon$ for some small
unknown $\epsilon$. By the Kraft inequality, Figure~\ref{fig:most},
for all distributions $P$, 
there exists some code on sequences of length $n$ such that for all
$x^n \in {\cal X}^n$, $L(x^n) = \lceil - \log P(x^n) \rceil$. The fact
that
the fraction of $1$s is approximately equal to $1/5$ suggests to model
$x^n$ as
independent outcomes of a coin with bias $1/5$-th. The corresponding
distribution $P_0$ satisfies
\begin{multline}
- \log P_0(x^n) = \log \biggl( \frac{1}{5} \biggr)^{n_{[1]}} 
\biggl( \frac{4}{5} \biggr)^{n_{[0]}}
= n \bigl[ - \bigl( \frac{1}{5} + \epsilon \bigr) \log \frac{1}{5}
- 
\bigl( \frac{4}{5} - \epsilon \bigr) \log \frac{4}{5}\bigr] = \\
n[ \log 5 - \frac{8}{5} + 2 \epsilon],
\nonumber
\end{multline}
where $n_{[j]}$ denotes the number of occurrences of symbol $j$ in
$x^n$. For small enough $\epsilon$, the part between brackets is
smaller than $1$, so that,
using the code $L_0$ with lengths $- \log P_0$, the sequence
can be encoded using $\alpha n$ bits were $\alpha$ 
satisfies $0 < \alpha < 1$. Thus, using the code $L_0$,
the sequence can be compressed by a linear amount, if we use a
specially designed code that assigns short codelengths to sequences
with about four times as many $0$s than $1$s.

We note that {\em after\/} the data (\ref{seqlc}) has been observed,
it is always possible to design a code which uses arbitrarily few bits
to encode $x^n$ - the actually observed sequence may be encoded as `1' for
example, and no other sequence is assigned a codeword. The point is
that with a code that has been designed {\em before\/} seeing the
actual sequence, given {\em only\/} the knowledge that the sequence will contain
approximately four times as many 0s than 1s, the sequence is
guaranteed to be compressed by an amount linear in $n$.
\end{example}
\end{ownquote}
\paragraph{Continuous Sample Spaces} How does the correspondence work
for continuous-valued ${\cal X}$? 
In this tutorial we only consider $P$ on $\xspace$ such that $P$ admits
a density\endnote{As understood in elementary probability, i.e.  with respect
to Lebesgue measure.}. Whenever in the following we make a general
statement about sample spaces $\xspace$ and distributions $P$,
$\xspace$ may be finite, countable or any subset of $\reals^l$, for
any integer $l \geq 1$, and $P(x)$ represents the probability mass
function or density of $P$, as the case may be. In the
continuous case, all sums should be read as integrals. The
correspondence between probability distributions and codes may be
extended to distributions on continuous-valued $\xspace$: we may think
of $L(x^n) := - \log P(x^n)$ as a code-length function corresponding
to $\zspace = \xspace^n$ encoding the values in $\xspace^n$ at unit
precision; here $P(x^n)$ is the density of $x^n$ according to $P$. We
refer to  \cite{CoverT91} for further details.
\subsection{The Information Inequality - Codelengths \&
  Probabilities, II}\label{sec:information}
In the previous subsection, we established the first fundamental
relation between probability distributions and codelength functions.
We now discuss the second relation, which is nearly as important.

In the correspondence to codelength functions,
probability distributions were treated as mathematical objects and {\em
nothing else}. That is, if we decide to use a code $C$ to encode our data, this
definitely does {\em not\/} necessarily
mean that we assume our data to be drawn according to
the probability distribution corresponding to $L$: we may have no idea
what distribution generates our data; or conceivably, such a
distribution may not even exist\endnote{Even if one adopts a Bayesian
  stance and postulates that an agent can come up with a 
(subjective) distribution for {\em every\/} conceivable domain, this
problem remains: in practice, the adopted distribution may be so
complicated that we cannot design the optimal code corresponding to
it, and have to use some {\em ad hoc\/}-instead.}.
Nevertheless, {\em if\/} the data are distributed according to some
distribution $P$, {\em then\/} the code corresponding to $P$ turns out
to be the optimal code to use, in an expected sense -- see
Figure~\ref{fig:second}.
\begin{figure}
\ownbox{The $P$ that corresponds to $L$ minimizes expected
  codelength}{
\label{box:second}
Let $P$ be a distribution on (finite, countable or continuous-valued)
$\zspace$ and let $L$ be defined by
\begin{equation}
\label{eq:infineq}
L \isbydefinition \myargmin{L \in \prefixlengths_{\zspace}} E_P[ L(Z)].
\end{equation}
Then $L$ exists, is unique, and is identical to the codelength
function corresponding to $P$, with lengths $L(z) = - \log P(z)$.  
}\vspace*{-0.5  cm}
\caption{\label{fig:second} The second most important observation of
  this tutorial.}
\end{figure}
This result may be recast as follows: for all distributions $P$ and
$Q$ with $Q \neq P$,
$$
%\begin{equationr}
%\label{eq:informationinequality}
E_P [ - \log Q(X)] > E_P [ - \log P(X) ].
$$
%\end{equationr}
In this form, the result is known as the {\em information inequality}.
\index{information inequality|textbf}
It is easily proved using concavity of the logarithm \cite{CoverT91}.

The information inequality says the following: suppose $Z$ is
distributed according to $P$ (`generated by $P$'). Then, among all
possible codes for ${\cal Z}$, the code with lengths $- \log P(Z)$
`on average' gives the shortest encodings of outcomes of $P$. Why
should we be interested in the average? The {\em law of large
  numbers\/} \cite{Feller68a} implies that, for large samples of data
distributed according to $P$, with high
$P$-probability, the code that gives the shortest expected lengths
will also give the shortest {\em actual\/} codelengths, which is what
we are really interested in. This will hold if data are i.i.d., but
also more generally if $P$ defines a `stationary and ergodic' process.
\begin{ownquote}
\begin{example}
\label{ex:lln} \rm
  Let us briefly illustrate this. Let $P^*$, $Q_A$ and $Q_B$ be three
  probability distributions on $\samplespace$, extended to $\zspace = \xspace^n$ by  independence. Hence $P^*(x^n) = \prod P^*(x_i)$
  and similarly for $Q_A$ and $Q_B$. Suppose we obtain a sample
  generated by $P^*$. Mr. A and Mrs. B both want to encode the sample
  using as few bits as possible, but neither knows that $P^*$ has
  actually been used to generate the sample. A decides to use the code
  corresponding to distribution $Q_A$ and B decides to use the code
  corresponding to $Q_B$. Suppose that $E_{P^*} [ - \log Q_A(X)] <
  E_{P^*} [ - \log Q_B(X)]$. Then, by the law of large numbers , with
  $P^*$-probability 1,
$
{n}^{-1} [ - \log Q_j(X_1, \ldots, X_n) ] \rightarrow E_{P^*} [ - \log Q_j(X)],
$ for both $j \in \{A,B\}$
(note $- \log Q_j(X^n) = -\sum_{i=1}^n \log Q_j(X_i)$). It follows
that, with probability 1, Mr. A will need less (linearly in $n$) bits
to encode $X_1, \ldots, X_n$ than Mrs. B.
\end{example}
\end{ownquote}
The qualitative content of this result is not so surprising: in a
large sample generated by $P$, the frequency of each $x \in \samplespace$ will
be approximately equal to the probability $P(x)$. In order to obtain a short
codelength for $x^n$, we should use a code that assigns a small codelength to
those symbols in $\samplespace$ with high frequency (probability), and a large
codelength to those symbols in $\samplespace$ with low frequency (probability).
\paragraph{Summary and Outlook}
In this section we introduced (prefix) codes and thoroughly discussed the relation between probabilities and codelengths. We are now almost ready to formalize a simple version of MDL -- but first we need to review some concepts of statistics.
\section{Statistical Preliminaries and Example Models}
\label{sec:statistics}
In the next section we will make precise the crude form of MDL
informally presented in Section~\ref{sec:modselintro}. We will freely
use some convenient statistical concepts which we review in this
section; for details see, for example, \cite{CasellaB90}. 
We also describe the model class of {\em Markov chains\/} of
arbitrary order, which we use as our running example.  These admit a simpler treatment than the polynomials, to which we return in Section~\ref{sec:regression}.
\paragraph{Statistical Preliminaries}
A {\em probabilistic model\/}\endnote{Henceforth, we simply use `model'
    to denote probabilistic models; we typically use $\cH$ to denote
    sets of hypotheses such as polynomials, and reserve $\M$ for
    probabilistic models.} $\M$ is a set of probabilistic
sources. Typically one uses the word `model' to denote sources of the same functional form.
%Examples of models are the family of all normal
%distributions or the family of all Markov chain distributions of
%arbitrary order.  
We often index the elements $P$ of a
model $\M$ using some parameter $\theta$. In that case we write $P$ as
$P(\cdot \mid \theta)$, and $\M$ as $\M = \{ P(\cdot \mid \theta) \mid
\theta \in \Theta \}$, for some {\em parameter space\/} $\Theta$.  If
$\M$ can be parameterized by some connected  
$\Theta \subseteq \reals^k$ for some $k \geq 1$ and the mapping $\theta
\rightarrow P(\cdot \mid \theta)$ is smooth (appropriately defined),
we call $\M$ a {\em parametric model\/} or {\em family}. For example,
the model $\M$ of all normal distributions on ${\cal X} = \reals$ 
is a parametric model that can be parameterized by
$\theta = (\mu,\sigma^2)$ where $\mu$ is the mean and
$\sigma^2$ is the variance of the distribution indexed by $\theta$. 
The family of all Markov chains of all orders
is a model, but not a parametric model. We call a model $\M$ an {\em
  i.i.d. model\/} \index{i.i.d.} if, according to all $P \in \M$,
$X_1, X_2, \ldots$ are i.i.d. We call $\M$ {\em $k$-dimensional\/} if
$k$ is the smallest integer $k$ so that $\M$ can be smoothly
parameterized by some $\Theta \subseteq 
\reals^k$. 
%In that case, we write $\dim(\M) = k$.  

For a given model $\M$ and sample $D= x^n$, the {\em maximum
  likelihood\/} (ML) $P$ is the $P \in \M$ maximizing $P(x^n)$. For a
parametric model with parameter space $\Theta$, the maximum likelihood
{\em estimator\/} $\hat{\theta}$ is the function that, for each $n$,
maps $x^n$ to the $\theta \in \Theta$ that maximizes the likelihood
$P(x^n \mid \theta)$.  \index{ML} \index{maximum likelihood} The ML
estimator may be viewed as a `learning algorithm'. This is a procedure that,
when input a sample $x^n$ of arbitrary length, outputs a parameter or 
hypothesis
$P_n \in \M$. We say a learning algorithm is {\em consistent\/} relative to distance
measure $d$ if for
all $P^* \in \M$, if data are distributed according to $P^*$, then the
output $P_n$ converges to $P^*$ in the sense that $d(P^*,P_n) \rightarrow 0$ with
$P^*$-probability $1$.  Thus, if
$P^*$ is the `true' state of nature, then given enough data, the
learning algorithm will learn a good approximation of $P^*$ 
with very high probability.  \index{consistency}
\begin{example}{\bf [Markov and Bernoulli models]}
\label{ex:markov}
\rm 
Recall that a $k$-th order Markov chain on 
$\xspace = \{0,1\}$ is a probabilistic source such
that for every $n > k$,
\begin{multline}
\label{eq:markovdef}
P(X_n = 1 \mid X_{n-1} = x_{n-1}, \ldots, X_{n-k} = x_{n-k}) = \\
P(X_n = 1 \mid X_{n-1} = x_{n-1}, \ldots, X_{n-k} = x_{n-k}, \ldots,
X_1 = x_1).
\end{multline}
That is, the probability distribution on $X_n$ depends only on the $k$
symbols preceding $n$. Thus, there are $2^k$ possible distributions of
$X_n$, and each such distribution is identified with a {\em state\/}
of the Markov chain. To fully identify the chain, we also need to
specify the {\em starting state}, defining the
first $k$ outcomes $X_1, \ldots, X_k$. The $k$-th order {\em Markov model\/} is the set
of all $k$-th order Markov chains, i.e. all sources satisfying
(\ref{eq:markovdef}) equipped with a starting state. 

The special case of the $0$-th order Markov model is the {\em
  Bernoulli\/} or {\em biased coin\/} model, which we denote by
$\BernoulliN{0}$ 
%(it is customary to index models by their number of
%free parameters, hence the superscript `$1$' rather than `$0$'). 
We can parameterize the
Bernoulli model by a parameter $\theta \in [0,1]$ representing the
probability of observing a $1$. Thus, $\BernoulliN{0} = \{ P(\cdot \mid
\theta) \mid \theta \in [0,1]\}$, with $P(x^n \mid \theta)$ by
definition equal to
$$
P(x^n \mid \theta) = \prod_{i=1}^n P(x_i \mid \theta) = \theta^{n_{[1]}} 
(1 - \theta)^{n_{[0]}},$$
where $n_{[1]}$ stands for the number of $1$s, and $n_{[0]}$ 
for the number of $0$s in the sample.
Note that the Bernoulli model is i.i.d. The log-likelihood is given by
\index{log-likelihood}
\begin{equation}
\label{eq:loglik}
\log P(x^n \mid \theta) = n_{[1]} \log \theta + n_{[0]} \log (1- \theta).
\end{equation}
Taking the derivative of (\ref{eq:loglik}) with respect to $\theta$,
we see that for fixed $x^n$, the log-likelihood is maximized by
setting the probability of $1$ equal to the observed frequency. Since
the logarithm is a monotonically increasing function, the likelihood
is maximized at the same value: the ML estimator is given by 
$\hat{\theta}(x^n) = n_{[1]}/n$. 

Similarly, the first-order Markov model $\BernoulliN{1}$ 
can be parameterized by a
vector $\theta = (\theta_{[1|0]}, \theta_{[1|1]}) \in [0,1]^2$ together with a starting state in $\{0,1\}$. Here
$\theta_{[1|j]}$ represents the probability of observing a $1$
following the symbol $j$. The log-likelihood is given by
%\begin{equationr}
%\label{eq:loglikb}
$$
\log P(x^n \mid \theta) = n_{[1|1]} \log \theta_{[1|1]} 
+ n_{[0|1]} \log (1- \theta_{[1|1]}) + 
n_{[1|0]} \log \theta_{[1|0]} 
+ n_{[0|0]} \log (1- \theta_{[1|0]}),
$$
%\end{equationr}
$n_{[i|j]}$ denoting the number of times outcome $i$ is observed in
state (previous outcome) $j$. This is maximized by
setting $\hat{\theta} = (\hat{\theta}_{[1|0]}, \hat{\theta}_{[1|1]})$,
with $\hat{\theta}_{[i|j]} = n_{[i|j]} = n_{[ji]}/n_{[j]}$ set to the conditional
  frequency of $i$ preceded by $j$. In general, a
  $k$-th order Markov chain has $2^k$ parameters and the corresponding likelihood is maximized by setting the parameter
  $\theta_{[i|j]}$ equal to the number of times $i$ was observed in
  state $j$ divided by the number of times the chain was in state $j$.
\end{example}
Suppose now we are given data $D = x^n$ and we want to find
the Markov chain that best explains $D$. Since we do not want to
restrict ourselves to chains of fixed order, we run a large risk of
overfitting: simply picking, among all Markov chains of each order,
the ML Markov chain that maximizes the probability of the data, we
typically end up with a chain of order $n-1$ with starting state given
by the sequence $x_1, \ldots, x_{n-1}$, and $P(X_n = x_n \mid X^{n-1}
= x^{n-1}) = 1$. Such a chain will assign probability 1 to $x^n$. Below we
show that MDL makes a more reasonable choice.
% in contrast to ML, MDL will select a Markov chain that
%is a compromise between complexity (Markov chain order) and
%fit.
\section{Crude MDL}
\label{sec:crude}
Based on the information-theoretic
(Section~\ref{sec:prob}) and statistical
(Section~\ref{sec:statistics}) preliminaries discussed before, we now formalize a first,
crude version of MDL. 

Let $\M$ be a class of probabilistic sources (not necessarily Markov
chains). Suppose we observe a sample $D = (x_1, \ldots, x_n) \in
\xspace^n$.  Recall `the crude\endnote{The terminology `crude MDL' is not standard.
      It is introduced here for pedagogical reasons, to make clear the
      importance of having a single, unified principle for designing
      codes.  It should be noted that Rissanen's and Barron's early
      theoretical papers on MDL already contain such principles,
      albeit in a slightly different form than in their recent papers.
      Early practical applications \cite{QuinlanR89,Grunwald96d} often
      do use {\em ad hoc\/} two-part codes which really are `crude' in
      the sense defined here.}
two-part code MDL Principle' from Section~\ref{sec:modselintro},
page~\pageref{box:hypsel}: 
 \ownbox{Crude\endnote{See the previous endnote.}, Two-part
  Version of MDL Principle}{Let $\HN{1}, \HN{2},
  \ldots$ be a set of candidate models. The best point hypothesis
  $H \in \HN{1} \cup \HN{2} \cup \ldots $ to explain data $D$
  is the one which minimizes the sum $L(\phyp) + L(D|\phyp)$, where
\label{box:hypselb}
\begin{itemize}
\item $L(\phyp)$ is the length, in bits, of the description of the
  hypothesis; and
\item $L(D|\phyp)$ is the length, 
in bits, of the description of the data when encoded
with the help of the hypothesis.
\end{itemize}
The best {\em model\/} to explain $D$ is the smallest model
containing the selected $\phyp$.}
%\caption{\label{fig:twopart} The two-part MDL Principle: first, crude
%  implementation of the MDL ideas.
In this section, we implement this crude MDL Principle by giving a
precise definition of the 
terms $L(H)$ and $L(D|H)$. To make the first
term precise, we must design a code $C_1$ for encoding hypotheses $H$
such that $L(H) = L_{C_1}(H)$. For the second term, we must design a
set of codes $C_{2,H}$ (one for each $H \in \M$) such that for all $D
\in \xspace^n$, $L(D|H) = L_{C_{2,H}}(D)$. We start by describing the
codes $C_{2,H}$.
\subsection{Description Length of Data given Hypothesis}
Given a sample of size $n$, each hypothesis $H$ may be viewed as a
probability distribution on $\xspace^n$. We denote the
corresponding probability mass function by $P(\cdot \mid H)$. We need
to associate with $P(\cdot \mid H)$ a code, or really, just a
codelength function for $\xspace^n$. We already know that there exists
a code with length function $L$ such that for all $x^n \in \xspace^n$,
$L(x^n) = - \log P(x^n \mid H)$. This is the code that we will pick.
It is a  natural choice for two reasons:
\begin{enumerate}
\item With this choice, the code length $L(x^n \mid H)$ is equal to
  minus the log-likelihood of $x^n$ according to $H$, which is a
  standard statistical notion of `goodness-of-fit'. 
\item {\em If\/} the data turn out to be distributed according to $P$,
  {\em then\/} the code $L(\cdot \mid H)$ will uniquely minimize the
  expected code length (Section~\ref{sec:prob}).
\end{enumerate}
\begin{ownquote}
The second item implies that our choice is, in a sense, the only
reasonable choice\endnote{but see \cite{Grunwald98b}, Chapter 5 for
  more discussion.}. To see this, suppose $\M$ is a finite i.i.d.
model containing, say, $M$ distributions. 
Suppose we assign an arbitrary but finite
code length $L(H)$ to each $H \in \M$. Suppose $X_1, X_2, \ldots$ are
actually distributed i.i.d. according to some `true' $H^* \in \M$. By
the reasoning of Example~\ref{ex:lln}, we see that MDL will select the
true distribution $P(\cdot \mid H^*)$ for all large $n$, with
probability 1.  This means that MDL is {\em consistent\/} for finite
$\M$.  If we were to assign codes to distributions in some other
manner not satisfying $L(D \mid H) = - \log P(D \mid H)$, then there
would exist distributions $P(\cdot \mid H)$ such that $L(D|H) \neq - \log
P(D|H)$. But by Figure~\ref{fig:most}, there {\em must\/} be some 
distribution $P(\cdot \mid H')$ with $L(\cdot|H) = - \log P(\cdot \mid
H')$. Now let $\M = \{H, H' \}$ and suppose data are distributed
according to $P(\cdot \mid H')$. Then, by the reasoning of
Example~\ref{ex:lln},  MDL would select $H$ rather than 
$H'$ for all large $n$!
Thus, MDL would be inconsistent even in this simplest of all
imaginable cases -- there would then be no hope for good performance
in the considerably more complex situations we intend to use it
for\endnote{See Section~\ref{sec:philosophy} of Chapter~\ref{chap:survey} 
for a discussion on the
  role of consistency in MDL.}.
\end{ownquote}
\subsection{Description Length of Hypothesis}
In its weakest and crudest form, the two-part code MDL Principle does
not give any guidelines as to how to encode hypotheses (probability
distributions). Every code for encoding hypotheses is allowed, 
{\em as long as such a code does not change with the
  sample size $n$}.
\begin{ownquote}
  To see the danger in allowing codes to depend on $n$, consider the
  Markov chain example: if we were allowed to use different codes for
  different $n$, we could use, for each $n$, a code assigning a
  uniform distribution to all Markov chains of order $n-1$ with all
  parameters equal to $0$ or $1$. Since there are only a finite number
  ($2^{n-1})$ of these, this is possible. But then, for each $n$, $x^n
  \in \xspace^n$, MDL would select the ML Markov chain of order $n-1$.
  Thus, MDL would coincide with ML and, no matter how large $n$, we
  would overfit.
\end{ownquote}
\paragraph{Consistency of Two-part MDL}
Remarkably, if we fix an arbitrary code for all hypotheses, identical
for all sample sizes $n$, this is sufficient to make MDL
consistent\endnote{See, for example \cite{BarronC91}, \cite{Barron85}} for a wide variety of models, including the Markov
chains. For example, let $L$ be the length
function corresponding to some code for the Markov chains.  Suppose
some Markov chain $P^*$ generates the data such that $L(P^*) < \infty$
under our coding scheme. Then, broadly speaking, for every $P^*$ of
every order, with probability 1 there exists some $n_0$ such that for all samples larger
than $n_0$,  two-part MDL will select $P^*$ -- here $n_0$ may depend on
$P^*$ and $L$. 

While this results indicates that MDL may be doing something sensible,
it certainly does not justify the use of arbitrary codes - different
codes will lead to preferences of different hypotheses, and it is not
at all clear how a code should be designed that leads to good
inferences with small, practically relevant sample sizes.

\citeN{BarronC91} have developed a precise theory of how to design
codes $C_1$ in a `clever' way, anticipating the developments of
`refined MDL'. Practitioners have often simply used `reasonable'
coding schemes, based on the following idea. Usually there exists some
`natural' \commentout{Doen in boek moet hier natuurlijk uitleg over de
  vraag wat 'natural' is !}  decomposition of the models under
consideration, $\M = \bigcup_{k > 0} \MN{k}$ where the dimension of
$\MN{k}$ grows with $k$ but is not necessarily equal to $k$. In the
Markov chain example, we have $\Bernoulli = \bigcup \BernoulliN{k}$
where $\BernoulliN{k}$ is the $k$-th order, $2^k$-parameter Markov
model. Then {\em within\/} each submodel $\MN{k}$, we may use a
fixed-length code for $\theta \in \Theta^{(k)}$.  Since the set
$\Theta^{(k)}$ is typically a continuum, we somehow need to discretize
it to achieve this.
\begin{example}{\bf [a Very Crude Code for the Markov Chains]}
\rm
%Let $\Bernoulli$ be the set of all Markov chain distributions of each order
%with rational-valued parameters. There are several ways to encode such
%Markov chains.  First note that f
We can describe a Markov chain of order $k$ by first describing $k$,
and then describing a parameter vector $\theta \in [0,1]^{k'}$ with
$k' = 2^k$. We describe $k$ using our simple code for the integers
(Example~\ref{ex:integers}). This takes $2 \log k + 1$ bits. We now
have to describe the $k'$-component parameter vector.  We saw in
Example~\ref{ex:markov} that for any $x^n$, the best-fitting (ML)
$k$-th order Markov chain can be identified with $k'$ frequencies. It
is not hard to see that these frequencies are uniquely determined by
the counts $n_{[1|0\ldots 00]}, n_{[1|0\ldots 01]}, \ldots,
n_{[1|1\ldots11]}$. Each individual count must be in the
$(n+1)$-element set $\{0, 1, \ldots, n \}$. Since we assume $n$ is
given in advance\endnote{Strictly speaking, the assumption that $n$ is
  given in advance (i.e., both encoder and decoder know $n$)
  contradicts the earlier requirement that the code to be used for
  encoding hypotheses is not allowed to depend on $n$. Thus, strictly
  speaking, we should first encode some $n$ explicitly, using $2 \log
  n + 1$ bits (Example~\ref{ex:integers}), and then pick the $n$
  (typically, but not necessarily equal to the actual sample size)
  that allows for the shortest three-part codelength of the data
  (first encode $n$, then $(k,\theta)$, then the data). In practice
  this will not significantly alter the chosen hypothesis, unless for
  some quite special data sequences.}, we may use a simple
fixed-length code to encode this count, taking $\log (n+1)$ bits
(Example~\ref{ex:uniform}).  Thus, once $k$ is fixed, we can describe
such a Markov chain by a uniform code using $k' \log (n+1)$ bits.
With the code just defined we get for any $P \in \Bernoulli$, indexed
by parameter $\Theta^{(k)}$,
$$
%  \begin{equationr}
%\label{eq:markovtp}
  L(P) = L(k, \Theta^{(k)}) = 2 \log k + 1 + k \log (n+1),
$$
%\end{equationr}
so that with these codes, MDL tells us to pick the $k, \theta^{(k)}$ minimizing
\begin{equation}
\label{eq:twopartmdl}
L(k, \theta^{(k)}) + L(D \mid k, \theta^{(k)}) = 2 \log k + 1 + k \log (n+1)
- \log P(D \mid k, \theta^{(k)}),
\end{equation}
where the $\theta^{(k)}$ that is chosen will be equal to the ML estimator for $\MN{k}$.
\end{example}
\paragraph{Why (not) this code?}
We may ask two questions about this code. First, why did we only
reserve codewords for $\theta$ that are potentially ML estimators for
the given data? The reason is that, given $k' = 2^k$, the codelength
$L(D \mid k,\theta^{(k)})$ is minimized by $\hat{\theta}^{(k)}(D)$,
the ML estimator within $\theta^{(k)}$. Reserving codewords for
$\theta \in [0,1]^{k'}$ that cannot be ML estimates would only serve
to lengthen $L(D \mid k, \theta^{(k)})$ and can never shorten $L(k,
\theta^{(k)})$. Thus, the total description length needed to encode
$D$ will increase. Since our stated goal is to minimize description
lengths, this is undesirable. 

However, by the same logic we may also ask whether we have not
reserved {\em too many\/} codewords for $\theta \in [0,1]^{k'}$. And
in fact, it turns out that we have: the distance between two adjacent
ML estimators is $O(1/n)$. Indeed, if we had used a coarser precision,
only reserving codewords for parameters with distances $O(1 /
\sqrt{n})$, we would obtain smaller code lengths -
(\ref{eq:twopartmdl}) 
would become 
\begin{equation}
\label{eq:twopartmdlb}
L(k, \theta^{(k)}) + L(D \mid k, \theta^{(k)}) = 
- \log P(D \mid k, \hat{\theta}^{(k)}) + \frac{k}{2} \log n + \const_k,
\end{equation}
where $\const_k$ is a small constant depending on $k$, but not $n$ \cite{BarronC91}. In Section~\ref{sec:refined} we show that (\ref{eq:twopartmdlb}) is in some sense `optimal'.
\paragraph{The Good News and the Bad News}
The good news is (1) we have found a principled, non-arbitrary manner
to encode data $D$ given a probability distribution $H$, namely, to
use the code with lengths $- \log P(D\mid H)$; and (2),
asymptotically, {\em any\/} code for hypotheses will lead to a
consistent criterion. The bad news is that we have not found clear
guidelines to design codes for hypotheses $H \in \M$. We found some
intuitively reasonable codes for Markov chains, and we then reasoned that these could be somewhat `improved', but what is
conspicuously lacking is a sound theoretical {\em principle\/} for
designing and improving codes.

We take the good news to mean that our idea may be worth pursing
further. We take the bad news to mean that we do have to modify or
extend the idea to get a meaningful, non-arbitrary and practically
relevant model selection method. Such an extension was already suggested in
Rissanen's early works \cite{Rissanen78,Rissanen83} 
and refined by \citeN{BarronC91}. 
However, in these works, the principle was still restricted to two-part
codes. To get a fully satisfactory solution, we need to move to
`universal codes', of which the two-part codes are merely a special case.
\section{Information Theory II: Universal Codes and Models}
% Doen GEBRUIK GEEN D voor data in deze sectie!
\label{sec:universal}
We have just indicated why the two-part code formulation of MDL needs
to be refined. It turns out that the key concept we need is that of
{\em universal coding}. Broadly speaking, a code $\luniv$ that is
universal relative to a set of candidate codes ${\cal L}$ allows us to
compress every sequence $x^n$ almost as well as the code in ${\cal L}$
that compresses that particular sequence most. Two-part codes are
universal (Section~\ref{sec:twopartuni}), but there exist other
universal codes such as the Bayesian mixture code
(Section~\ref{sec:unimodel}) and the Normalized Maximum Likelihood
(NML) code (Section~\ref{sec:optimal}). We also discuss {\em universal
  models}, which are just the probability distributions corresponding
to universal codes.
In this section, we are not concerned with learning from data; we only care
about compressing data as much as possible. We reconnect our findings
with learning in Section~\ref{sec:refined}.
\paragraph{Coding as Communication}
\index{messages} \index{encoder} \index{decoder} Like many other topics in coding, `universal coding' can best be
explained if we think of descriptions as {\em messages\/}: 
we can always view a
description as a message that some sender or {\em encoder}, say Mr. A,
sends to some receiver or {\em decoder}, say Mr. B. Before sending any
messages, Mr. A and Mr. B meet in person. They agree on the set of
messages that A may send to B. Typically, this will be the set $\xspace^n$ 
of sequences $x_1, \ldots, x_n$, where each $x_i$ is an outcome in the space 
$\xspace$ . They also agree upon a (prefix) code that will be used by
A to send his messages to B. Once this has been done, A and B go back
to their respective homes and A sends his messages to B in the form of
binary strings. The unique decodability property of prefix codes
implies that, when B receives a message, he should always be able to
decode it in a unique manner.
\paragraph{Universal Coding}
Suppose our encoder/sender is about to observe a sequence $x^n \in
\xspace^n$ which he plans to compress as much as possible.
Equivalently, he wants to send an encoded version of $x^n$ to the
receiver using as few bits as possible. Sender and receiver have a set
of {\em candidate codes\/} $\codelengths$ for $\xspace^n$
available\endnote{As explained in Figure~\ref{fig:new}, we identify
  these codes with their length functions, which is the only aspect we
  are interested in.}. They believe or hope that one of these codes
will allow for substantial compression of $x^n$. However, they must
decide on a code for $\xspace^n$ before sender observes the actual
$x^n$, and they do not know {\em which\/} code in $\codelengths$ will
lead to good compression of the actual $x^n$. What is the best thing
they can do?  They may be tempted to try the following: upon seeing
$x^n$, sender simply encodes/sends $x^n$ using the $L \in
\codelengths$ that minimizes $L(x^n)$ among all $L \in \codelengths$.
But this naive scheme will not work: since decoder/receiver does not
know what $x^n$ has been sent before decoding the message, he does not
know which of the codes in $\codelengths$ has been used by
sender/encoder.  Therefore, decoder cannot decode the message: the
resulting protocol does not constitute a uniquely decodable, let alone
a prefix code.  Indeed, as we show below, in general {\em no\/} code
$\luniv$ exists such that for all $x^n \in \xspace^n$, $\luniv(x^n)
\leq \min_{L \in \codelengths} L(x^n)$: in words, there exists no code
which, no matter what $x^n$ is, always mimics the best code for $x^n$.
\begin{example}
\label{ex:nineber}
\rm
Suppose we think that our sequence can be reasonably well-compressed
by a code corresponding to some biased coin model. For simplicity, we restrict ourselves to a finite number of such models. Thus, let 
$\codelengths = \{L_1, \ldots, L_9 \}$ where $L_1$ is the code length function 
corresponding to the Bernoulli model $P(\cdot \mid \theta)$ with
parameter $\theta = 0.1$, $L_2$ corresponds to $\theta = 0.2$ and so on.
From (\ref{eq:loglik}) we see that, for example,
$$
\begin{array}{l}
L_8(x^n) = - \log P(x^n | 0.8) = - n_{[0]} \log 0.2 -
n_{[1]} \log 0.8  \\ 
L_9(x^n) = - \log P(x^n | 0.9) = - n_{[0]} \log 0.1 - n_{[1]} \log
0.9. 
\end{array}
$$
Both $L_8(x^n)$ and $L_9(x^n)$ are linearly increasing in the
number of $1$s in $x^n$. However, if the frequency $n_1/n$ is
approximately $0.8$, then $\min_{L \in \codelengths} L(x^n)$ will be
achieved for $L_8$. If $n_1/n \approx 0.9$ then $\min_{L \in
  \codelengths} L(x^n)$ is achieved for $L_9$. More generally, if
$n_1/n \approx j/10$ then $L_j$ achieves the minimum\endnote{The
  reason is that, in the full Bernoulli model with parameter $\theta
  \in [0,1]$, the maximum likelihood estimator is given by $n_1/n$,
  see Example~\ref{ex:markov}. Since the likelihood $\log P(x^n \mid
  \theta)$ is a continuous function of $\theta$, this implies that if
  the frequency $n_1/n$ in $x^n$ is approximately (but not precisely)
  $j/10$, then the ML estimator in the restricted model $\{ 0.1,
  \ldots, 0.9 \}$ is still given by $\hat{\theta} = j/10$. Then $\log
  P(x^n |\theta)$ is maximized by $\hat{\theta} = j/10$, so that the
  $L \in \codelengths$ that minimizes codelength corresponds to
  $\theta = j/10$.}.  We would like to send $x^n$ using a code
$\luniv$ such that for all $x^n$, we need at most $\hat{L}(x^n)$ bits,
where $\hat{L}(x^n)$ is defined as $\hat{L}(x^n) := \min_{L \in
  \codelengths} L(x^n)$. Since $- \log$ is monotonically decreasing,
$\hat{L}(x^n) = - \log P(x^n \mid \hat{\theta}(x^n))$.  We already
gave an informal explanation as to why a code with lengths
$\hat{L}$ does not exist. We can now explain this more formally as
follows: if such a code were to exist, it would correspond to some
distribution $\puniv$. Then we would have for all $x^n$, $\luniv(x^n)
= - \log \puniv(x^n)$. But, by definition, for all $x^n \in {\cal
  X}^n$, $\luniv(x^n) \leq \hat{L}(x^n) = - \log P(x^n |
\hat{\theta}(x^n))$ where $\hat{\theta}(x^n) \in \{0.1, \ldots, 0.9
\}$. Thus we get for all $x^n$, $
%$
- \log \puniv(x^n) \leq 
- \log P(x^n \mid \hat{\theta}(x^n)) \text{\ or \ }
\puniv(x^n) \geq P(x^n \mid \hat{\theta}(x^n)), 
%$
$
so that, since $|{\cal L}| > 1$,
% contains more than one code, 
\begin{equation}
\label{eq:largerone}
\sum_{x^n} \puniv(x^n) \geq \sum_{x^n} P(x^n \mid \hat{\theta}(x^n)) =
\sum_{x^n} \max_{\theta} P(x^n \mid \theta) > 1,
\end{equation}
where the last inequality follows because for any two $\theta_1,
\theta_2$ with $\theta_1 \neq \theta_2$, there is at least one $x^n$
with $P(x^n \mid \theta_1) > P(x^n \mid \theta_2)$.
(\ref{eq:largerone}) says that $\puniv$ is not a probability
distribution. It follows that $\luniv$ cannot be a codelength function.
The argument can be extended beyond the Bernoulli model of the example
above: as long as $|{\cal L}| > 1$, 
%contains more than one code, 
and all codes in ${\cal L}$ correspond to a non-defective
distribution, (\ref{eq:largerone}) must still hold, so that there
exists no code $\luniv$ with $\luniv(x^n) = \hat{L}(x^n)$ for all
$x^n$. The underlying reason that no such code exists is the fact that
probabilities must sum up to something $\leq 1$; or equivalently, that
there exists no coding scheme assigning short code words to many
different messages -- see Example~\ref{ex:nocompression}.
\end{example}
Since there exists no code
which, no matter what $x^n$ is, always mimics the best code for $x^n$,
it may make sense to look for the next best thing: does
there exist a code which, for all $x^n \in \xspace^n$, is `nearly' (in
some sense) as good as $\hat{L}(x^n)$? 
It turns out that in many cases, the answer is {\em
  yes}: there typically exists codes
$\luniv$ such that no matter what $x^n$ arrives, $\luniv(x^n)$
is not much larger than $\hat{L}(x^n)$, which may be viewed as the code that is best `with hindsight' (i.e., after
seeing $x^n$). Intuitively, codes which
satisfy this property are called universal codes - a more precise definition
follows below. 
The first (but perhaps not foremost) example of a
universal code is the {\em two-part code\/} that we have encountered
in Section~\ref{sec:crude}.
\subsection{Two-part Codes as simple Universal Codes} 
\label{sec:twopartuni}
\begin{example}{\bf [finite $\codelengths$]}
\label{ex:nineber2}
\rm Let $\codelengths$ be as in Example~\ref{ex:nineber}. We can
devise a code $\lunivtp$ for all $x^n \in \samplespace^n$ as follows:
to encode $x^n$, we first encode the $j \in \{1, \ldots, 9\}$ such
that $L_j(x^n) = \min_{L \in \codelengths} L(x^n)$, 
%(if there is more than one such $j$, we take the smallest) 
using a uniform code. This takes $\log 9 $ bits.  We then encode $x^n$
itself using the code indexed by $j$. This takes $L_j$ bits. Note that
in contrast to the naive scheme discussed in Example~\ref{ex:nineber},
the resulting scheme properly defines a prefix code: a decoder can
decode $x^n$ by first decoding $j$, and then decoding $x^n$ using
$L_j$. Thus, for {\em every possible\/} $x^n \in \samplespace^n$, we
obtain
$$
\lunivtp(x^n) = \min_{L \in \codelengths} L(x^n) + \log 9.
$$
%It is easy to see that in our example, 
For all $L \in
\codelengths$, $\min_{x^n} L(x^n)$ grows linearly in $n$:
$\min_{\theta,x^n} - \log P(x^n \mid \theta) = - n \log 0.9 \approx
0.15 n$.  Unless $n$ is {\em very\/} small, no matter what $x^n$
arises, the extra number of bits we need using $\lunivtp$
compared to $\hat{L}(x^n)$ is negligible.
\end{example}
More generally, let $\codelengths = \{L_1, \ldots, L_M \}$ where $M$
can be arbitrarily large, and the $L_j$ can be any codelength
functions we like; they do not necessarily represent Bernoulli
distributions any more.  By the reasoning of
Example~\ref{ex:nineber2}, there exists a (two-part) code such that
for {\em all\/} $x^n \in \samplespace^n$,
\begin{equation}
\label{eq:finite}
\lunivtp(x^n) = \min_{L \in \codelengths} L(x^n) + \log M.
\end{equation}
In most applications $\min L(x^n)$ grows linearly in
$n$, and we see from (\ref{eq:finite}) that, as soon as $n$ becomes
substantially larger than $\log M$, the relative difference in
performance between our universal code and $\hat{L}(x^n)$ becomes negligible.
In general, we do not always want to use a uniform code for the
elements in ${\cal L}$; note that any arbitrary code on ${\cal L}$
will give us an analogue of (\ref{eq:finite}), but with a worst-case
overhead larger than $\log M$ - corresponding to the largest
codelength of any of the elements in ${\cal L}$.
\begin{example}{\bf [Countably Infinite $\codelengths$]}
\label{ex:countable}
  \rm We can also construct a 2-part code for arbitrary countably
  infinite sets of codes $\codelengths = \{ L_1, L_2, \ldots \}$: we
  first encode some $k$ using our simple code for the integers
  (Example~\ref{ex:integers}).  With this code we need $2 \log k +1 $
  bits to encode integer $k$. We then encode $x^n$ using the code
  $L_k$. $\lunivtp$ is now defined as the code we get if, for any
  $x^n$, we encode $x^n$ using the $L_k$ minimizing the total two-part
  description length $ 2 \log k +1 + L_k(x^n)$.
  
  In contrast to the case of finite $\codelengths$, there does {\em
    not\/} exist a constant $\const$  any more
  such that for all $n, x^n \in \xspace^n$, $\lunivtp(x^n) \leq
  \inf_{L \in \codelengths} L(x^n) + c$. Instead we have the following
  weaker, but still remarkable property: for all $k$, all $n$, all
  $x^n$, $\lunivtp(x^n) \leq L_k(x^n) + 2 \log k + 1$, so that also,
$$\lunivtp(x^n) \leq \inf_{L \in \{L_1, \ldots, L_k \}} L(x^n) + 2
\log k +1.$$
For any $k$, as $n$ grows larger, the code
$\lunivtp$ starts to mimic whatever $L \in \{ L_1, \ldots, L_k \}$
compresses the data most. However, the larger $k$, the larger $n$ has
to be before this happens.
\end{example}
\subsection{From Universal Codes to Universal Models}
\label{sec:unimodel}
Instead of postulating a set of candidate codes
$\codelengths$, we may equivalently postulate a set $\M$ of
candidate probabilistic sources, such that $\codelengths$ is the set
of codes corresponding to $\M$. We already implicitly did this in
Example~\ref{ex:nineber}.

The reasoning is now as follows: we think that one of the $P \in \M$
will assign a high likelihood to the data to be observed. Therefore we
would like to design a code that, for all $x^n$ we might observe,
performs essentially as well as the code corresponding to the
best-fitting, maximum likelihood (minimum codelength) $P \in \M$ for
$x^n$. Similarly, we can think of universal codes such as the two-part
code in terms of the (possibly defective, see Section~\ref{sec:prob}
and Figure~\ref{fig:most})) {\em distributions\/}
corresponding to it. Such distributions corresponding to universal
codes are called {\em universal models}. The use of mapping universal codes back to distributions is illustrated by the {\em Bayesian universal model} which we now introduce.
\begin{ownquote}{\bf Universal model: Twice Misleading Terminology}
  The words `universal' and `model' are somewhat of a misnomer:
  first, these codes/models are only `universal' relative to a
  restricted `universe' $\M$. Second, the
  use of the word `model' will be very confusing to statisticians, who
  (as we also do in this paper) call a family of distributions such as
  $\M$ a 'model'. But the phrase originates from information theory,
  where a `model' often refers to a single distribution rather than a
  family. Thus, a `universal model' is a single distribution,
  representing a statistical `model' $\M$.
\end{ownquote}
\begin{example}{ \bf [Bayesian Universal Model]}
\label{ex:bayes}
\rm  Let $\M$ be a finite or countable set of probabilistic sources,
  parameterized by some parameter set $\Theta$. Let $W$ be a
  distribution on $\Theta$. Adopting terminology from Bayesian
  statistics, $W$ is usually called a {\em prior distribution}.  We
  can construct a new probabilistic source $\punivbayes$ by taking a
  weighted (according to $W$) average or mixture over the
  distributions in $\M$. That is, we define for all $n$, $x^n \in \xspace$, 
\begin{equation}
\label{eq:bayes}
\punivbayes(x^n) := \sum_{\theta \in \Theta} P(x^n \mid \theta) W(\theta).
\end{equation}
\index{prior distribution}It is easy to check that $\punivbayes$ is a
probabilistic source according to our definition. In case $\Theta$ is
continuous, the sum gets replaced by an integral but otherwise nothing
changes in the definition. 
In Bayesian statistics, $\punivbayes$ is called the {\em
  Bayesian marginal likelihood\/} or {\em Bayesian mixture\/}
\cite{BernardoS94}. \index{marginal likelihood} \index{Bayes evidence}
To see that $\punivbayes$ is a universal model, note that
for all $\theta_0 \in \Theta$,
\begin{equation}
\label{eq:bayesuniversal}
- \log \punivbayes(x^n) := 
- \log \sum_{\theta \in \Theta} P(x^n \mid \theta) W(\theta) \leq
- \log  P(x^n \mid \theta_0) + \const_{\theta_0}
\end{equation}
where the inequality follows because a sum is at least as large as
each of its terms, and $\const_{\theta} = - \log W(\theta)$ depends on
$\theta$ but not on $n$. Thus, $\punivbayes$ is a universal model or
equivalently, the code with lengths $- \log \punivbayes$ is a
universal code. Note that the derivation in 
(\ref{eq:bayesuniversal}) only works if $\Theta$ is finite or
countable; the case of continuous $\Theta$ is treated in Section~\ref{sec:refined}.
\end{example}
\paragraph{Bayes is Better than Two-part} The Bayesian model is in a sense superior to the two-part code.
Namely, in the two-part code we first encode an element of $\M$ or its
parameter set $\Theta$ using some code $L_0$. Such a
code must correspond to some `prior' distribution $W$ on $\M$ so that
the two-part code gives codelengths
\begin{equation}
\label{eq:mdlworse}
\lunivtp(x^n) = \min_{\theta \in \Theta} 
- \log P(x^n |\theta) - \log W(\theta) 
\end{equation}
where $W$ depends on the specific code $L_0$ that was used. Using the
Bayes code with prior $W$, we get as in (\ref{eq:bayesuniversal}),
$$
%\begin{equation}
%\label{eq:bayesbetter}
- \log \punivbayes(x^n) = - \log \sum_{\theta \in \Theta} P(x^n \mid \theta) W(\theta) \leq \min_{\theta \in \Theta} 
- \log P(x^n |\theta) - \log W(\theta).
%\end{equation}
$$
The inequality becomes strict whenever
$P(x^n | \theta) > 0$ for more than one value of $\theta$. Comparing
to (\ref{eq:mdlworse}), we see that 
in general the Bayesian code is preferable over the two-part code: for
all $x^n$ it never assigns codelengths  larger than $\lunivtp(x^n)$, and
in many cases it assigns strictly shorter codelengths for some $x^n$.
But this raises two
important issues: (1) what exactly do we mean by `better' anyway? (2)
can we say that `some prior distributions are better than others'?
These questions are answered below.
\subsection{NML as an {\em Optimal Universal Model}}
\label{sec:optimal}
We can measure the performance of universal
models relative to a set of candidate sources $\M$ using the {\em
  regret\/}:
\begin{definition}{ \bf [Regret]}
\label{def:regret}
Let $\M$ be a class of probabilistic sources.
Let $\puniv$ be a probability distribution on $\xspace^n$ ($\puniv$ 
is not necessarily in $\M$).  For given $x^n$, the {\em regret\/} of
$\puniv$ relative to $\M$ is defined as
\begin{equation}
\label{eq:regret}
- \log \puniv(x^n) - \min_{P \in \M} \{ - \log P(x^n ) \}.
\end{equation}
%The {\em average regret-per-outcome\/} is given by 
%$$
%\frac{1}{n} \biggl\{ 
%- \log \puniv(x^n) - \{ - \log P(x^n) \} \biggr\}.
%$$ 
\end{definition}
The regret of $\puniv$ relative to $\M$ for $x^n$ is the
additional number of bits needed to encode $x^n$ using the
code/distribution $\puniv$, as compared to the number of bits that had
been needed if we had used code/distribution in $\M$ that was {\em
  optimal (`best-fitting') with hind-sight}. For simplicity, from now on 
we tacitly assume that for all the
models $\M$ we work with, there is a single $\hat{\theta}(x^n)$
maximizing the likelihood for every $x^n \in \xspace^n$. In that case (\ref{eq:regret}) simplifies to
$$
- \log \puniv(x^n) - \{ - \log P(x^n \mid \hat{\theta}(x^n)) \}.
$$
We would like to measure the quality of a
universal model $\puniv$ in terms of its regret. However, $\puniv$ may
have small (even $< 0$) regret for some $x^n$, and very large regret
for other $x^n$. We must somehow find a measure of quality that takes into
account {\em all\/} $x^n \in {\cal X}^n$. We take a
worst-case approach, and look for universal models $\puniv$ with small {\em worst-case\/} regret, where the worst-case is over all sequences. Formally,
the {\em maximum\/} or {\em worst-case regret\/} 
of $\puniv$ relative to $\M$ is defined as
$$
\maxregret(\puniv ) \isbydefinition \max_{x^n \in \xspace^n} \bigl\{ 
- \log \puniv(x^n) -  \{ - \log P(x^n \mid \hat{\theta}(x^n)) \} \bigr\}.
$$
If we use $\maxregret$ as our quality measure, then the `optimal' universal model relative to $\M$, for given sample size $n$, is the distribution minimizing
\begin{equation}
\label{eq:minimax}
\min_{\puniv} \  \maxregret(\puniv) =
\min_{\puniv} \ \max_{x^n \in \xspace^n} \ 
 \bigl\{ 
- \log \puniv(x^n) -  \{ - \log P(x^n \mid \hat{\theta}(x^n)) \} 
\bigr\}
\end{equation}
where the minimum is over {\em all\/} defective distributions on
$\xspace^n$. The $\puniv$ minimizing (\ref{eq:minimax})
corresponds to the code minimizing the additional number of bits
compared to code in $\M$ that is best in hindsight in the
worst-case over all possible $x^n$. It turns out that we can solve for $\puniv$ in (\ref{eq:minimax}). To this end, we first define the {\em complexity\/} of a given model $\M$ as
\begin{equation}
\label{eq:complexity}
\complexity_n(\M) \isbydefinition 
\log {\sum_{x^n \in \xspace^n} P(x^n \mid \hat{\theta}(x^n))}.
\end{equation}
\index{complexity}This quantity plays a fundamental role in refined
MDL, Section~\ref{sec:refined}.  To get a first idea of why
$\complexity_n$ is called model complexity, note that the more sequences
$x^n$ with large $P(x^n \mid \hat{\theta}(x^n))$, the larger
$\complexity_n(\M)$. In other words, the more sequences that can be fit
well by an element of $\M$, the larger $\M$'s  complexity.
\begin{proposition}{\bf 
%[Minimax, {Shtarkov} or {NML\/} Universal Code]
\cite{Shtarkov87}}
Suppose that $\complexity_n(\M)$ is finite.
Then the minimax regret (\ref{eq:minimax}) is uniquely achieved for the distribution $\punivnml$ given by
\begin{equation}
\label{eq:punivnml}
\punivnml(x^n) \isbydefinition \frac{P(x^n \mid \hat{\theta}(x^n))}{\sum_{y^n \in 
\xspace^n} P(y^n \mid \hat{\theta}(y^n))}.
\end{equation}
{\em The distribution $\punivnml$ is known as the {\em Shtarkov
    distribution\/} or the {\em normalized maximum likelihood\/} (NML)
  distribution}.
\end{proposition}
\paragraph{Proof} Plug in $\punivnml$ in (\ref{eq:minimax}) and notice that for
all $x^n \in {\cal X}^n$,
\begin{equation}
\label{eq:minimaxb}
- \log \puniv(x^n) -  \{ - \log P(x^n \mid \hat{\theta}(x^n)) \}
= \maxregret(\puniv) = \complexity_n(\M),
\end{equation}
so that $\punivnml$ achieves the {\em same\/} regret, equal to
$\complexity_n(\M)$, {\em no matter what $x^n$ actually obtains}.
Since every distribution $P$ on ${\cal X}^n$ with $P \neq
\punivnml$ must satisfy $P(z^n) < \punivnml(z^n)$ for at least one
$z^n \in {\cal X}^n$, it follows that 
\begin{multline}
\maxregret(P) \geq  - \log P(z^n) + \log P(z^n \mid \hat{\theta}(z^n)) >
\\
- \log \punivnml(z^n) + \log P(z^n \mid \hat{\theta}(z^n)) 
= \maxregret(\punivnml). \nonumber
\index{Shtarkov, Y.}  \index{normalized maximum likelihood}
\index{minimax regret} 
\end{multline} \mbox{$\Box$}

\ 
\\
\noindent
$\punivnml$ is quite literally a `normalized
maximum likelihood' distribution: it tries to assign to each $x^n$
the probability of $x^n$ according to the ML distribution for
$x^n$. By (\ref{eq:largerone}), this is not possible: the
resulting `probabilities' add to something larger than $1$. But we can
normalize these `probabilities' by dividing by their sum $\sum_{y^n \in 
\xspace^n} P(y^n \mid \hat{\theta}(y^n))$, and then 
we obtain a probability distribution on ${\cal X}^n$ after all. 

Whenever ${\cal X}$ is finite, the sum $\complexity_n(\M)$ is finite
so that the NML distribution is well-defined. If ${\cal X}$ is countably
infinite or continuous-valued, the sum $\complexity_n(\M)$ may be
infinite and then the NML distribution may be undefined. In that case, 
there exists {\em no\/} universal model
achieving constant regret as in (\ref{eq:minimaxb}). If $\M$ is
parametric, then $\punivnml$ is typically well-defined as long as we
suitably restrict the
parameter space. The parametric case forms the basis of `refined MDL'
and will be discussed at length in the next section. 
%A summary of the
%present section can be found in the box on page~\pageref{box:universal}.
\begin{figure}[h]
\ownbox{Summary: Universal Codes and Models}{
\label{box:universal}
Let $\M$ be a family of
  probabilistic sources.  A {\em universal
model in an individual sequence sense\endnote{
What we call `universal model' in this text is known in the literature as 
a `universal model in the individual sequence sense' -- there also
exist universal models in an `expected sense', see
Section~\ref{sec:whatis}. These lead to
slightly different versions of MDL.} relative to $\M$\/}, in this
text simply called a `universal model for $\M$', is a sequence
of distributions $\puniv^{(1)}, \puniv^{(2)}, \ldots$ on ${\cal X}^1, {\cal
  X}^2, \ldots$ respectively, such that for all $P \in \M$, 
for all $\epsilon > 0$,
$$
%\begin{equationr}
%\label{eq:universal}
\max_{x^n \in {\cal X}^n} \frac{1}{n} \biggl\{
- \log \puniv^{(n)}(x^n) - [ - \log P(x^n)]  \biggr\}
\leq \epsilon \text{\ as \ } n \rightarrow \infty.
%\end{equationr}
$$
Multiplying both sides with $n$ we see that $\puniv$ is universal if
for every $P \in \M$, the codelength difference $- \log \puniv(x^n) +
\log P(x^n)$ increases sublinearly in $n$. If $\M$ is finite, then the
two-part, Bayes and NML distributions are universal in a very strong
sense: rather than just increasing sublinearly, the codelength difference
is bounded by a constant.

We already discussed two-part,
Bayesian and 
minimax optimal (NML) universal models, but there several other types. 
We mention
prequential universal models (Section~\ref{sec:prequential}), the {\em
  Kolmogorov\/} universal model, {\em
  conditionalized\/} two-part codes \cite{Rissanen01} \/ and
Cesaro-average codes \cite{BarronRY98}.}\vspace*{-0.7 cm}
\end{figure}
% We will treat the parametric
%case in detail in the next section when discussing refined MDL. 
%We end this section with a formal definition and summary:
\commentout{Barron's type of universal models: 2 subtypes
1: expected regret (compare to $\hat{\theta}$), 2. redundancy (compare
to true $\theta$). In both cases (? see Barron Rissanen Yu 98) the
Bayes universal model with a certain prior is now minimax! This model
has the advantage that it is a random process, and that it is more
direct generalization of Shannon.} 
\section{Simple Refined MDL and its Four Interpretations}
\label{sec:refined}
In Section~\ref{sec:crude}, we indicated that `crude' MDL needs to be refined. In Section~\ref{sec:universal} we introduced universal models. 
We now show how universal models, in particular the minimax optimal universal model $\punivnml$, can be used to define a refined version of MDL model selection. Here we only discuss the simplest case: 
suppose we are given data $D = (x_1, \ldots, x_n)$ and two
models $\MN{1}$ and $\MN{2}$ such that $\complexity_n(\MN{1})$ and
$\complexity_n(\MN{2})$ (Equation~\ref{eq:complexity}) are both
finite. For example, we could have some binary data and $\MN{1}$ and
$\MN{2}$ are the first- and second-order Markov models
(Example~\ref{ex:markov}), both considered possible explanations for
the data.  We show how to deal
with an infinite number of models and/or models with infinite
$\complexity_n$ in Section~\ref{sec:unify}.

Denote by
$\punivnml(\cdot \mid \MN{j})$ the NML distribution on ${\cal X}^n$
corresponding to model $\MN{j}$.  Refined MDL tells us to pick the
model $\MN{j}$ maximizing the {\em normalized\/} maximum likelihood
$\punivnml(D\mid \MN{j})$, or, by (\ref{eq:punivnml}), equivalently, minimizing
\begin{equation}
\label{eq:refinedmdl}
- \log \punivnml(D \mid \MN{j}) =  - \log P(D \mid \hat{\theta}^{(j)}(D)) 
+ \complexity_n (\MN{j})
\end{equation}
From a coding theoretic point of view, we associate with each $\MN{j}$
a code with lengths $\punivnml(\cdot \mid \MN{j})$, and we pick the model
minimizing the codelength of the data. The codelength 
$- \log \punivnml(D \mid \MN{j})$ 
has been called the {\em stochastic complexity of the data $D$
  relative to model $\MN{j}$\/} \cite{Rissanen87}, whereas
$\complexity_n(\MN{j})$ is called the {\em parametric complexity\/} or
{\em model cost\/} of $\MN{j}$ (in this survey we simply call it
`complexity'). We have already indicated in the previous section that
$\complexity_n(\MN{j})$ measures something like the `complexity' of
model $\MN{j}$. On the other hand, $- \log P(D \mid \hat{\theta}^{(j)}(D))$ is minus the maximized log-likelihood of the data, so it measures something like (minus) fit or {\em error\/} -- in the linear regression case, it can be directly related to the mean squared error, Section~\ref{sec:regression}. Thus, (\ref{eq:refinedmdl}) embodies a trade-off between lack of fit
(measured by minus log-likelihood) and complexity (measured by
$\complexity_n(\MN{j})$). 
\commentout{Waarom heet het ook al weer 'parametric'?}
\index{stochastic complexity}
\index{model cost}
\index{model complexity}
The {\em confidence\/} in the decision is given by the codelength difference
$$
\biggl| - \log \punivnml(D \mid \MN{1}) - [  - \log \punivnml(D \mid \MN{2})] \biggr|.
$$
In general, $- \log \punivnml(D \mid\M)$ can only be evaluated
numerically -- the only exception this author is aware of is when $\M$
is the Gaussian family, Example~\ref{ex:boundary}. In many cases even
numerical evaluation is computationally problematic. But the re-interpretations of
$\punivnml$ we provide below also indicate that in many cases, $- \log
\puniv(D \mid\M)$ is relatively easy to approximate.
\begin{example}{\bf [Refined MDL and GLRT]}
\label{ex:GLRT}
\rm  \index{likelihood ratio test}
Generalized likelihood ratio testing \cite{CasellaB90}
\commentout{Doen CHECK REFERENCE}
%adopt the model
%whose best-fitting distribution best fits the data. That is, we 
tells us to 
pick the $\MN{j}$ maximizing $\log P(D \mid \hat{\theta}^{(j)}(D)) + \const$ where $\const$ is determined by the desired type-I and type-II errors. In practice one often applies a naive variation\endnote{To be fair, we should add  that 
this naive version of GLRT is introduced here for educational purposes only. It is not recommended by any serious statistician!}, simply picking the model  $\MN{j}$ maximizing $\log P(D \mid \hat{\theta}^{(j)}(D))$. 
This amounts to ignoring the complexity terms $\complexity_n(\MN{j})$
in (\ref{eq:refinedmdl}): MDL tries to avoid overfitting by picking
the model maximizing the {\em normalized\/} rather than the ordinary
likelihood. The more distributions in $\M$ that fit the data well, the
larger the normalization term. \end{example} The hope is that the
normalization term $\complexity_n(\MN{j})$ strikes the right balance
between complexity and fit.  Whether it really does depends on whether
$\complexity_n$ is a `good' measure of complexity. In the remainder of
this section we shall argue that it is, by giving four different
interpretations of $\complexity_n$ and of the resulting trade-off
(\ref{eq:refinedmdl}):
\begin{enumerate}
\item Compression interpretation.
\item Counting interpretation.
\item Bayesian interpretation.
\item Prequential (predictive) interpretation.
\end{enumerate}
\subsection{Compression Interpretation}
\label{sec:compression}
Rissanen's original goal was to select the model that detects the most
regularity in the data; he identified this with the `model that allows
for the most compression of data $x^n$'. To make this precise, 
a code is associated with each model. The NML code with lengths $- \log
\punivnml(\cdot \mid \MN{j})$ seems to be a very reasonable choice 
for such a code because of the following two properties:
\begin{enumerate}
\item The better the best-fitting distribution in $\MN{j}$ fits the
  data, the shorter the codelength $ - \log \punivnml(D \mid \MN{j})$.
\item No distribution in $\MN{j}$ is given a prior preference over any
  other distribution, since the regret of $\punivnml(\cdot \mid
  \MN{j})$ is the same for all $D \in {\cal X}^n$
  (Equation~(\ref{eq:minimaxb})). $\punivnml$ is the {\em only\/} complete
  prefix code with this property, which may be restated as:
$\punivnml$ treats {\em all distributions
    within each $\MN{j}$ on the same footing!}
\end{enumerate} 
Therefore, if one is willing to accept the basic ideas
underlying MDL as {\em first principles}, then the  use of NML in
model selection is now justified to some extent.  Below we
give additional justifications that are not directly based on
data compression; but we first provide some further interpretation of
$- \log \punivnml$.
\paragraph{Compression and Separating Structure from Noise}
We present the following ideas in an imprecise fashion --
\citeN{RissanenT04} recently showed how to make them precise.
The stochastic complexity of data $D$ relative to $\M$, given by
(\ref{eq:refinedmdl}) can be interpreted as the amount of information in the data
relative to $\M$, measured in bits. Although a one-part codelength, it
still consists of two terms: a term
$\complexity_n(\M)$ measuring the amount of {\em structure\/} or {\em
  meaningful information\/} in the data (as `seen through $\M$'), and
a term 
$- \log P(D \mid \hat{\theta}(D))$ measuring the amount of {\em noise\/}
or {\em accidental information\/} in the data. To see that this second
term measures noise, consider the regression example, 
Example~\ref{ex:modsel}, again. As
will be seen in Section~\ref{sec:beyond}, Equation (\ref{eq:gauss}), in that case $- \log P(D
\mid \hat{\theta}(D))$ becomes equal to a linear function of the
mean squared error of the best-fitting polynomial in the set of $k$-th
degree polynomials. To see that the first
term measures structure, we reinterpret it below as the number of
bits needed to specify a `distinguishable' distribution in $\M$, using
a uniform code on all `distinguishable'
distributions. 
\subsection{Counting Interpretation}
\label{sec:counting}
The parametric complexity can be interpreted as measuring (the log of)
the {\em number of distinguishable distributions in the model}. Intuitively,
the more distributions a model contains, the more patterns it can fit
well so the larger the risk of overfitting. However, if two
distributions are very `close' in the sense that they assign high
likelihood to  the same patterns, they do not contribute so
much to the complexity of the overall model. It seems that we should
measure complexity of a model in terms of the number of distributions
it contains that are `essentially different' (distinguishable), and we now show that
$\complexity_n$ measures something like this.  Consider a finite model
$\M$ with parameter set $\paraset = \{ \theta_1, \ldots, \theta_M \}$.
Note that
\begin{multline} \sum_{x^n \in \samplespace^n}
P(x^n | \hat{\theta}(x^n)) = \sum_{j= 1.. M}  
\sum_{x^n: \hat{\theta}(x^n) = \theta_j} P(x^n | \theta_j) = \\
\sum_{j= 1.. M} \bigl( 1- \sum_{x^n: \hat{\theta}(x^n) \neq  \theta_j} 
P(x^n | \theta_j) \bigr) = M -
\sum_{j} P(\hat{\theta}(x^n) \neq \theta_j | \theta_j). \nonumber
\end{multline}
We may think of $ P(\hat{\theta}(x^n) \neq \theta_j | \theta_j)$ as the probability, according to $\theta_j$, that the data look as if they come from some $\theta \neq \theta_j$. Thus, it is the probability that $\theta_j$ is {\em mistaken\/} for another distribution in $\Theta$. Therefore, for finite $\M$,
Rissanen's model complexity is the logarithm of the {\em number of distributions minus
the  summed probability that some $\theta_j$ is `mistaken' for some
$\theta \neq \theta_j$}. Now suppose $\M$ is i.i.d.
By the law of large numbers \cite{Feller68a},  we immediately see
that the `sum of mistake probabilities' $\sum_{j}
P(\hat{\theta}(x^n) \neq \theta_j | \theta_j) $ tends to 0 as $n$ grows. 
%By using the 
%Chernoff/Hoeffding
%bounds \cite{Hoeffding63} one can even show that the sum converges to 0
%exponentially fast. 
It follows that for large $n$, the model complexity
converges to $\log M$. For large $n$, the distributions in $\M$
are `perfectly distinguishable' (the probability that a sample coming from one
is more representative of another is negligible), and then the parametric
complexity $\complexity_n(\M)$ of $\M$ is simply the $\log$ of the number of
distributions in $\M$. 

\begin{example}{ \bf \ [NML vs.  Two-part Codes]}
\rm
Incidentally, this shows that for finite i.i.d. $\M$, the two-part code with uniform prior $W$ on 
$\M$ is asymptotically minimax optimal: for all $n$, the regret of the
2-part code is $\log M$ (Equation~\ref{eq:finite}), 
whereas we just showed that for
$n \rightarrow 
\infty$, $\regret(\punivnml) = \complexity_n(\M) \rightarrow \log
M$. However, for small $n$, some distributions in $\M$ may be
mistaken for one another; the number of {\em distinguishable\/} distributions
in $\M$ is then smaller than the actual number of distributions, and this is
reflected in $\complexity_n(\M)$ being (sometimes much) smaller than $\log M$.
\end{example}
For the more interesting case of parametric models, containing
infinitely many distributions, Balasubramanian
\citeyear{Balasubramanian97,Balasubramanian04} has a somewhat
different counting interpretation of $\complexity_n(\M)$ as a ratio
between two volumes. \citeN{RissanenT04} give a more direct counting
interpretation of $\complexity_n(\M)$. These extensions are both based
on the asymptotic expansion of $\punivnml$, which we now discuss.
\paragraph{Asymptotic Expansion of $\punivnml$ and $\complexity_n$}
Let $\M$ be a $k$-dimensional parametric model. Under 
regularity conditions on $\M$ and the parameterization $\Theta
\rightarrow \M$, to be detailed below, we obtain the following
asymptotic expansion of $\complexity_n$ \cite{Rissanen96,TakeuchiB97,TakeuchiB98,Takeuchi00}:
\begin{equation}
\label{eq:asnml}
\complexity_n(\M) = \frac{k}{2} \log \frac{n}{2 \pi} + \log
\int_{\theta \in \Theta} \sqrt{ |\fishermatrix(\theta )|}  d \theta + o(1).
\end{equation}
Here $k$ is the number of parameters (degrees of freedom) in model
$\M$, $n$ is the sample size, and $o(1) \rightarrow 0$ as $n
\rightarrow \infty$.  $|\fishermatrix(\theta)|$ is the determinant of
the $k \times k$ {\em Fisher information matrix\endnote{The standard
    definition of Fisher information \cite{KassV97} is in terms of
    first derivatives of the log-likelihood; for most parametric models of
    interest, the present definition coincides with the standard
    one.}} $\fishermatrix$ evaluated at $\theta$. In case $\M$ is an
  i.i.d. model, $\fishermatrix$ is given by
$$
%\begin{equationr}
%\label{eq:fisher}
\fishermatrix_{ij}(\theta^*) \isbydefinition 
\Exp_{\theta^*} \biggl\{- \frac{\partial^2}{\partial \theta_i \partial \theta_j}
 \log P(X | \theta) \biggr\}_{\theta = \theta^*} .
%\end{equationr}
$$
This is generalized to non-i.i.d. models as follows:
$$\fishermatrix_{ij}(\theta^*) \isbydefinition 
\lim_{n \rightarrow \infty} \frac{1}{n} 
\Exp_{\theta^*} \biggl\{- \frac{\partial^2}{\partial \theta_i \partial \theta_j}
 \log P(X^n | \theta) \biggr\}_{\theta = \theta^*}.$$
(\ref{eq:asnml}) only holds if the model $\M$, its parameterization
$\Theta$ 
and the sequence
$x_1, x_2, \ldots$ all satisfy certain conditions. Specifically, we
require:
\begin{enumerate}
\item $\complexity_n(\M) < \infty$ and $\int \sqrt{
    |\fishermatrix(\theta )|} d \theta < \infty$;
\item $\hat{\theta}(x^n)$ does not come
  arbitrarily close to the boundary of $\Theta$: for some $\epsilon >
  0$, for all large $n$, $\hat{\theta}(x^n)$ remains farther than
  $\epsilon$ from the boundary of $\Theta$.
\item $\M$ and $\Theta$ satisfy certain further conditions. A simple
  sufficient condition is that $\M$ be an exponential family
  \cite{CasellaB90}. Roughly, 
this is a family that can be parameterized so that for all $x$, 
$P(x \mid \beta) = \exp(\beta t(x)) f(x) g(\beta)$, where $t: {\cal
  X} \rightarrow \reals$ is a function of $X$.  The Bernoulli model is
an exponential family, as can be seen by setting $\beta := \ln ( 1-
\theta) - \ln \theta$ and $t(x) = x$. Also the 
multinomial, Gaussian, Poisson, Gamma, exponential, Zipf and many
other models are exponential families; but, for example, mixture models are not.
\end{enumerate}
More general conditions are given by
\citeANP{TakeuchiB97} \citeyear{TakeuchiB97,TakeuchiB98,Takeuchi00}. Essentially, if $\M$ behaves `asymptotically' like an exponential
family, then (\ref{eq:asnml}) 
still holds. For example, (\ref{eq:asnml}) holds for the Markov models
and for AR and ARMA processes.
\begin{example}{\bf [Complexity of the Bernoulli Model]}
\rm 
The Bernoulli model $\BernoulliN{0}$ can be parameterized in a
  1-1 way by the unit interval (Example~\ref{ex:markov}). Thus, $k=1$. 
An easy calculation shows that the Fisher information is given by 
$\theta (1- \theta)$.
Plugging this into
  (\ref{eq:asnml}) and calculating $\int \sqrt{
    |\theta(1- \theta) | } d \theta$ gives 
$$\complexity_n(\BernoulliN{0}) = \frac{1}{2} \log n 
+ \frac{1}{2} \log \frac{\pi}{2} -3 + o(1)
=\frac{1}{2} \log n -2.674251935 + o(1).$$
Computing the integral of the Fisher determinant is not
easy in general. 
\citeN{HansonF04} compute it for several practically relevant models.
\end{example}
Whereas for finite $\M$, $\complexity_n(\M)$ remains finite, for
parametric models it generally grows logarithmically in $n$. Since
typically  $- \log P(x^n \mid \hat{\theta}(x^n))$
grows linearly in $n$, it is still the case that for {\em fixed\/}
dimensionality $k$ (i.e. for a fixed $\M$ that is $k$-dimensional)
and large $n$, the part of the codelength $- \log \punivnml(x^n \mid
\M)$ due to the complexity of $\M$ is very small
compared to the part needed to encode data $x^n$ with
$\hat{\theta}(x^n)$. The term $\curvelength{\paraset}$ may be
interpreted as the contribution of the {\em functional form\/} of $\M$
to the  model complexity \cite{Balasubramanian04}. It does not grow
with $n$ so that, when selecting between two models, it becomes irrelevant and can be ignored for {\em very\/} large $n$. But for
small $n$, it can be important, as can be seen from Example~\ref{ex:modselb}, 
Fechner's and Stevens' model. Both models have two parameters, yet the
$\curvelength{\paraset}$-term is much larger for Fechner's than for
Stevens' model. In the experiments of \citeN{MyungBP00}, the parameter
set was restricted to $0 < a < \infty, 0 < b < 3$ for Stevens' model
and $0 < a < \infty, 0 < b < \infty$ for Fechner's model. The variance
of the error $Z$ was set to $1$ in both models. With these values, the
difference in $\curvelength{\paraset}$ is 3.804, which is
non-negligible for small samples. Thus, Stevens' model contains more
distinguishable distributions than Fechner's, and is better able to
capture random noise in the data -- as \citeN{Townsend75} already
speculated almost 30 years ago. Experiments suggest that for
regression models such as Stevens' and Fechner's', as well as for
Markov models and general exponential families, the approximation
(\ref{eq:asnml}) is reasonably accurate already for small samples. But
this is certainly not true for general models:
\ownbox{The Asymptotic Expansion of
  $\complexity_n$ Should Be Used with Care!}{
(\ref{eq:asnml}) does {\em
    not\/} hold for all parametric models; and for some models for
  which it does hold, the $o(1)$ term may only converge to $0$ only
  for quite large sample sizes. Foster and Stine
\citeyear{FosterS99,FosterS04} show that 
the approximation (\ref{eq:asnml}) is, in general,
 only valid if $k$   is much smaller than $n$. 
%Navarro exhibits two models
%$\MN{1}$ and $\MN{2}$ so that $\MN{1}$ is
%nested in $\MN{2}$, and hence automatically $\complexity_n(\MN{2}) >
%\complexity_n(\MN{1})$. Yet for small $n$, (\ref{eq:asnml}) is smaller
%for $\MN{1}$ than for $\MN{2}$!
} 
\paragraph{Two-part codes and $\complexity_n(\M)$}
We now have a clear guiding principle (minimax regret) which we can
use to construct `optimal' two-part codes, that achieve the minimax
regret among all {\em two-part codes}. How do such optimal two-part
codes compare to the NML codelength? Let $\M$ be a $k$-dimensional
model. By slightly adjusting the arguments of
\cite[Appendix]{BarronC91}, one can show that, under regularity
conditions, the minimax optimal two-part code $\punivtp$ achieves
regret
$$
- \log \punivtp(x^n \mid \M ) + \log P(x^n \mid \hat{\theta}(x^n))
= \frac{k}{2} \log \frac{n}{2 \pi} + \log
\int_{\theta \in \Theta} \sqrt{ |\fishermatrix(\theta )|}  d \theta
+ f(k) + o(1),
$$
where $f: {\mathbb N} \rightarrow \reals$ 
is a bounded positive function satisfying 
$\lim_{k \rightarrow \infty} f(k) = 0$. Thus, for large $k$, {\em
  optimally designed\/} two-part codes are about as good as NML. The
problem with two-part code MDL is that {\em in practice}, people often
use much cruder codes with much larger minimax regret. 
\subsection{Bayesian Interpretation}
\label{sec:bayes1}
The Bayesian method of statistical inference provides several alternative
approaches to model selection. The most popular of these is based on {\em Bayes factors} \cite{KassR95}.
The Bayes factor method is very closely related
to the refined MDL approach. 
\commentout{
Suppose we want to select between models
$\MN{1}$ and $\MN{2}$. According to the
Bayesian approach, we should specify a prior distribution on $\MN{1}$
and $\MN{2}$, and then, conditioned on each model, prior densities
$w_1$ and $w_2$. As in (\ref{eq:bayesuniversal}), 
the Bayesian marginal likelihood is defined as
$$
%\begin{equationr}
%\label{eq:bayesint}
\punivbayes(x^n \mid \MN{j}) = \int P(x^n \mid \theta) w_j(\theta) d\theta.
%\end{equation} 
$$
Bayes then tells us to select the model with the maximum a posteriori
probability $P(\MN{j} \mid D)$. This can be computed from Bayes' rule as
$$ 
P(\MN{j} \mid D) = 
\frac{\punivbayes(D \mid \MN{j}) W(j)}{\sum_{k = 1,2} P(D \mid \MN{k}) W(k).
}  
$$
Since the denominator does not depend on $j$, Bayes tells us to
pick the $\MN{j}$ maximizing $\punivbayes(D \mid \MN{j}) W(j)$.
Typically, the contribution of the prior $W(j)$ is negligible for all
but the smallest $n$, so that the model with largest $\punivbayes(D
\mid \MN{j})$ is chosen.  } Assuming uniform priors on models $\MN{1}$
and $\MN{2}$, it tells us to select the model with largest marginal
likelihood $\punivbayes(x^n \mid \MN{j})$, where $\punivbayes$ is as
in (\ref{eq:bayes}), with the sum replaced by an integral, and
$w^{(j)}$ is the density of the prior distribution on $\MN{j}$:
\begin{equation}
\label{eq:bayesint}
\punivbayes(x^n \mid \MN{j}) = \int P(x^n \mid \theta) w^{(j)}(\theta) d\theta.
\end{equation} 
\paragraph{$\M$ is Exponential Family}
Let now $\punivbayes = \punivbayes(\cdot \mid \M)$ for some fixed
model $\M$. 
Under regularity conditions on $\M$, we can perform a {\em Laplace
approximation\/} of the integral in (\ref{eq:bayes}). For the special
case that $\M$ is an exponential family, we obtain the following
expression for the regret
\cite{Jeffreys61,Schwarz78,KassR95,Balasubramanian97}:
\commentout{I couldn't find it in Jeffreys 61, but what can you do!
- Stone 79 gives a reference but I cant read Jeffreys' math}
\begin{equation}
\label{eq:asbayes1}
- \log \punivbayes(x^n ) - [ - \log P(x^n \mid \hat{\theta}(x^n))] = 
\frac{k}{2} \log \frac{n}{2 \pi} - \log w(\hat{\theta}) + \log 
\sqrt{{| \fishermatrix(\hat{\theta}) |}}  + o(1).
\end{equation}
Let us compare this with (\ref{eq:asnml}). 
Under the regularity conditions needed for (\ref{eq:asnml}), the
quantity on the right of (\ref{eq:asbayes1}) is within $O(1)$ of
$\complexity_n(\M)$. Thus, the code length achieved with $\punivbayes$ is
within a constant of the minimax optimal $- \log \punivnml(x^n)$.
Since $ - \log P(x^n \mid \hat{\theta}(x^n))$ increases linearly in
$n$, this means that if we compare two models $\MN{1}$ and $\MN{2}$,
then for large enough $n$, Bayes and refined MDL
select the same model. If we equip the Bayesian universal model with a
special prior known as the {\em Jeffreys-Bernardo\/} prior
\cite{Jeffreys46,BernardoS94},
\begin{equation}
\label{eq:jeffreys}
w_{\text{Jeffreys}}(\theta) =
\frac{\sqrt{|\fishermatrix(\theta)|}}{\int_{\theta \in \Theta} \sqrt{
    |\fishermatrix(\theta )|} d \theta},
\end{equation}
then Bayes and refined NML become even more closely related: plugging
in (\ref{eq:jeffreys}) into (\ref{eq:asbayes1}), we find that the
right-hand side of (\ref{eq:asbayes1}) now simply {\em coincides\/}
with (\ref{eq:asnml}). A concrete example of Jeffreys' prior is given
in Example~\ref{ex:laplace}. Jeffreys introduced his prior as a `least
informative prior', to be used when no useful prior knowledge about
the parameters is available \cite{Jeffreys46}. As one may expect from
such a prior, it is invariant under continuous 1-to-1
reparameterizations of the parameter space. The present analysis shows
that, when $\M$ is an exponential family, then it also leads to
asymptotically minimax codelength regret: for large $n$, {\em refined
  NML model selection becomes indistinguishable from
Bayes factor model selection with
  Jeffreys' prior}.

\paragraph{$\M$ is not an Exponential Family}
Under weak conditions on
$\M$, $\Theta$ and the sequence $x^n$, we get the following generalization of (\ref{eq:asbayes1}):
\begin{multline}
\label{eq:asbayes2}
- \log \punivbayes(x^n \mid \M) = \\
- \log P(x^n \mid \hat{\theta}(x^n)) + \frac{k}{2} \log \frac{n}{2 \pi} - \log w(\hat{\theta}) 
+  \log  \sqrt{\bigl| \observedmatrix(x^n) \bigr| }  + o(1). 
\end{multline}
Here $\observedmatrix(x^n)$ is the so-called {\em observed
  information}, 
sometimes also called {\em observed Fisher
  information\/}; see \cite{KassV97} for a definition.
\commentout{
\endnote{\rm To be fully precise, what we call the
    `observed information (matrix)' is $n^{-1}$ times the standard
    definition of `observed information (matrix)' \cite{KassV97}.}, defined as
\begin{equation}
\label{eq:empiricalfisher}
\observedmatrix(x^n) \isbydefinition  \biggl\{ \frac{1}{n} \frac{d^2}{d \theta^2} -  \log
P(x^n | \theta) \biggr\}_{\theta = \hat{\theta}}.
\end{equation}
}
If $\M$ is an exponential family, then 
the observed Fisher information at $x^n$ coincides with the Fisher
information at $\hat{\theta}(x^n)$, leading to (\ref{eq:asbayes1}). If
$\M$ is not exponential, then if data are distributed according
  to one of the distributions in $\M$, the observed Fisher information
  still converges with probability 1 to the expected Fisher
  information. If $\M$ is
  neither exponential, nor are the data actually generated by a
  distribution in $\M$, then there may be $O(1)$-discrepancies between $-
  \log \punivnml$ and $- \log \punivbayes$ even for large $n$.
\subsection{Prequential Interpretation}
\label{sec:prequential}
\paragraph{Distributions as Prediction Strategies}
Let $P$ be a distribution on ${\cal X}^n$. Applying the definition of
conditional probability, we can write for every $x^n$:
\begin{equation}
\label{eq:chainrule}
P(x^n) = \prod_{i=1}^n \frac{P(x^i)}{P(x^{i-1})} = \prod_{i=1}^n P(x_i \mid x^{i-1}),
\end{equation}
so that also
\begin{equation}
\label{eq:aclogloss}
- \log P(x^n) = \sum_{i=1}^n - \log P(x_i \mid x^{i-1})
\end{equation}
Let us abbreviate $P(X_i = \cdot \mid X^{i-1} = x^{i-1})$ to $P(X_i
\mid x^{i-1})$. $P(X_i \mid x^{i-1})$ (capital $X_i$) is the {\em
  distribution\/} (not a single number) of $X_i$ given $x^{i-1}$;
$P(x_i \mid x^{i-1})$ (lower case $x_i$) is the {\em probability\/} (a
single number) of actual outcome $x_i$ given $x^{i-1}$. We can think
of $- \log P(x_i \mid x^{i-1})$ as the {\em loss\/} incurred when
predicting $X_i$ based on the conditional distribution $P(X_i \mid
x^{i-1})$, and the actual outcome turned out to be $x_i$. Here `loss'
is measured using the so-called {\em logarithmic score}, also known
simply as `log loss'. Note that the more likely $x$ is judged to be,
the smaller the loss incurred when $x$ actually obtains. The log loss
has a natural interpretation in terms of sequential {\em gambling\/}
\cite{CoverT91}, but its main interpretation is still in terms of
coding: by (\ref{eq:aclogloss}), the codelength needed to encode $x^n$
based on distribution $P$ is just the {\em accumulated log loss
  incurred when $P$ is used to sequentially predict the $i$-th outcome
  based on the past $(i- 1)$-st outcomes}.

(\ref{eq:chainrule}) gives a fundamental
re-interpretation of probability distributions as prediction
strategies, \index{prediction} \index{loss function} mapping each
individual sequence of {\em past observations\/} $x_1, \ldots,
x_{i-1}$ to a { \em probabilistic prediction of the next outcome\/}
$P(X_i \mid x^{i-1})$.  Conversely, (\ref{eq:chainrule}) also shows
that every probabilistic prediction strategy for sequential prediction
of $n$ outcomes may be thought of as a probability distribution on
${\cal X}^n$: a strategy is identified with a function mapping all
potential initial segments $x^{i-1}$ to the prediction that is made for
the next outcome $X_i$, after having seen $x^{i-1}$. Thus, it is a
function $S: \cup_{0 \leq i < n} {\cal X}^i \rightarrow {\cal P}_{\cal
  X}$, where ${\cal P}_{\cal X}$ is the set of distributions on ${\cal
  X}$. We can now define, for each $i < n$, all $x^i \in {\cal X}^i$,
$P(X_i \mid x^{i-1}) \isbydefinition S(x^{i-1})$. We can turn these
partial distributions into a distribution on ${\cal X}^n$ by
sequentially plugging them into (\ref{eq:chainrule}). 

\paragraph{Log Loss for Universal Models}
Let $\M$ be some parametric model and let $\puniv$ be some universal
model/code relative to $\M$. 
\commentout{Doen NIET NODIG DENK IK WANT IK ZEG HET STRAKS AL IN GLRT EXAMPLE 
We can think of $- \log P(x^n \mid
\hat{\theta}(x^n))$ as the minimum accumulated log loss, achieveable
with any of the distributions in $\M$, with hindsight, after seeing
data $x^n$. $-\log \puniv(x^n)$ is the accumulated log loss
achieved with the universal model. The prediction strategy $\puniv$ can be applied without hindsight - it does not require knowledge of the ML estimator 
$\hat{\theta}(x^n)$. }
What do the individual predictions $\puniv(X_i \mid x^{i-1})$ look
like? Readers familiar with Bayesian statistics will realize that for
i.i.d. models, the Bayesian predictive distribution $\punivbayes(X_i
\mid x^{i-1})$ converges to the ML distribution $P(\cdot \mid
\hat{\theta}(x^{i-1}))$; Example~\ref{ex:laplace} provides a concrete
case.  It seems reasonable to assume that something similar holds
not just for $\punivbayes$ but for universal models in general.
This in turn suggests that we may approximate the conditional
distributions $\puniv(X_i \mid x^{i-1})$ of any `good' universal model
by the maximum likelihood predictions $P(\cdot \mid
\hat{\theta}(x^{i-1}))$. Indeed, we can recursively define the {\em `maximum likelihood plug-in'\/} distribution 
$\punivplugin$ by setting, for $i = 1$ to $n$,
\begin{equation}
\label{eq:preqdef}
\punivplugin(X_i = \cdot \mid x^{i-1}) \isbydefinition P(X = \cdot \mid
\hat{\theta}(x^{i-1})).
\end{equation}
Then
\begin{equation}
\label{eq:preqsnack}
 - \log \punivplugin(x^n) \isbydefinition 
\sum_{i=1}^n - \log P(x_i \mid \hat{\theta}(x^{i-1})).
\end{equation}
Indeed, it turns out that under regularity conditions on $\M$ and $x^n$, 
\begin{equation}
\label{eq:preq}
- \log \punivplugin(x^n) = - \log P(x^n \mid
\hat{\theta}(x^n)) + \frac{k}{2} \log n + O(1),
\end{equation}
showing that $\punivplugin$ acts as a universal model relative to
$\M$, its performance being within a constant of the minimax optimal
$\punivnml$.  The construction of $\punivplugin$ can be easily
extended to non-i.i.d.  models, and then, under regularity conditions,
(\ref{eq:preq}) still holds; we omit the details.
\begin{ownquote}
  We note that all general proofs of (\ref{eq:preq}) that we are aware
  of show that (\ref{eq:preq}) holds with probability 1 or in
  expectation for sequences generated by some distribution in $\M$
  \cite{Rissanen84,Rissanen86a,Rissanen89}.  Note that the expressions
  (\ref{eq:asnml}) and (\ref{eq:asbayes2}) for the regret of
  $\punivnml$ and $\punivbayes$ hold for a much wider class of
  sequences; they also hold with probability 1 for i.i.d. sequences
  generated by sufficiently regular distributions {\em outside\/}
  $\M$. Not much is known about the regret obtained by $\punivplugin$
  for such sequences, except for some special cases such as if $\M$ is the
  Gaussian model.
\end{ownquote}
In general, there is no need to use the ML estimator
$\hat{\theta}(x^{i-1})$ in the definition (\ref{eq:preqdef}). Instead,
we may try some other estimator which asymptotically converges to the
ML estimator -- it turns out that some estimators considerably
outperform the ML estimator in the sense that (\ref{eq:preqsnack})
becomes a much better approximation of $- \log \punivnml$,  
see Example~\ref{ex:laplace}. Irrespective of
whether we use the ML estimator or something else, we call model
selection based on (\ref{eq:preqsnack}) the
{\em prequential\/} form of MDL in honor of A.P. Dawid's `prequential
analysis', Section~\ref{sec:others}.
It is also known as `predictive MDL'. The validity of
(\ref{eq:preq}) was discovered independently by Rissanen
\citeyear{Rissanen84} and Dawid \citeyear{Dawid84}.

The prequential view gives us a fourth interpretation of refined MDL
model selection: given models $\MN{1}$ and $\MN{2}$, MDL tells us to
pick the model that minimizes the accumulated prediction error
resulting from sequentially predicting future outcomes given all the
past outcomes.
\begin{example}{\bf [GLRT and Prequential Model Selection]}
\rm
How does this differ from the naive version of the 
generalized likelihood ratio test (GLRT)
that we introduced in Example~\ref{ex:GLRT}? In GLRT, we
associate with each model the log-likelihood (minus log loss) that can be
obtained by the ML estimator. This is the predictor within the model
that minimizes log loss {\em with hindsight}, {\em after\/} having
seen the data.
%; if the model $\M$ is large enough, the fact that such a
%ML estimator achieves small log loss tells us nothing about its
%quality on future data, whence GLRT is prone to overfitting.  
In contrast, prequential model selection associates with each model
the log-likelihood (minus log loss) that can be obtained by using a
sequence of ML estimators $\hat{\theta}(x^{i-1})$ to predict data
$x_{i}$. Crucially, the data on which ML estimators are evaluated
has not been used in constructing the ML estimators themselves. This
makes the prediction scheme `honest' (different data are used for 
training and testing) and explains why it automatically protects us against overfitting.
\end{example}
\begin{example}{\bf [Laplace and Jeffreys]}
\label{ex:laplace}
\rm Consider the
prequential distribution for the Bernoulli model,
Example~\ref{ex:markov}, defined as in
(\ref{eq:preqdef}). We show that if we take $\hat{\theta}$ in
(\ref{eq:preqdef}) equal to the ML estimator $n_{[1]} / n$, then the
resulting $\punivplugin$ is not a universal model; but a slight
modification of the ML estimator makes $\punivplugin$ a very good
universal model.  Suppose that $n \geq 3 $ and $(x_1,x_2,x_3) =
(0,0,1)$ -- a not-so-unlikely initial segment according to most $\theta$.
Then $\punivplugin(X_3 = 1 \mid x_1,x_2) = P(X = 1 \mid
\hat{\theta}(x_1,x_2)) = 0$, so that by (\ref{eq:preqsnack}), 
$$- \log \punivplugin(x^n) \geq -
\log \punivplugin(x_3 \mid x_1,x_2) = \infty,$$
whence $\punivplugin$
is not universal.  Now let us consider the modified ML estimator
\begin{equation}
\label{eq:laplace}
\hat{\theta}_{\lambda}(x^n) \isbydefinition \frac{n_{[1]} +
  \lambda}{n + 2 \lambda}.
\end{equation}
If we take $\lambda = 0$, we get the ordinary ML estimator. If we take
$\lambda = 1$, then an exercise involving beta-integrals shows that,
for all $i, x^i$, $P(X_i \mid \hat{\theta}_1(x^{i-1})) =
\punivbayes(X_i \mid x^{i-1})$, where $\punivbayes$ is defined
relative to the uniform prior $w(\theta) \equiv 1$.  Thus
$\hat{\theta}_1(x^{i-1})$ corresponds to the Bayesian predictive
distribution for the uniform prior. This prediction rule was advocated
by the great probabilist P.S. de Laplace, co-originator of Bayesian
statistics. It may be interpreted as ML estimation based on an {\em
  extended\/} sample, containing some `virtual' data: an extra $0$ and
an extra $1$.

Even better, a similar calculation shows that if we take $\lambda =
1/2$, the resulting estimator is equal to $\punivbayes(X_i \mid
x^{i-1})$ defined relative to {\em Jeffreys' prior}. Asymptotically,
$\punivbayes$ with Jeffreys' prior achieves the same codelengths as
$\punivnml$ (Section~\ref{sec:bayes1}). It follows that $\punivplugin$
with the slightly modified ML estimator is asymptotically
indistinguishable from the optimal universal model $\punivnml$!
\commentout{ As an aside, the prequential interpretation gives us a
  simple way to efficiently compute $- \log \punivbayes(x^n)$ for the
  Bernoulli model: we simply add, for $i = 1, \ldots, n$, $- \log
  \punivbayes(x_i \mid x^{i-1})$, which can be easily computed by
  (\ref{eq:laplace}). By (\ref{eq:aclogloss}), the result must be
  equal to $ - \log \punivbayes(x^n)$.  \commentout{Doen In boek,
    extra section on 'how to compute it'}

The example can be extended to the multinomial model, the set of distributions on ${\cal X} = \{1, \ldots, m \}$ extended to $n$ outcomes by independence.
In that case, (\ref{eq:laplace}) becomes
\begin{equation}
\label{eq:laplaceb}
\hat{\theta}_{\lambda}(x^n) \isbydefinition (\frac{n_{[1]} +
  \lambda}{n + m \lambda}, \ldots,\frac{n_{[m]} +
  \lambda}{n + m \lambda}),
\end{equation}
where $\lambda = 1$ corresponds to the uniform, and $\lambda = 1/2$ to the Jeffreys' prior.
}

For more general models $\M$, such simple modifications of the ML
estimator usually do not correspond to a Bayesian predictive
distribution; for example, if $\M$ is not convex (closed under taking mixtures) then a point estimator (an element of $\M$) typically does not correspond to the Bayesian predictive distribution (a mixture of elements of $\M$). 
Nevertheless, modifying the ML estimator by adding some
virtual data $y_1, \ldots, y_m$ and replacing $P(X_i \mid
\hat{\theta}(x^{i-1}))$ by $P(X_i \mid \hat{\theta}(x^{i-1}, y^m))$ in
the definition (\ref{eq:preqdef}) may still lead to good universal
models. This is of great practical importance, since, using
(\ref{eq:preqsnack}), $- \log \punivplugin(x^n)$ is often much easier to
compute than $- \log \punivbayes(x^n)$.  \index{Laplace, P.S. de}
\index{Jeffreys, H.}
\end{example}
\paragraph{Summary}
We introduced the refined MDL Principle for model selection in a
restricted setting. Refined MDL
amounts to selecting the model under which the data achieve the
smallest {\em stochastic complexity}, which is the codelength
according to the minimax optimal universal model. We gave an
asymptotic expansion of stochastic and parametric complexity, and
interpreted these concepts in four different ways.
\section{General Refined MDL: Gluing it All Together}
\label{sec:unify}
\index{glue}
In the previous section we introduced a `refined' MDL principle based
on minimax regret. Unfortunately, this principle can be applied only
in very restricted settings. We now show how to extend refined MDL,
leading to a general MDL Principle, applicable to a wide variety of
model selection problems. In doing so we glue all our previous
insights (including `crude MDL') together, thereby uncovering a single
general, underlying principle, formulated in Figure~\ref{fig:unified}. Therefore, {\em if one understands the material in this section, then
  one understands the Minimum Description Length Principle.}

First, Section~\ref{sec:uncountable}, we show how to compare
infinitely many models. Then, Section~\ref{sec:boundary}, we show how
to proceed for models $\M$ for which the parametric complexity is
undefined. Remarkably, a single, general idea resides behind our solution of both problems, and this leads us to formulate, in Section~\ref{sec:general}, a single, general refined MDL Principle.
\subsection{Model Selection with Infinitely Many Models}
\label{sec:uncountable}
Suppose we want to compare more than two models for the same data.  If
the number to be compared is finite, we can proceed as before and pick
the model $\MN{k}$ with smallest $- \log \punivnml(x^n \mid \MN{k})$.
If the number of models is infinite, we have to be more careful. Say
we compare models $\MN{1},\MN{2}, \ldots$ for data $x^n$. We may be
tempted to pick the model minimizing $- \log \punivnml(x^n \mid
\MN{k})$ over all $k \in \{1,2, \ldots \}$, but in some cases this
gives unintended results. To illustrate, consider the extreme case
that every $\MN{k}$ contains just one distribution. For example, let
$\MN{1} = \{P_1\}, \MN{2} = \{P_2\}, \ldots$ where $\{P_1, P_2, \ldots
\}$ is the set of {\em all\/} Markov chains with rational-valued
parameters. In that case, $\complexity_n(\MN{k}) = 0$ for all $k$, and 
we would always select the maximum likelihood Markov chain that
assigns probability 1 to data $x^n$. Typically this will be a chain of
very high order, severely overfitting the data. This cannot be right!
A better idea is to pick the model minimizing
\begin{equation}
\label{eq:twopartnml}
- \log \punivnml(x^n \mid \MN{k}) + L(k), 
\end{equation}
where $L$ is the codelength function of some code for encoding model
indices $k$. We would typically choose the standard prior for the
integers, $L(k) = 2 \log k+1$, Example~\ref{ex:integers}.  By using
(\ref{eq:twopartnml}) we avoid the overfitting problem mentioned
above: if $\MN{1} = \{P_1\}, \MN{2} = \{P_2\}, \ldots$ where $P_1,P_2,
\ldots$ is a list of all the rational-parameter Markov chains,
(\ref{eq:twopartnml}) would reduce to two-part code MDL
(Section~\ref{sec:crude}) which is asymptotically consistent. On the
other hand, if $\MN{k}$ represents the set of $k$-th order Markov chains,
the term $L(k)$ is typically negligible compared to
$\complexity_n(\MN{k})$, the complexity term associated with $\MN{k}$
that is hidden in $- \log \punivnml(\MN{k})$: thus, the complexity of
$\MN{k}$ comes from the fact that for large $k$, $\MN{k}$ contains
many distinguishable distributions; not from the much smaller term
$L(k) \approx 2 \log k$.

To make our previous approach for a finite set of models compatible with
(\ref{eq:twopartnml}), we can reinterpret it as follows: we assign
uniform codelengths (a uniform prior) to the $\MN{1}, \ldots, \MN{M}$
under consideration, so that for $k = 1, \ldots, M$, $L(k) = \log M$.
We then pick the model minimizing (\ref{eq:twopartnml}). Since $L(k)$
is constant over $k$, it plays no role in the minimization and can be
dropped from the equation, so that our procedure reduces to our
original refined MDL model selection method. We shall henceforth
assume that we {\em always\/} encode the model index, either
implicitly (if the number of models is finite) or explicitly. The
general principle behind this is explained in
Section~\ref{sec:general}. 
\commentout{Acknowledge in book that MDL does not deal
  well with subjectivity. Say {\bf I} have no problems with it, as
  long as it is only used INSIDE the logarithm!}
\subsection{The Infinity Problem}
\label{sec:boundary}
For some of the most commonly used models, the parametric complexity 
$\complexity(\M)$ is undefined. A prime example is the Gaussian location model, which we discuss below. As we will see, we can `repair' the situation using the same general idea as in the previous subsection. 
\begin{example}{\bf \ Parametric Complexity of the Normal Distributions}
\label{ex:boundary}
\rm 
\index{Normal distribution}
\index{Gaussian distribution}
Let $\M$ be the family of Normal distributions with fixed variance $\sigma^2$ 
and varying mean $\mu$, identified by their densities
$$
P(x| \mu) = \frac{1}{\sqrt{2\pi} \sigma} e^{-\frac{(x-\mu)^2}{2 \sigma^2}},
$$
extended to sequences $x_1, \ldots, x_n$ by taking product
densities. As is well-known  \cite{CasellaB90}, the ML estimator
$\hat{\mu}(x^n)$ is equal to the sample mean: $\hat{\mu}(x^n) =  n^{-1} \sum_{i=1}^n x_i$. An easy calculation shows that 
$$\complexity_n(\M) = 
\int_{x^n} P(x^n \mid \hat{\mu}(x^n)) dx^n = \infty,$$
where we abbreviated $d x_1 \ldots dx_n$ to $d x^n$.
\commentout{As is well-known \cite{CasellaB90}, for every sample $x^n$, the
ML estimator $\hat{\theta}(x^n)= (\hat{\mu}(x^n),\hat{\sigma}(x^n))$ is
equal to the empirical mean and variance, respectively:
\begin{equation}
\hat{\mu}(x^n) =  \frac{1}{n} \sum_{i=1}^n x_i \ \ \ ; \ \ \  \hat{\theta}(x^n) = \frac{1}{n}
\sum_{i=1}^n (x_i - \hat{\mu}(x^n))^2 =  \frac{1}{n}
\sum_{i=1}^n (x_i)^2 - (\frac{1}{n}
\sum_{i=1}^n (x_i))^2 
\end{equation}
}
%For the normal family, 
%\begin{equation}
%\label{eq-normalcomp}
%\complexity_n(\M) = \int_{x^n \in \samplespace^n} f(x^n |
%\hat{\mu}(x^n), \hat{\theta}(x^n)) dx_1 \ldots dx_n 
%= \infty
%\end{equation}  as an easy calculation shows. Therefore, 
Therefore, we cannot use basic MDL model selection. It also turns out
that 
$I(\mu) = \sigma^{-2}$ so that
$$\curvelength{\paraset} = \int_{\mu \in \reals} \sqrt{|I(\mu)|} d
\mu = \infty.$$ 
Thus, the Bayesian universal model approach with
Jeffreys' prior cannot be applied either. Does this mean that our MDL
model selection and complexity definitions break down even in such a
simple case? Luckily, it turns out that they can be repaired, as we now show.
\citeN{BarronRY98} and \citeN{FosterS01} show that, for all intervals $[a,b]$,
\begin{equation}
\label{eq:condcomp}
\int_{x^n: \hat{\mu}(x^n) \in [a,b]} P(x^n \mid \hat{\mu}(x^n)) d x^n =
\frac{b-a}{\sqrt{2\pi} \sigma} \cdot \sqrt{n}.
\end{equation}
Suppose for the moment that it is known that 
$\hat{\mu}$ lies in some set $[-K, K]$ for some fixed $K$. Let $\M_K$
be the set of conditional distributions thus obtained: $\M_K = \{
P'(\cdot \mid \mu) \mid \mu \in \reals \}$, where $P'(x^n \mid \mu)$
is the density of $x^n$ according to the normal distribution with mean
$\mu$, conditioned on $| n^{-1} \sum x_i| \leq K$.
By (\ref{eq:condcomp}), the 
`conditional' minimax regret distribution $\punivnml(\cdot \mid \M_K)$
is well-defined for all $K> 0$. That is,  
for all $x^n$ with $| \hat{\mu}(x^n) | \leq K$,
$$ 
\punivnml(x^n \mid \M_K) = \frac{P'(x^n \mid \hat{\mu}(x^n))}{
\int_{x^n \; : \; |\hat{\mu}(x^n)| < K}  
P'(x^n \mid \hat{\mu}(x^n)) d x^n, 
}
$$
with regret (or `conditional' complexity),
$$
\complexity_n(\M_K) =
\log \int_{|\hat{\mu}(x^n)| < K}  P'(x^n \mid \hat{\mu}(x^n)) d x^n =
\log K + \frac{1}{2} \log \frac{n}{2 \pi} - \log \sigma  +1.
$$
This suggests to redefine the complexity of the full model $\M$ so
that its regret depends on the area in which $\hat{\mu}$ falls.  The
most straightforward way of achieving this is to define a {\em
  meta-universal model\/} for $\M$, combining the NML with a two-part
code: we encode data by first encoding some value for $K$.
% using a code with lengths $L(K)$ where $L$ represents some code for $K$. 
We then encode the actual data $x^n$ using the code
$\punivnml(\cdot |\M_{K})$. 
The resulting code $\punivmeta$ is a universal code for $\M$ with lengths
\begin{equation}
\label{eq:firstsol}
- \log \punivmeta(x^n |\M) \isbydefinition 
\min_{K} \ \bigl\{ 
- \log \punivmeta(x^n \mid \M_K) +L(K) \bigr\}.
\end{equation}
The idea is now to base MDL model selection on $\punivmeta(\cdot |\M)$ as
in (\ref{eq:firstsol}) rather than on the (undefined) $\punivnml(\cdot
| \M)$.  To make this work, we need to choose $L$ in a clever manner.
%Note that, had we known the bound $K$ beforehand, we would have
%%incurred a regret of $\log \complexity_n(\M_K)$ bits. The idea is to
%choose a code with lengths $L$ so that the relative additional regret
%we incur on top of this is very small, {\em no matter what $K$ is}.  
A good choice is to encode $K' = \log K $ as an integer, 
using the standard code for the integers. To see why, note that the regret of $\punivmeta$ now becomes:
\begin{multline}
\label{eq:firstsolb}
- \log \punivmeta(x^n \mid \M) - [ - \log P(x^n \mid \hat{\mu}(x^n))]
=
\\
\min_{K: \log K \in \{1,2, \ldots \} } \ \bigl\{ \log K + \frac{1}{2}
\log \frac{n}{2 \pi} - \log \sigma +1 +
2 \log \lceil \log K \rceil \ \bigr\}+1 \leq \\
\log | \hat{\mu}(x^n)| +   2 \log \log | \hat{\mu}(x^n) | + \frac{1}{2} \log \frac{n}{2 \pi} - \log
\sigma + 4 \leq \\
\complexity_n(\M_{|\hat{\mu}|}) + 2 \log \complexity_n(\M_{|\hat{\mu}|}) + 3.
\end{multline}
If we had known a good bound $K$ on $| \hat{\mu} |$ {\em a priori}, we
could have used the NML model $\punivnml(\cdot \mid \M_K)$. With
`maximal' a priori knowledge, we would have used the model
$\punivnml(\cdot \mid \M_{|\hat{\mu}|})$, leading to regret
$\complexity_n(\M_{|\hat{\mu}|})$. The regret achieved by $\punivmeta$ is {\em
  almost\/} as good as this `smallest possible regret-with-hindsight'
$\complexity_n(\M_{|\hat{\mu}|})$: the difference is much smaller
than, in fact logarithmic in, $\complexity_n(\M_{|\hat{\mu}|})$
itself, {\em no matter what $x^n$ we observe}. 
This is the underlying reason why we choose to encode $K$ with log-precision: the basic idea in refined MDL was to minimize worst-case
regret, or {\em additional code-length\/} compared to the code that achieves
the minimal code-length with hindsight. Here, we use this basic idea
on a meta-level: we design a code such that the {\em additional 
  regret\/} is minimized, compared to the code that achieves the
minimal regret with hindsight.
\end{example}
\commentout{Give rissanen's neat trick for exponential families in the
  book at this place! It is ok to integrate over \hat{\theta} rather
  than the full data!}  This meta-two-part coding idea was introduced
by \citeN{Rissanen96}. It can be extended to a wide range of models
with $\complexity_n(\M) = \infty$; for example, if the $X_i$ represent
outcomes of a Poisson or geometric distribution, one can encode a
bound on $\mu$ just like in Example~\ref{ex:boundary}. If $\M$ is the
full Gaussian model with both $\mu$ and $\sigma^2$ allowed to vary,
one has to encode a bound on $\hat{\mu}$ and a bound on
$\hat{\sigma}^2$. Essentially the same holds for linear regression
problems, Section~\ref{sec:regression}. 

\paragraph{Renormalized Maximum Likelihood}
Meta-two-part coding is just one possible solution to the problem of
undefined $\complexity_n(\M)$. It is suboptimal, the main
reason being the use of 2-part codes.  Indeed, these 2-part codes
are not complete (Section~\ref{sec:prob}): they reserve several codewords for the same data $D = (x_1, \ldots, x_n)$
(one for each integer value of $\log K$); therefore, there must exist more
efficient (one-part) codes $\punivmeta'$ such that for all $x^n \in
{\cal X}^n$, $\punivmeta'(x^n) > \punivmeta(x^n)$; in keeping with the
idea that we should minimize description length, such alternative
codes are preferable.  This realization has led to a search for more
efficient and intrinsic solutions to the problem.  \citeN{FosterS01}
consider the possibility of restricting the {\em parameter values\/}
rather than the data, and develop a general framework for comparing
universal codes for models with undefined $\complexity(\M)$.
\citeN{Rissanen01} suggests the following elegant solution.  He
defines the {\em Renormalized Maximum Likelihood (RNML) distribution}
$\punivrnml$.  In our Gaussian example, this universal model would be
defined as follows. Let $\hat{K}(x^n)$ be the bound on
$\hat{\mu}(x^n)$ that maximizes $\punivnml(x^n \mid \M_K)$ for the
actually given $K$. That is, $\hat{K}(x^n) = | \hat{\mu}(x^n)|$. Then
$\punivrnml$ is defined as, for all $x^n \in {\cal X}^n$,
\begin{equation}
\label{eq-rnml}
\punivrnml(x^n | \M) =  
\frac{\punivnml(x^n | \M_{\hat{K}(x^n)})}
{\int_{x^n \in \reals^n} \punivnml(x^n \mid \M_{\hat{K}(x^n)}) dx^n}.
\end{equation}
Model selection between a finite set of models now proceeds by
selecting the model maximizing the {\em re-\/}normalized likelihood
(\ref{eq-rnml}).  \index{Rissanen renormalization} \index{Rissanen
  complexity} \index{renormalized Rissanen complexity} 

\paragraph{Region Indifference} All the
approaches considered thus far slightly prefer some regions of the
parameter space over others.  In spite of its elegance, even the
Rissanen renormalization is slightly `arbitrary' in this way: had we
chosen the origin of the real line differently, the same sequence
$x^n$ would have achieved a different codelength $- \log
\punivrnml(x^n \mid \M)$. In recent work, Liang and Barron
\citeyear{LiangB04,LiangB03} consider a novel and quite different
approach for dealing with infinite $\complexity_n(\M)$ that partially
addresses this problem.  They make use of the fact that, while
Jeffreys' {\em prior\/} is improper ($\curvelength{}$ is infinite),
using Bayes' rule we can still compute Jeffreys' {\em posterior\/}
based on the first few observations, and this posterior turns out to
be a proper probability measure after all.  Liang and Barron use
universal models of a somewhat different type than $\punivnml$, so it
remains to be investigated whether their approach can be adapted to
the form of MDL discussed here.
\subsection{The General Picture}
\label{sec:general}
Section~\ref{sec:uncountable} illustrates that, in {\em all\/}
applications of MDL, we first define a {\em single\/} universal model
that allows us to code all sequences with length equal to the given
sample size.  If the set of models is finite, we use the uniform
prior.  We do this in order to be as `honest' as possible, treating
all models under consideration on the same footing. But if the set of
models becomes infinite, there exists no uniform prior any more.
Therefore, we must choose a non-uniform prior/non-fixed length code to
encode the model index.  In order to treat all models still `as
equally as possible', we should use some code which is `close' to
uniform, in the sense that the codelength increases only very slowly
with $k$.  We choose the standard prior for the integers
(Example~\ref{ex:integers}), but we could also have chosen different
priors, for example, a prior $P(k)$ which is uniform on $k = 1..M$ for
some large $M$, and $P(k) \propto k^{-2}$ for $k > M$. Whatever prior
we choose, we are forced to encode a slight preference of some models
over others; see Section~\ref{sec:perceived}.
\begin{figure}
\ownbox{GENERAL `REFINED' MDL PRINCIPLE for Model Selection}{
Suppose 
we plan to select between models $\MN{1},\MN{2}, \ldots$ for data $D= (x_1, \ldots, x_n)$.
MDL tells us to design a universal code $\puniv$ 
for ${\cal X}^n$, in which the index $k$ of $\MN{k}$ is encoded explicitly. The resulting code has two parts, the two sub-codes being defined such that
\begin{enumerate}
\item All models $\MN{k}$ are treated on the same footing, as far as
  possible: we assign a uniform prior  to these models, or, if that is not a possible, a prior `close to' uniform. 
\item All distributions within each $\MN{k}$ are treated on the same footing, as far as possible:
we use the minimax regret universal model $\punivnml(x^n \mid
\MN{k})$. If this model is undefined or too hard to compute, we
instead  use a different universal model that achieves regret `close
to' the minimax regret for each submodel of $\MN{k}$ in the sense of (\ref{eq:firstsolb}).
\end{enumerate}
In the end, we encode data $D$ using a hybrid two-part/one-part 
universal model, explicitly encoding the models we want to select between and implicitly encoding any distributions contained in those models.
}\vspace*{-0.5 cm}
\caption{\label{fig:unified} The Refined MDL Principle.}
\end{figure}
  
Section~\ref{sec:boundary} applies the same idea, but implemented at a
meta-level: we try to associate with $\MN{k}$ a code for encoding
outcomes in ${\cal X}^n$ that achieves uniform (= minimax) regret for
every sequence $x^n$. If this is not possible, we still try to assign
regret as `uniformly' as we can, by carving up the parameter space in
regions with larger and larger minimax regret, and devising a
universal code that achieves regret not much larger than the minimax
regret achievable within the smallest region containing the ML
estimator. Again, the codes we used encoded a slight preference of
some regions of the parameter space over others, but our aim was to
keep this preference as small as possible. The general idea is
summarized in Figure~\ref{fig:unified}, which provides an (informal)
definition of MDL, but only in a restricted context. If we go beyond
that context, these prescriptions cannot be used literally -- but
extensions in the same spirit suggest themselves.  Here is a first
example of such an extension:
\begin{example}
\label{ex:local}
{\bf \ [MDL and Local Maxima in the Likelihood]} \rm In practice we often work 
with models for which the ML estimator cannot
be calculated efficiently; or at least, no algorithm for efficient
calculation of the ML estimator is known. Examples are finite and Gaussian
mixtures and Hidden Markov models. In
such cases one typically resorts to methods such as EM or gradient
descent, which find a {\em local\/} maximum of the likelihood surface
(function) $P(x^n \mid \theta)$, leading to a {\em local\/} maximum
likelihood estimator (LML) $\dot{\theta}(x^n)$. Suppose we need to
select between a finite number of such models. We may be tempted to
pick the model $\M$ maximizing the normalized likelihood $\punivnml(x^n
\mid \M)$.  However, if we then plan to use the local estimator
$\dot{\theta}(x^n)$ for predicting future data, this is {\em not\/}
the right thing to do. To see this, note
that, if suboptimal estimators $\dot{\theta}$ are to be used, the
ability of model $\M$ to fit arbitrary data patterns may be severely
diminished! 
%As an extreme case, if the estimator $\dot{\theta}$ is so
%bad that $\dot{\theta}(x^n)$ takes on the same value for each $x^n$,
%then $\M$ reduces to a single distribution with minimal overfitting
%capability. 
Rather than using $\punivnml$, 
we should redefine it to take into account the fact that $\dot{\theta}$ is not the global ML estimator:
%\begin{equationr}
$$
\punivnml'(x^n) \isbydefinition \frac{P(x^n \mid \dot{\theta}(x^n))}
{\sum_{x^n \in {\cal X}^n} P(x^n \mid \dot{\theta}(x^n))},
%\end{equationr}
$$
leading to an adjusted parametric complexity
\begin{equation}
\label{eq:localnml}
\complexity'_n(\M) \isbydefinition 
\log \sum_{x^n \in {\cal X}^n} P(x^n \mid \dot{\theta}(x^n)), 
\end{equation}
which, for every estimator $\dot{\theta}$ different from
$\hat{\theta}$ {\em must\/} be strictly smaller than $\complexity_n(\M)$. 
\end{example}
\paragraph{Summary}
We have shown how to extend refined MDL beyond the restricted settings
of Section~\ref{sec:refined}. This uncovered the general principle
behind refined MDL for model selection, given in
Figure~\ref{fig:unified}. General as it may be, it only applies to
model selection -- in the next section we briefly discuss extensions
to other applications.
\section{Beyond Parametric Model Selection}
\label{sec:beyond}
The general principle as given in
Figure~\ref{fig:unified} only applies to model selection. It can be
extended 
in several
directions. These range over many
different tasks of inductive inference -- we mention {\em prediction},
{\em transduction\/} (as defined in \cite{Vapnik98}), {\em
  clustering\/} \cite{KontkanenMBRT04} and {\em similarity
  detection\/} \cite{LiCLV03}. In these areas there has been less
research and a `definite' MDL approach has not yet been formulated. 

MDL {\em has\/} been developed in some detail for some other inductive
tasks: {\em non-parametric\/} inference, {\em parameter
  estimation\/} and {\em regression\/} and {\em classification\/}
problems.  We give a very brief overview of these - for details we
refer to \cite{BarronRY98,HansenY01} and, for the classification case,
\cite{GrunwaldL04}.

\paragraph{Non-Parametric Inference}
Sometimes the model class $\M$ is so large that it cannot be finitely
parameterized. For example, let ${\cal X} = [0,1]$ be the unit
interval and let $\M$ be the i.i.d. model consisting of {\em all\/}
distributions on ${\cal X}$ with densities $f$ such that $-\log f(x)$
is a continuous function on ${\cal X}$.  $\M$ is clearly
`non-parametric': it cannot be meaningfully parameterized by a
connected finite-dimensional parameter set $\Theta^{(k)} \subseteq
\reals^k$. We may still try to learn a distribution from $\M$ in
various ways, for example by {\em histogram density estimation}
\cite{RissanenSY92} or {\em kernel density estimation\/}
\cite{Rissanen89}. MDL is quite suitable for such applications, in
which we typically select a density $f$ from a class $\MN{n} \subset
\M$, where $\MN{n}$ grows with $n$, and every $P^* \in \M$ can be
arbitrarily well approximated by members of $\MN{n}, \MN{n+1}, \ldots$
in the sense that $\lim_{n \rightarrow \infty} \inf_{P \in \MN{n}}
D(P^* \| P) = 0$ \cite{BarronRY98}.  Here $D$ is the {\em
  Kullback-Leibler\/} divergence \cite{CoverT91} between $P^*$ and
$P$.
\paragraph{MDL Parameter Estimation: Three Approaches} The `crude' MDL method
(Section~\ref{sec:crude}) was a means of doing model selection and
parameter estimation at the same time. `Refined' MDL only dealt with
selection of {\em models}. If instead, or at the same time, parameter
estimates are needed, they may be obtained in three different ways.
Historically the first way \cite{Rissanen89,HansenY01} was to simply
use the refined MDL Principle to pick a parametric model $\MN{k}$, and
then, within $\MN{k}$, pick the ML estimator $\hat{\theta}^{(k)}$.
After all, we associate with $\MN{k}$ the distribution $\punivnml$
with codelengths `as close as possible' to those achieved by the ML
estimator. This suggests that within $\MN{k}$, we should prefer the ML
estimator. But upon closer inspection, Figure~\ref{fig:unified}
suggests to use a two-part code also to select $\theta$ within
$\MN{k}$; namely, we should discretize the parameter space in such a
way that the resulting 2-part code achieves the minimax regret among
all two-part codes; we then pick the (quantized) $\theta$ minimizing
the two-part code length. Essentially this approach has been worked
out in detail by \citeN{BarronC91}.  The resulting estimators may be
called {\em two-part code MDL estimators}.  A third possibility is to
define {\em predictive\/} MDL estimators such as the Laplace and
Jeffreys estimators of Example~\ref{ex:laplace}; once again, these can
be understood as an extension of Figure~\ref{fig:unified}
\cite{BarronRY98}. These second and third possibilities are more
sophisticated than the first. However, if {\em the model $\M$ is
  finite-dimensional parametric and $n$ is large}, then both the
two-part and the predictive MDL estimators will become
indistinguishable from the maximum likelihood estimators. For this
reason, it has sometimes been claimed that MDL parameter estimation is
just ML parameter estimation.  Since for small samples, the estimates
can be quite different, this statement is misleading.

\commentout{
two-part code parameter estimation:
not the same as maximum likelihood (but precision chosen so that it
will not do something very different) 

My own view: predictive parameter estimation and `model selection'
with prequential codes or their cesaro average: unified view, can show
consistency if code is universal model, avoid diaconis freedman.
non-predictive paramter/model estimation: more in Rissanen's style, to
'gain insight', there will always be the problem that you assign
probability $0$ to something which then happens.  }
\paragraph{Regression}
\label{sec:regression}
In regression problems we are interested in learning how the values $y_1, \ldots, y_n$ of a {\em
  regression\/} variable $Y$ depend on the values $x_1, \ldots, x_n$
of the {\em regressor\/} variable $X$. We assume or hope that there exists
some function $h: {\cal X} \rightarrow {\cal Y}$ so that $h(X)$
predicts the value $Y$ reasonably well, and we want to learn such an
$h$ from data. To this end, we assume a set of {\em candidate
  predictors\/} (functions) ${\cal H}$.  In
Example~\ref{ex:modsel}, we took ${\cal H}$ to be the set of all
polynomials. In the standard
formulation of this problem, we take $h$ to express that
\begin{equation}
\label{eq:funcreg}
Y_i = h(X_i) + Z_i,
\end{equation}
where the $Z_i$ are i.i.d. Gaussian random variables with mean $0$ and
some variance $\sigma^2$, independent of $X_i$. That is, we assume
Gaussian noise: (\ref{eq:funcreg})
implies that the conditional density of $y_1, \ldots, y_n$, 
given $x_1, \ldots, x_n$, is equal to
the product of $n$ Gaussian densities:
\begin{equation}
\label{eq:densreg}
P(y^n \mid x^n, \sigma, h) = \biggl( \frac{1}{\sqrt{2\pi} \sigma} \biggr)^n 
\exp \biggl(- \frac{\sum_{i=1}^n (y_i - h(x_i))^2}{2 \sigma^2} \biggr).
\end{equation}
With this choice, the log-likelihood becomes a linear function of the
squared error:
\begin{equation}
\label{eq:gauss}
- \log P(y^n \mid x^n, \sigma, h)= \frac{1}{2 \sigma^2} \sum_{i=1}^n (y_i -
  h(x_i))^2 + \frac{n}{2} \log 2 \pi \sigma^2.
\end{equation}
Let us now assume that ${\cal H} = \cup_{k \geq 1} {\cal H}^{(k)}$
where for each $k$, ${\cal H}^{(k)}$ is a set of functions $h: {\cal
  X} \rightarrow {\cal Y}$. For example, ${\cal H}^{(k)}$ may be the
set of $k$-th degree polynomials.

With each model ${\cal H}^{(k)}$ we can associate a set of densities
(\ref{eq:densreg}), one for each $(h,\sigma^2)$ with $h \in {\cal
  H}^{(k)}$ and $\sigma^2 \in \reals^+$. Let $\MN{k}$ be the resulting
set of conditional distributions. Each $P(\cdot \mid h,\sigma^2) \in
\MN{k}$ is identified by the parameter vector
$(\alpha_0,\ldots,\alpha_{k},\sigma^2)$ so that $h(x) \isbydefinition
\sum_{j=0}^{k} \alpha_j x^j$. 
%Note that this makes $\MN{k}$ is a $k+2$-dimensional model.
By Section~\ref{sec:uncountable}, (\ref{eq:twopartmdl}) MDL tells us to select
the model minimizing
\begin{equation}
\label{eq:basicreg}
- \log \puniv(y^n \mid \MN{k},x^n) + L(k)
\end{equation}
where we may take $L(k) = 2 \log k +1$, and $\puniv(\cdot \mid \MN{k},
\cdot)$ is now a {\em conditional\/} universal model with small minimax
regret.  (\ref{eq:basicreg}) ignores the
codelength of $x_1, \ldots, x_n$. Intuitively, this is
because we are only interested in learning how $y$ {\em depends\/} on $x$;
therefore, we do not care how many bits are needed to encode $x$.
Formally, this may be understood as follows: we really {\em are\/}
encoding the $x$-values as well, but we do so using a fixed code that
does not depend on the hypothesis $h$ under consideration. Thus, we
are really trying to find the model $\MN{k}$ minimizing
$$
%\begin{equationr}
%\label{eq:ignorex}
- \log \puniv(y^n \mid \MN{k},x^n) + L(k) + L'(x^n)
%\end{equationr}
$$
where $L'$ represents some code for ${\cal X}^n$. 
Since this codelength does not involve $k$,
it can be dropped from the minimization; see Figure~\ref{fig:ignore}.
\begin{figure}
\ownbox{When the Codelength for $x^n$ Can Be Ignored}{
If all models under consideration represent {\em conditional\/}
densities or probability mass functions $P(Y\mid X)$, then the
codelength for $X_1, \ldots, X_n$ can be ignored in model and
parameter selection. Examples are applications of MDL in {\em
  classification\/} and {\em regression}.
}\vspace*{-0.5 cm}
\caption{\label{fig:ignore} Ignoring codelengths.}
\end{figure}
We will not go into the precise definition of $\puniv(y^n \mid
\MN{k},x^n)$. Ideally, it should be an NML distribution, but just as
in Example~\ref{ex:boundary}, this NML distribution is not well-defined. We can
get reasonable alternative universal models after all using any of the
methods described in Section~\ref{sec:boundary}; see \cite{BarronRY98}
and \cite{Rissanen00} for details.
\paragraph{`Non-probabilistic' Regression and Classification}
In the approach we just described, we modeled the noise as being
normally distributed. Alternatively, it has been tried to {\em
  directly\/} try to learn functions $h \in {\cal H}$ from the data,
without making any probabilistic assumptions about the noise
\cite{Rissanen89,Barron91,Yamanishi98,Grunwald98b,Grunwald99a}.  The
idea is to learn a function $h$ that leads to good predictions of
future data from the same source in the spirit of Vapnik's \citeyear{Vapnik98} {\em statistical learning theory}. Here prediction quality is measured
by some fixed loss function; different loss functions lead to
different instantiations of the procedure. Such a
version of MDL is meant to be more robust, leading to inference of a `good'
$h \in {\cal H}$ irrespective of the details of the noise
distribution. This loss-based approach has also been the method of
choice in applying MDL to {\em classification\/} problems. Here ${\cal
  Y}$ takes on values in a finite set, and the goal is to match each
{\em feature\/} $X$ (for example, a bit map of a handwritten digit)
with its corresponding {\em label\/} or {\em class\/} (e.g., a digit).
While several versions of MDL for classification have been proposed
\cite{QuinlanR89,Rissanen89,KearnsMNR97}, most of these can be reduced to
the same approach based on a $0/1$-valued loss function 
\cite{Grunwald98b}. In recent work
\cite{GrunwaldL04} we show that this MDL approach to classification
without making assumptions about the noise may behave suboptimally: we
exhibit situations where no matter how large $n$, MDL keeps
overfitting, selecting an overly complex model with suboptimal
predictive behavior. Modifications of MDL suggested by
\citeN{Barron91} and \citeN{Yamanishi98} do not suffer from this
defect, but they do not admit a natural coding interpretation any
longer. All in all, current versions of MDL that avoid 
probabilistic assumptions are still in their infancy, and more research
is needed to find out whether they can be modified to perform well in
more general and realistic settings.
\paragraph{Summary}
In the previous sections, we have covered basic refined MDL
(Section~\ref{sec:refined}), general refined
MDL (Section~\ref{sec:unify}), and several
extensions of refined MDL (this section). This concludes our technical
description of refined MDL. It only remains to place MDL in its proper
context: what does it {\em do\/} compared to other methods of
inductive inference?  And how {\em
  well\/} does it {\em perform}, compared to other methods? The next
two sections are devoted to these questions.
\section{Relations to Other Approaches to Inductive Inference}
How does MDL compare to
other model selection and statistical inference methods? In order to
answer this question, we first have to be precise about what we mean
by `MDL'; this is done in Section~\ref{sec:whatis}. We then continue
in Section~\ref{sec:bayes} by summarizing MDL's relation to {\em Bayesian inference}, Wallace's {\em
  MML Principle}, Dawid's {\em prequential model validation}, {\em
  cross-validation\/} and an `idealized' version of MDL based on Kolmogorov complexity. 
The literature has also established connections between MDL and
Jaynes' \citeyear{Jaynes03}
{\em Maximum Entropy Principle\/} 
\cite{Feder86,LiV97,Grunwald98b,Grunwald00a,GrunwaldD04} and Vapnik's  \citeyear{Vapnik98} {\em structural risk
  minimization principle} \cite{Grunwald98b}, but there is no space here to discuss
these. Relations between MDL and Akaike's {\em AIC\/} \cite{BurnhamA02} are subtle. They are discussed by, for example, \citeN{SpeedY93}.

\commentout{Uiteindelijk elk onderwerp: in 3 stappen: (1) the method;
  how does it compare to MDL in practice? ; how general is it? (AIC
  and CV: only for model selection; Bayes, like MDL, for everything
  that can be handled probabilistic -- and Bayes would even claim that
  everything CAN be handled probabilistically ; SRM: everything that
  is 'not' probabilistic (2) is the procedure
  admissible according to Rissanen's views (as espoused in the little
  green book) and (3) my own view ;

make clear that Barron's universal models are NOT contradicting
Rissanen's philosophy (model = code) . 

Of course we WANT our approach to be consistent and have small
prediction errors in many different ways, but we cannot design the
approach on that. (and if our approach is not consistent, then maybe
there is something 'anomalous' about MDL...)  }
\label{sec:others}
\subsection{What is `MDL'?}
\label{sec:whatis}
`MDL' is used by different authors in somewhat different meanings. 
Some authors use MDL as a broad
umbrella term for all types of inductive inference based on data
compression.  This would, for example, include the `idealized'
versions of MDL based on Kolmogorov complexity and Wallaces's MML
Principle, to be discussed below. On the other extreme, for
historical reasons, some authors use the {\em MDL Criterion\/} to
describe a very specific (and often not very successful) model
selection criterion equivalent to BIC, discussed further below.

Here we adopt the meaning of the term that is embraced in the survey
\cite{BarronRY98}, written by arguably the three most important
contributors to the field: we use MDL for general {\em inference based
  on universal models}. These include, but are not limited to
approaches in the spirit of Figure~\ref{fig:unified}. For example,
some authors have based their inferences on `expected' rather than
`individual sequence' universal models \cite{BarronRY98,LiangB04}.
Moreover, if we go beyond model selection (Section~\ref{sec:beyond}),
then the ideas of Figure~\ref{fig:unified} have to be modified to some
extent. In fact, one of the main strengths of ``MDL'' in
this broad sense is that it can be applied to ever more exotic
modeling situations, in which the models do not resemble anything that
is usually encountered in statistical practice.  An example is the
model of context-free grammars, already suggested by
\citeN{Solomonoff64}. In this tutorial, we call applications of MDL
that strictly fit into the scheme of Figure~\ref{fig:unified} {\em
  refined MDL for model/hypothesis selection\/}; when we simply say `MDL', we mean `inductive inference based on
universal models'. This form of inductive inference goes hand in hand with Rissanen's radical MDL  {\em
  philosophy}, which views learning as finding
useful properties of the data, not necessarily related to the
existence of a `truth' underlying the data. This view was
outlined in Chapter~\ref{chap:survey}, Section~\ref{sec:philosophy}.
Although MDL practitioners and theorists are usually sympathetic to
it, the different interpretations of MDL listed in
Section~\ref{sec:refined} make clear that MDL applications can also be
justified without adopting such a radical philosophy.
\subsection{MDL and Bayesian Inference}
\label{sec:bayes}
\commentout{ One important part that is missing here: AGAIN, we are
  not allowed to take expectations, now over the prior!  (and this
  also invalidates mml!)  Reason: in the special case where Pnml is
  defined, we really don't know what generates the data so taking
  expectation is misleading.  In the case where it is not defined
  then, it least in some subcases of this case(for example, when a
  location parameter needs to be restricted (! actually a good example
  !), some subjectivity is INEVITABLE. Nevertheless, if we prefer a
  value \mu_0 over \mu_1, this means either 'we don't know anything
  else to do' or 'if the experiment were to be repeated, then *mean*
  in the area around mu_1 would turn up less often than *mean* in the
  area around \mu_0. It does NOT mean that the *normal distribution*
  with mean close to mu_1 turns up less than the *normal* distribution
  with mean close to mu_0, and this is what would be implied if were
  to take expectations over the prior! If we only use mu_0 and mu_1 in
  coding, they already work well if SOME distribution with these means
  generates data!

MY OWN OPINION: I have never found Savage's, De Finetti's
etc. characterizations convincing as saying that 'Bayes is the only
reasonable procedure'. They are impressive as characterizations of
Bayes - but if Bayes also has some bad properties (e.g. inconsistency)
then it is not clear why the pro-Bayes arguments should be more
compelling.
viz. Wasserman and Robins.
Also, does not explain the fact that Bayes often still works if the
prior is NOT right.
}
{\em Bayesian statistics\/} \cite{Lee97,BernardoS94} is one of the most
well-known, frequently and successfully applied paradigms of
statistical inference. It is often claimed that `MDL is really just
a special case of Bayes\endnote{The author has heard many people say this at many conferences. The reasons are probably historical: while the underlying philosophy has always been different, until Rissanen introduced the use of $\punivnml$, most actual implementations of MDL `looked' quite Bayesian.}'.  Although there are close similarities, this
is simply not true.  To see this quickly, consider the basic quantity in
refined MDL: the NML distribution $\punivnml$,
Equation~(\ref{eq:punivnml}). While $\punivnml$ -- although defined in a completely different manner -- turns out to be closely related to the Bayesian marginal likelihood, this is no longer the case for its `localized'
version (\ref{eq:localnml}). There is no mention of anything like this
code/distribution in any Bayesian textbook!  Thus, it must be the case
that Bayes and MDL are somehow different.
\paragraph{MDL as a Maximum Probability Principle}
For a more detailed analysis, we need to distinguish between the two
central tenets of modern Bayesian statistics: (1) Probability
distributions are used to represent uncertainty, and to serve as a
basis for making predictions; rather than standing for some imagined
`true state of nature'. (2) All inference and decision-making is done
in terms of prior and posterior distributions.  MDL sticks with (1)
(although here the `distributions' are primarily interpreted as
`codelength functions'), but not (2): MDL allows the use of arbitrary
universal models such as NML and prequential universal models; the
Bayesian universal model does not have a special status among these.
In this sense, Bayes offers the statistician {\em less\/} freedom in
choice of implementation than MDL.  In fact, MDL may be reinterpreted
as a {\em maximum probability principle}, where the maximum is
relative to some given model, in the worst-case over all sequences
(Rissanen \citeyear{Rissanen87,Rissanen89} uses the phrase `{\em
  global\/} maximum likelihood principle').  Thus, whenever the
Bayesian universal model is used in an MDL application, a prior should
be used that minimizes worst-case codelength regret, or equivalently,
maximizes worst-case relative probability. There is no comparable
principle for choosing priors in Bayesian statistics, and in this
respect, Bayes offers a lot {\em more\/} freedom than MDL.
\begin{ownquote}
\begin{example}
\label{ex:bayescraze}
\rm There is a conceptual problem with Bayes' use of prior
distributions: in practice, we very often want to use models which we
{\em a priori\/} know to be wrong -- see Example~\ref{ex:sanity}. If
we use Bayes for such models, then we are forced to put a prior
distribution on a set of distributions which we know to be wrong -
that is, we have degree-of-belief $1$ in something we know not to be
the case. From an MDL viewpoint, these priors are interpreted as tools
to achieve short codelengths rather than degrees-of-belief and there
is nothing strange about the situation; but from a Bayesian viewpoint,
it seems awkward.  To be sure, Bayesian inference often gives good
results even if the model $\M$ is known to be wrong; the point is that
(a) if one is a strict Bayesian, one would never apply Bayesian
inference to such misspecified $\M$, and (b), the Bayesian theory
offers no clear explanation of why Bayesian inference might still give
good results for such $\M$. MDL provides both codelength and 
predictive-sequential interpretations of Bayesian inference, which
help explain why Bayesian inference may do something reasonable even if
$\M$ is misspecified.  
To be fair, we should add that there exists variations of the Bayesian
philosophy (e.g. \citeN{DeFinetti74}'s) which avoid the conceptual
problem we just described.
\end{example}
\end{ownquote}
\paragraph{MDL and BIC}
\index{BIC}
\index{AIC}
\index{Bayesian Information Criterion}
In the first paper on MDL, \citeN{Rissanen78} used a two-part code and showed that, asymptotically, and under regularity conditions,
the two-part codelength of $x^n$ based on a $k$-parameter model $\M$ with an optimally discretized parameter space is given by 
\begin{equation}
\label{eq:BIC}
- \log P(x^n \mid \hat{\theta}(x^n)) + \frac{k}{2} \log n,
\end{equation}
thus ignoring $O(1)$-terms, which, as we have already seen, can be
quite important. In the same year \citeN{Schwarz78} showed that, for
large enough $n$, Bayesian model selection between two exponential
families amounts to selecting the model minimizing (\ref{eq:BIC}),
ignoring $O(1)$-terms as well. As a result of Schwarz's paper, model
selection based on (\ref{eq:BIC}) became known as the {\em BIC
  (Bayesian Information Criterion)}. Not taking into account the
functional form of the model $\M$, it often does not work very well in practice. 

It has sometimes been claimed that MDL = BIC; for example, \cite[page
286]{BurnhamA02} write ``Rissanen's result is equivalent to BIC''.
This is wrong, even for the 1989 version of MDL that Burnham and
Anderson refer to -- as pointed out by \citeN{FosterS04}, the BIC
approximation only holds if the number of parameters $k$ is kept fixed
and $n$ goes to infinity.  If we select between nested families of
models where the maximum number of parameters $k$ considered is either
infinite or grows with $n$, then model selection based on both
$\punivnml$ and on $\punivbayes$ tends to select quite different
models than BIC - if $k$ gets closer to $n$, the contribution to
$\complexity_n(\M)$ of each additional parameter becomes much smaller
than $0.5 \log n$ \cite{FosterS04}. However, researchers who claim MDL
= BIC have a good excuse: in early work, Rissanen himself
has used the phrase `MDL criterion' to refer to (\ref{eq:BIC}), and
unfortunately, the phrase has stuck.
\paragraph{MDL and MML}
MDL shares some ideas with the {\em Minimum Message Length (MML)
  Principle\/} which predates MDL by 10 years.  
Key references are \cite{WallaceB68,WallaceB75}
and \cite{WallaceF87}; a long list is in \cite{ComleyD04}.  \index{Minimum
  Message Length} Just as in MDL, MML chooses the hypothesis
minimizing the code-length of the data. But the {\em codes\/} that
are used are quite different from those in MDL. First of all, in MML
one {\em always\/} uses two-part codes, so that MML automatically selects both a
model family and parameter values. Second, while the MDL codes such as
$\punivnml$ minimize {\em worst-case relative code-length\/} (regret),
the two-part codes used by MML are designed to minimize {\em expected
  absolute\/} code-length. Here the expectation is taken over a
subjective prior distribution defined on the collection of models and
parameters under consideration. While this approach
contradicts Rissanen's philosophy, in practice it
often leads to similar results.

Indeed, Wallace and his co-workers stress that their approach is fully
(subjective) {\em Bayesian}. Strictly speaking, a Bayesian should
report his findings by citing the full posterior distribution. But
sometimes one is interested in a single model, or hypothesis for the
data. A good example is the inference of phylogenetic trees in
biological applications: the full posterior would consist of a mixture
of several of such trees, which might all be quite different from each
other. Such a mixture is almost impossible to interpret -- to get
insight in the data we need a single tree.  In that case, Bayesians
often use the MAP (Maximum A Posteriori) hypothesis which maximizes
the posterior, or the posterior mean parameter value. The first
approach has some unpleasant properties, for example, it is not
invariant under reparameterization. The posterior mean approach cannot
be used if different model families are to be compared with each
other. The MML method provides a theoretically sound way of proceeding
in such cases.
\begin{figure}
%\vspace*{-1.5 cm}
\centerline{
\epsfxsize=6cm
\epsfbox{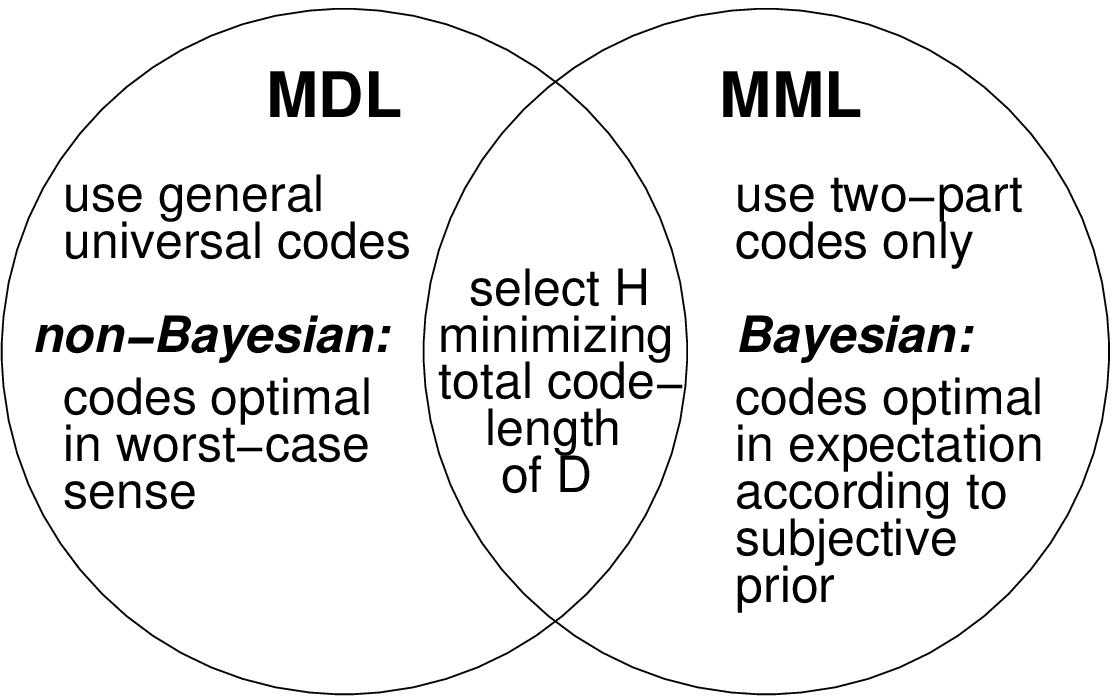}
\epsfxsize=6cm
%\hspace*{-3 cm}
\epsfbox{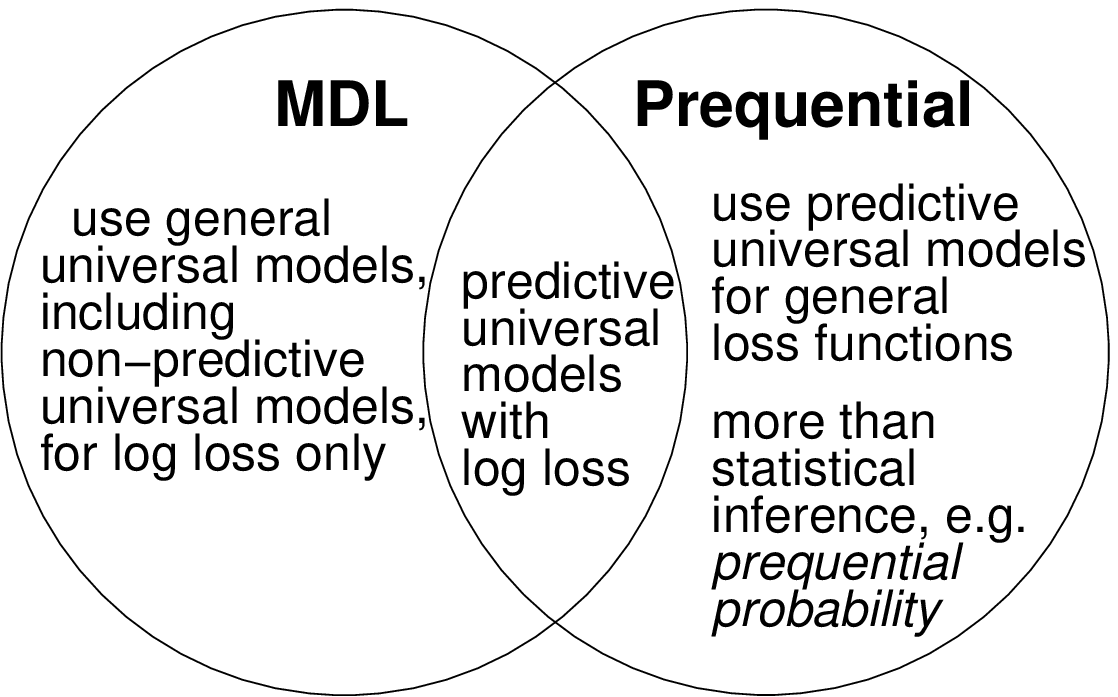}} 
\caption{\label{fig:mml} Rissanen's MDL, Wallace's MML and Dawid's
  Prequential Approach.}
\end{figure}
\index{MML|see{Minimum Message Length}}
\subsection{MDL, Prequential Analysis and Cross-Validation}
\index{prequential analysis} In a series of papers, A.P. Dawid
\citeyear{Dawid84,Dawid92,Dawid97} put forward a methodology for
probability and statistics based on sequential prediction which he
called the {\em prequential approach}. When applied to model selection
problems, it is closely related to MDL - Dawid proposes to construct,
for each model $\MN{j}$ under consideration, a `probability
forecasting system' (a sequential prediction strategy) where the
$i+1$-st outcome is predicted based on either the Bayesian posterior
$\punivbayes(\theta | x^i)$ or on some estimator $\hat{\theta}(x^i)$.
Then the model is selected for which the associated sequential
prediction strategy minimizes the accumulated prediction error.
Related ideas were put forward by \citeN{Hjorth82} under the name {\em
  forward validation} \index{forward validation} and
\citeN{Rissanen84}. From Section~\ref{sec:prequential} we see
that this is just a form of MDL - strictly speaking, {\em every\/}
universal code can be thought of as as prediction strategy, but for
the Bayesian and
  the plug-in universal models (Sections~\ref{sec:bayes1},
  \ref{sec:prequential}) the interpretation is much more natural than for
others\endnote{The reason is that the Bayesian and plug-in models
can be
  interpreted as probabilistic sources. The NML and the
  two-part code models are no probabilistic sources, since 
 $\puniv^{(n)}$ and $\puniv^{(n+1})$
  are not compatible in the sense of Section~\ref{sec:prob}.}. Dawid mostly talks about such `predictive' universal models.   On
the other hand, Dawid's framework allows to adjust the
prediction loss to be measured in terms of arbitrary loss functions,
not just the log loss. In this sense, it is more general than MDL.
Finally, the prequential idea goes beyond statistics: there is
also a `prequential approach' to probability theory developed by Dawid
\cite{DawidV99} and \citeN{ShaferV01}.

Note that the prequential approach is similar in spirit to
cross-validation. In this sense MDL is related to cross-validation as
well. The main differences are that in MDL and the prequential
approach, 
(1) all predictions are done
{\em sequentially\/} (the future is never used to predict the past),
and (2) each outcome is predicted {\em exactly once}.

\subsection{Kolmogorov Complexity and Structure Function; Ideal MDL}
\label{sec:kolmogorov}
Kolmogorov complexity \cite{LiV97} has played a large but mostly
inspirational role  in Rissanen's development of MDL. Over the last fifteen years, several `idealized' versions
of MDL have been proposed, which are more directly based on Kolmogorov
complexity theory
\cite{Barron85,BarronC91,LiV97,VereshchaginV02}. These are all based
on two-part codes, where hypotheses are described using a universal
programming language such as {\sc C} or {\sc Pascal}. 
For example, in one proposal \cite{BarronC91}, given data $D$ one
picks the distribution minimizing
\begin{equation}
\label{eq:twopideal}
K(P) + \bigl[ - \log P(D) \bigr],
\end{equation}
where the minimum is taken over {\em all\/} computable probability
distributions, and $K(P)$ is the length of the shortest computer
program that, when input $(x,d)$, outputs $P(x)$ to $d$ bits
precision.  While such a procedure is mathematically well-defined, it
cannot be used in practice. The reason is that in general, the $P$
minimizing (\ref{eq:twopideal}) cannot be effectively computed.
Kolmogorov himself used a variation of (\ref{eq:twopideal}) in which
one adopts, among all $P$ with $K(P) - \log P(D) \approx K(D)$, the
$P$ with smallest $K(P)$. Here $K(D)$ is the Kolmogorov complexity of
$D$, that is, the length of the shortest computer program that prints
$D$ and then halts. This approach is known as the Kolmogorov {\em
  structure function\/} or {\em minimum sufficient statistic\/}
approach \cite{Vitanyi04}. In this approach, the idea of separating
data and noise (Section~\ref{sec:compression}) is taken as basic, and
the hypothesis selection procedure is defined in terms of it.  The
selected hypothesis may now be viewed as capturing all structure
inherent in the data - given the hypothesis, the data cannot be
distinguished from random noise.  Therefore, it may be taken as a
basis for {\em lossy\/} data compression -- rather than sending the
whole sequence, one only sends the hypothesis representing the
`structure' in the data. The receiver can then use this hypothesis to
generate `typical' data for it - this data should then `look just the
same' as the original data D.  Rissanen views this separation idea as
perhaps the most fundamental aspect of `learning by compression'.
Therefore, in recent work he has tried to relate MDL (as defined here,
based on lossless compression) to the Kolmogorov structure function,
thereby connecting it to lossy compression, and, as he puts it,
`opening up a new chapter in the MDL theory'
\cite{VereshchaginV02,Vitanyi04,RissanenT04}.
\paragraph{Summary and Outlook}
We have shown that MDL is closely related to, yet distinct from,
several other methods for inductive inference. In the next section we
discuss how well it {\em performs\/} compared to such other methods.
\section{Problems for MDL?}
\label{sec:problems}
Some authors have criticized MDL either on conceptual grounds (the
idea makes no sense) \cite{Webb96,Domingos99} or on practical grounds (sometimes it does not
work very well in practice) \cite{KearnsMNR97,Pednault03}. Are these
criticisms justified? Let us consider them in turn.
\subsection{Conceptual Problems: Occam's Razor}
\label{sec:perceived}
The most-often heard 
conceptual criticisms are invariably related to Occam's razor. We
have already discussed in Section~\ref{sec:philosophy} of the previous
chapter why we regard these criticisms as being entirely
mistaken. Based on our newly acquired technical knowledge of MDL, 
let us discuss these criticisms a little bit further:
\paragraph{1. `Occam's Razor (and MDL) is arbitrary' (page~\pageref{page:arbitrary})}
If we restrict ourselves to refined MDL for
comparing a finite number of models for which the NML distribution is
well-defined, then there is {\em nothing\/} arbitrary about MDL - it
is exactly clear what codes we should use for our inferences. 
The NML distribution and its close cousins, the Jeffreys' prior
marginal likelihood $\punivbayes$ and the asymptotic expansion
(\ref{eq:asnml}) are all invariant to continuous 1-to-1
reparameterizations of the model: parameterizing our model in a
different way (choosing a different `description language') does not
change the inferred description lengths.

If we go beyond models for which the NML distribution is defined,
and/or we compare an infinite set of models at the same time, then
some `subjectivity' {\em is\/} introduced -- while there are still
tough restrictions on the codes that we are allowed to use, all such
codes prefer some hypotheses in the model over others. If one does not
have an a priori preference over any of the hypotheses, one may
interpret this as some arbitrariness being added to the procedure.
But this `arbitrariness' is of an infinitely milder sort than the
arbitrariness that can be introduced if we allow completely arbitrary codes for
the encoding of hypotheses as in crude two-part code MDL,
Section~\ref{sec:crude}. 

%However, the choice of codes we are allowed to use is still very much
%restricted: we are only allowed to use universal codes, and we have a
%clear preference for universal codes where the regret increases as
%slowly as possible as we move on to larger regions of the parameter
%space. 
\begin{ownquote}
\rm
%Moreover, we add that {\em every\/} method for selecting from
%such a large model suffers from this type of `subjectivity': if 
Things get more subtle if we are interested not in model selection
(find the best order Markov chain for the data) but in
infinite-dimensional 
estimation (find the best Markov chain parameters for the data, among
the set $\Bernoulli$ of all Markov chains of each order). In the
latter case, if we are to apply MDL, we somehow have to carve up
$\Bernoulli$ into subsets $\MN{0} \subseteq \MN{1} \subseteq \ldots
\subseteq \Bernoulli$. Suppose that we have already chosen $\MN{1} =
\BernoulliN{1}$ as the set of $1$-st order Markov chains. We normally
take $\MN{0} = \BernoulliN{0}$ , the set of $0$-th order Markov chains
(Bernoulli distributions). But we could also have defined $\MN{0}$ as
the set of all $1$-st order Markov chains with $P(X_{i+1} = 1 \mid
X_{i} = 1) = P(X_{i+1} = 0 \mid X_{i} = 0)$.  This defines a
one-dimensional subset of $\BernoulliN{1}$ that is {\em not\/} equal
to $\BernoulliN{0}$. While there are several good reasons\endnote{For
  example, $\BernoulliN{0}$ is better interpretable.} for choosing
$\BernoulliN{0}$ rather than $\MN{0}$, there may be no indication that
$\BernoulliN{0}$ is somehow a priori more likely than $\MN{0}$. While
MDL tells us that we somehow have to carve up the full set
$\Bernoulli$, it does not give us precise guidelines on how to do this
-- different carvings may be equally justified and lead to different
inferences for small samples.  In this sense, there is indeed some
form of arbitrariness in this type of MDL applications.  But this is
unavoidable: we stress that this type of arbitrariness is enforced by
{\em all\/} combined model/parameter selection methods - whether they
be of the Structural Risk Minimization type \cite{Vapnik98}, AIC-type
\cite{BurnhamA02}, cross-validation or any other type. The only
alternative is treating all hypotheses in the huge class $\Bernoulli$ on the same footing,
which amounts to maximum likelihood estimation and extreme
overfitting.
\end{ownquote}
\commentout{
\begin{ownquote}
  In the past, MDL has often been identified with the `idealized'
  version based on Kolmogorov complexity \cite{LiV97} or the two-part
  version of Section~\ref{sec:crude} without any further specification
  of the code-length function for the hypothesis. In these versions,
  `complexity' is a property of a {\em single\/} hypothesis, rather
  than a family of hypotheses, and as such there is a some
  arbitrariness to the definition of complexity: different codes will
  lead to different complexities.  However, both in the Kolmogorov
  complexity approach, and in the versions of the two-part approach
  described by \citeN{Rissanen78,Rissanen83} this arbitrariness is
  within strong limits - for any two codes that are allowed, $|
  L_1(x^n) - L_2(x^n) | < \const$ for some constant $\const$ not
  depending on $n$ Doen WAAR IS HYPOTHESIS IN CODE? Nevertheless, the
  constant $\const$ is typically unknown and can be quite large, so
  that for practical sample sizes there is a real problem. Once again,
  we stress that there are no `unknown' constants in refined MDL;
  quite on the contrary, the constant term in (\ref{eq:asnml}) may be
  decisive in choosing a model for the data.
\end{ownquote}
}
\paragraph{2. `Occam's razor is false' (page~\pageref{page:false})}
We often try to model real-world situations that can be arbitrarily
complex, so why should we favor simple models? We gave in informal
answer on page~\pageref{page:false} where we claimed that {\em even if
  the true data generating machinery is very complex, it may be a good
  strategy to prefer simple models for small sample sizes.}

We are now in a position to give one formalization of this informal
claim: it is simply the fact that MDL procedures, with their built-in
preference for `simple' models with small parametric complexity, are
typically {\em statistically consistent\/} achieving {\em good rates
  of convergence\/} (page~\pageref{page:consistent}), whereas methods
such as maximum likelihood which do not take model complexity into
account are typically {\em in\/}-consistent whenever they are applied
to complex enough models such as the set of polynomials of each degree
or the set of Markov chains of all orders. This has implications for the quality of predictions: with complex enough models, no matter how many
training data we observe, if we use the maximum likelihood
distribution to predict future data from the same source, the
prediction error we make will not converge to the prediction error
that could be obtained if the true distribution were known; if we use an MDL submodel/parameter estimate (Section~\ref{sec:beyond}), the prediction error {\em will\/} converge to this optimal achieveable error.

Of course, consistency is not the only desirable property of a
learning method, and it may be that in some particular settings, and
under some particular performance measures, some alternatives to MDL
outperform MDL. Indeed this can happen -- see below. Yet it remains
the case that all methods the author knows of that successfully deal
with models of arbitrary complexity have a built-in preference for
selecting simpler models at small sample sizes -- methods such as
Vapnik's \citeyear{Vapnik98} Structural Risk Minimization, penalized
minimum error estimators \cite{Barron91} and the Akaike criterion
\cite{BurnhamA02} all trade-off complexity with error on the data, the
result invariably being that in this way, good convergence properties
can be obtained. While these approaches measure `complexity' in a
manner different from MDL, and attach different relative weights to
error on the data and complexity, the fundamental idea of finding a
{\em trade-off\/} between `error' and `complexity' remains.
\subsection{Practical Problems with MDL}
We just described some perceived problems about MDL. Unfortunately,
there are also some real ones: MDL is not a perfect method. While in
many cases, the methods described here perform very well\endnote{We
  mention \cite{HansenY00,HansenY01} reporting excellent behavior of
  MDL in regression contexts; and
  \cite{VanAllenMG03,KontkanenMST99,ModhaM98} reporting excellent
  behavior of predictive (prequential) coding in Bayesian network
  model selection and regression. Also, 
`objective Bayesian' model selection methods are frequently and
successfully used in practice \cite{KassW96}. Since these are based on
  non-informative priors such as Jeffreys', they often coincide with a
  version of refined MDL and thus indicate successful performance of
  MDL.}
there are also cases where they perform suboptimally compared to other
state-of-the-art methods. Often this is due to one of two reasons: 
\begin{enumerate}
\item An asymptotic formula like (\ref{eq:asnml}) was used and the
  sample size was not large enough to justify this
\cite{Navarro04}.
\item $\punivnml$ was undefined for the models under consideration,
  and this was solved by cutting off the parameter ranges at {\em ad
    hoc\/} values \cite{Lanterman04}.
\end{enumerate}
In these cases the problem probably lies with the use of invalid
approximations rather than with the MDL idea itself. More research is
needed to find out when the asymptotics and other approximations 
can be trusted, and what is
the `best' way to deal with undefined $\punivnml$. For the time being,
we suggest to avoid using (\ref{eq:asnml}) whenever possible, and to
never cut off the parameter ranges at arbitrary values -- instead, if
$\complexity_n(\M)$ becomes infinite, then some of the methods
described in Section~\ref{sec:boundary} should be used. Given these
restrictions, $\punivnml$ and Bayesian inference with Jeffreys' prior
are the preferred methods, since they both achieve the minimax regret.
If they are either ill-defined or computationally prohibitive for the
models under consideration, one can use a prequential method or a
sophisticated two-part code such as described by \citeN{BarronC91}.

\paragraph{MDL and Misspecification} 
However, there is a class of problems where MDL is problematic in a
more fundamental sense. Namely, if none of the distributions under
consideration represents the data generating machinery very well, then
both MDL and Bayesian inference may sometimes do a bad job in finding
the `best' approximation within this class of not-so-good hypotheses.
This has been observed in practice\endnote{But see
  \citeN{ViswanathanWDK99} who point out that the problem of
  \cite{KearnsMNR97} disappears if a more reasonable coding scheme is
  used.} \cite{KearnsMNR97,Clarke02,Pednault03}. \citeN{GrunwaldL04}
show that MDL can behave quite unreasonably for some classification
problems in which the true distribution is not in $\M$. This is
closely related to the problematic behavior of MDL for classification
tasks as mentioned in Section~\ref{sec:regression}.  All this is a bit
ironic, since MDL was explicitly designed {\em not\/} to depend on the
untenable assumption that some $P^* \in \M$ generates the data. But
empirically we find that while it generally works quite well if some
$P^* \in \M$ generates the data, it may sometimes fail if this is not
the case.
\section{Conclusion}
\label{sec:conclusion}
MDL is a versatile method for inductive inference: it can be
interpreted in at least four different ways, all of which indicate
that it does something reasonable. It is typically asymptotically
consistent, achieving good rates of convergence. It achieves all this
{\em without\/} having been designed for consistency, being based on a
philosophy which makes no metaphysical assumptions about the existence
of `true' distributions.  All this strongly suggests that it is a good
method to use in practice.  Practical evidence shows that in many
contexts it is, in other contexts its behavior can be problematic. In
the author's view, the main challenge for the future is to improve MDL
for such cases, by somehow extending and further refining MDL
procedures in a non ad-hoc manner. I am confident that this can be
done, and that MDL will continue to play an important role in the
development of statistical, and more generally, inductive inference.
\paragraph{Further Reading}
MDL can be found on the web at {\tt www.mdl-research.org}.
Good places to start further exploration of MDL are \cite{BarronRY98}
and \cite{HansenY01}. Both papers provide excellent introductions, but
they are geared towards a more specialized audience of information
theorists and statisticians, respectively. Also worth reading is  Rissanen's \citeyear{Rissanen89}
monograph.  While
outdated as an introduction to MDL {\em methods}, this famous `little
green book' still serves as a great introduction to Rissanen's radical
but appealing {\em philosophy}, which is described very eloquently.
\paragraph{Acknowledgments}
The author would like to thank Jay Myung, Mark Pitt, Steven de
Rooij and Teemu Roos, who read a preliminary version of this chapter and suggested
several improvements. 
   \begingroup
     \parindent 0pt
     \theendnotes
     \endgroup

\ \\ \ \newpage
\setlength{\textheight}{20.8 cm}
\bibliography{grunwald/bib/master,grunwald/bib/MDL,grunwald/bib/peter,grunwald/bib/thisvolume}

\begin{thebibliography}{}

\bibitem[\protect\citeauthoryear{Akaike}{Akaike}{1973}]{Akaike73}
Akaike, H. (1973).
\newblock Information theory as an extension of the maximum likelihood
  principle.
\newblock In B.~N. Petrov and F.~Csaki (Eds.), {\em Second International
  Symposium on Information Theory}, Budapest, pp.\  267--281. Akademiai Kiado.

\bibitem[\protect\citeauthoryear{Allen, Madani, and Greiner}{Allen
  et~al.}{2003}]{VanAllenMG03}
Allen, T.~V., O.~Madani, and R.~Greiner (2003, July).
\newblock Comparing model selection criteria for belief networks.
\newblock Under submission.

\bibitem[\protect\citeauthoryear{Balasubramanian}{Balasubramanian}{1997}]{Bala%
subramanian97}
Balasubramanian, V. (1997).
\newblock Statistical inference, {O}ccam's {R}azor, and statistical mechanics
  on the space of probability distributions.
\newblock {\em Neural Computation\/}~{\em 9}, 349--368.

\bibitem[\protect\citeauthoryear{Balasubramanian}{Balasubramanian}{2004}]{Bala%
subramanian04}
Balasubramanian, V. (2004).
\newblock {MDL}, {B}ayesian inference and the geometry of the space of
  probability distributions.
\newblock In P.~D. Gr\"unwald, I.~J. Myung, and M.~A. Pitt (Eds.), {\em
  Advances in Minimum Description Length: Theory and Applications}. MIT Press.

\bibitem[\protect\citeauthoryear{Barron}{Barron}{1990}]{Barron91}
Barron, A. (1990).
\newblock Complexity regularization with application to artificial neural
  networks.
\newblock In G.~Roussas (Ed.), {\em Nonparametric Functional Estimation and
  Related Topics}, Dordrecht, pp.\  561--576. Kluwer Academic Publishers.

\bibitem[\protect\citeauthoryear{Barron and Cover}{Barron and
  Cover}{1991}]{BarronC91}
Barron, A. and T.~Cover (1991).
\newblock Minimum complexity density estimation.
\newblock {\em IEEE Transactions on Information Theory\/}~{\em 37\/}(4),
  1034--1054.

\bibitem[\protect\citeauthoryear{Barron}{Barron}{1985}]{Barron85}
Barron, A.~R. (1985).
\newblock {\em Logically Smooth Density Estimation}.
\newblock Ph.\ D. thesis, Department of EE, Stanford University, Stanford, Ca.

\bibitem[\protect\citeauthoryear{Barron, Rissanen, and Yu}{Barron
  et~al.}{1998}]{BarronRY98}
Barron, A.~R., J.~Rissanen, and B.~Yu (1998).
\newblock The {M}inimum {D}escription {L}ength {P}rinciple in coding and
  modeling.
\newblock {\em IEEE Transactions on Information Theory\/}~{\em 44\/}(6),
  2743--2760.
\newblock Special Commemorative Issue: Information Theory: 1948-1998.

\bibitem[\protect\citeauthoryear{Bernardo and Smith}{Bernardo and
  Smith}{1994}]{BernardoS94}
Bernardo, J. and A.~Smith (1994).
\newblock {\em Bayesian theory}.
\newblock John Wiley.

\bibitem[\protect\citeauthoryear{Burnham and Anderson}{Burnham and
  Anderson}{2002}]{BurnhamA02}
Burnham, K.~P. and D.~R. Anderson (2002).
\newblock {\em Model Selection and Multimodel Inference}.
\newblock New York: Springer-Verlag.

\bibitem[\protect\citeauthoryear{Casella and Berger}{Casella and
  Berger}{1990}]{CasellaB90}
Casella, G. and R.~Berger (1990).
\newblock {\em Statistical Inference}.
\newblock Wadsworth.

\bibitem[\protect\citeauthoryear{Chaitin}{Chaitin}{1966}]{Chaitin66}
Chaitin, G. (1966).
\newblock On the length of programs for computing finite binary sequences.
\newblock {\em Journal of the ACM\/}~{\em 13}, 547--569.

\bibitem[\protect\citeauthoryear{Chaitin}{Chaitin}{1969}]{Chaitin69}
Chaitin, G. (1969).
\newblock On the length of programs for computing finite binary sequences:
  statistical considerations.
\newblock {\em Journal of the ACM\/}~{\em 16}, 145--159.

\bibitem[\protect\citeauthoryear{Clarke}{Clarke}{2002}]{Clarke02}
Clarke, B. (2002).
\newblock Comparing {B}ayes and non-{B}ayes model averaging when model
  approximation error cannot be ignored.
\newblock Under submission.

\bibitem[\protect\citeauthoryear{Comley and Dowe}{Comley and
  Dowe}{2004}]{ComleyD04}
Comley, J.~W. and D.~L. Dowe (2004).
\newblock {M}inimum {M}essage {L}ength and generalised {B}ayesian nets with
  asymmetric languages.
\newblock In P.~D. Gr\"unwald, I.~J. Myung, and M.~A. Pitt (Eds.), {\em
  Advances in Minimum Description Length: Theory and Applications}. MIT Press.

\bibitem[\protect\citeauthoryear{Cover and Thomas}{Cover and
  Thomas}{1991}]{CoverT91}
Cover, T. and J.~Thomas (1991).
\newblock {\em Elements of Information Theory}.
\newblock New York: Wiley Interscience.

\bibitem[\protect\citeauthoryear{Csisz{\'a}r and Shields}{Csisz{\'a}r and
  Shields}{2000}]{CsiszarS00}
Csisz{\'a}r, I. and P.~Shields (2000).
\newblock The consistency of the {BIC} {M}arkov order estimator.
\newblock {\em The Annals of Statistics\/}~{\em 28}, 1601--1619.

\bibitem[\protect\citeauthoryear{Dawid}{Dawid}{1984}]{Dawid84}
Dawid, A. (1984).
\newblock Present position and potential developments: Some personal views,
  statistical theory, the prequential approach.
\newblock {\em Journal of the Royal Statistical Society, Series {A}\/}~{\em
  147\/}(2), 278--292.

\bibitem[\protect\citeauthoryear{Dawid}{Dawid}{1992}]{Dawid92}
Dawid, A. (1992).
\newblock Prequential analysis, stochastic complexity and {B}ayesian inference.
\newblock In J.~Bernardo, J.~Berger, A.~Dawid, and A.~Smith (Eds.), {\em
  {B}ayesian Statistics}, Volume~4, pp.\  109--125. Oxford University Press.
\newblock Proceedings of the Fourth Valencia Meeting.

\bibitem[\protect\citeauthoryear{Dawid}{Dawid}{1997}]{Dawid97}
Dawid, A. (1997).
\newblock Prequential analysis.
\newblock In S.~Kotz, C.~Read, and D.~Banks (Eds.), {\em Encyclopedia of
  Statistical Sciences}, Volume 1 (Update), pp.\  464--470. Wiley-Interscience.

\bibitem[\protect\citeauthoryear{Dawid and Vovk}{Dawid and
  Vovk}{1999}]{DawidV99}
Dawid, A.~P. and V.~G. Vovk (1999).
\newblock Prequential probability: Principles and properties.
\newblock {\em Bernoulli\/}~{\em 5}, 125--162.

\bibitem[\protect\citeauthoryear{{D}e Finetti}{{D}e
  Finetti}{1974}]{DeFinetti74}
{D}e Finetti, B. (1974).
\newblock {\em Theory of Probability. A critical introductory treatment}.
\newblock London: John Wiley \& Sons.

\bibitem[\protect\citeauthoryear{Domingos}{Domingos}{1999}]{Domingos99}
Domingos, P. (1999).
\newblock The role of {O}ccam's razor in knowledge discovery.
\newblock {\em Data Mining and Knowledge Discovery\/}~{\em 3\/}(4), 409--425.

\bibitem[\protect\citeauthoryear{Feder}{Feder}{1986}]{Feder86}
Feder, M. (1986).
\newblock Maximum entropy as a special case of the minimum description length
  criterion.
\newblock {\em IEEE Transactions on Information Theory\/}~{\em 32\/}(6),
  847--849.

\bibitem[\protect\citeauthoryear{Feller}{Feller}{1968}]{Feller68a}
Feller, W. (1968).
\newblock {\em An Introduction to Probability Theory and Its Applications},
  Volume~1.
\newblock Wiley.
\newblock Third edition.

\bibitem[\protect\citeauthoryear{Foster and Stine}{Foster and
  Stine}{1999}]{FosterS99}
Foster, D. and R.~Stine (1999).
\newblock Local asymptotic coding and the minimum description length.
\newblock {\em IEEE Transactions on Information Theory\/}~{\em 45}, 1289--1293.

\bibitem[\protect\citeauthoryear{Foster and Stine}{Foster and
  Stine}{2001}]{FosterS01}
Foster, D. and R.~Stine (2001).
\newblock The competitive complexity ratio.
\newblock In {\em Proceedings of the 2001 Conference on Information Sciences
  and Systems}.
\newblock WP8 1-6.

\bibitem[\protect\citeauthoryear{Foster and Stine}{Foster and
  Stine}{2004}]{FosterS04}
Foster, D.~P. and R.~A. Stine (2004).
\newblock The contribution of parameters to stochastic complexity.
\newblock In P.~D. Gr\"unwald, I.~J. Myung, and M.~A. Pitt (Eds.), {\em
  Advances in Minimum Description Length: Theory and Applications}. MIT Press.

\bibitem[\protect\citeauthoryear{G\'acs, Tromp, and Vit\'{a}nyi}{G\'acs
  et~al.}{2001}]{GacsTV01}
G\'acs, P., J.~Tromp, and P.~Vit\'{a}nyi (2001).
\newblock Algorithmic statistics.
\newblock {\em IEEE Transactions on Information Theory\/}~{\em 47\/}(6),
  2464--2479.

\bibitem[\protect\citeauthoryear{Gr\"unwald}{Gr\"unwald}{1996}]{Grunwald96d}
Gr\"unwald, P.~D. (1996).
\newblock A minimum description length approach to grammar inference.
\newblock In G.~S. S.~Wermter, E.~Riloff (Ed.), {\em Connectionist, Statistical
  and Symbolic Approaches to Learning for Natural Language Processing}, Number
  1040 in Springer Lecture Notes in Artificial Intelligence, pp.\  203--216.

\bibitem[\protect\citeauthoryear{Gr\"unwald}{Gr\"unwald}{1998}]{Grunwald98b}
Gr\"unwald, P.~D. (1998).
\newblock {\em The Minimum Description Length Principle and Reasoning under
  Uncertainty}.
\newblock Ph.\ D. thesis, University of Amsterdam, The Netherlands.
\newblock Available as ILLC Dissertation Series 1998-03.

\bibitem[\protect\citeauthoryear{Gr\"unwald}{Gr\"unwald}{1999}]{Grunwald99a}
Gr\"unwald, P.~D. (1999).
\newblock Viewing all models as `probabilistic'.
\newblock In {\em Proceedings of the Twelfth Annual Workshop on Computational
  Learning Theory (COLT' 99)}, pp.\  171--182.

\bibitem[\protect\citeauthoryear{Gr\"unwald}{Gr\"unwald}{2000}]{Grunwald00a}
Gr\"unwald, P.~D. (2000).
\newblock Maximum entropy and the glasses you are looking through.
\newblock In {\em Proceedings of the Sixteenth Conference on Uncertainty in
  Artificial Intelligence (UAI 2000)}, pp.\  238--246. Morgan Kaufmann
  Publishers.

\bibitem[\protect\citeauthoryear{Gr\"unwald and Dawid}{Gr\"unwald and
  Dawid}{2004}]{GrunwaldD04}
Gr\"unwald, P.~D. and A.~P. Dawid (2004).
\newblock Game theory, maximum entropy, minimum discrepancy, and robust
  {Bayesian} decision theory.
\newblock {\em Annals of Statistics\/}~{\em 32\/}(4).

\bibitem[\protect\citeauthoryear{Gr\"unwald and Langford}{Gr\"unwald and
  Langford}{2004}]{GrunwaldL04}
Gr\"unwald, P.~D. and J.~Langford (2004).
\newblock Suboptimal behaviour of {B}ayes and {MDL} in classification under
  misspecification.
\newblock In {\em Proceedings of the Seventeenth Annual Conference on
  Computational Learning Theory (COLT' 04)}.

\bibitem[\protect\citeauthoryear{Gr\"unwald, Myung, and Pitt}{Gr\"unwald
  et~al.}{2004}]{GrunwaldMP04}
Gr\"unwald, P.~D., I.~J. Myung, and M.~A. Pitt (Eds.) (2004).
\newblock {\em Advances in Minimum Description Length: Theory and
  Applications}.
\newblock MIT Press.

\bibitem[\protect\citeauthoryear{Hansen and Yu}{Hansen and
  Yu}{2000}]{HansenY00}
Hansen, M. and B.~Yu (2000).
\newblock Wavelet thresholding via {MDL} for natural images.
\newblock {\em IEEE Transactions on Information Theory\/}~{\em 46}, 1778--1788.

\bibitem[\protect\citeauthoryear{Hansen and Yu}{Hansen and
  Yu}{2001}]{HansenY01}
Hansen, M. and B.~Yu (2001).
\newblock Model selection and the principle of minimum description length.
\newblock {\em Journal of the American Statistical Association\/}~{\em
  96\/}(454), 746--774.

\bibitem[\protect\citeauthoryear{Hanson and Fu}{Hanson and
  Fu}{2004}]{HansonF04}
Hanson, A.~J. and P.~C.-W. Fu (2004).
\newblock Applications of {MDL} to selected families of models.
\newblock In P.~D. Gr\"unwald, I.~J. Myung, and M.~A. Pitt (Eds.), {\em
  Advances in Minimum Description Length: Theory and Applications}. MIT Press.

\bibitem[\protect\citeauthoryear{Hjorth}{Hjorth}{1982}]{Hjorth82}
Hjorth, U. (1982).
\newblock Model selection and forward validation.
\newblock {\em Scandinavian Journal of Statistics\/}~{\em 9}, 95--105.

\bibitem[\protect\citeauthoryear{Jaynes}{Jaynes}{2003}]{Jaynes03}
Jaynes, E. (2003).
\newblock {\em Probability Theory: the logic of science}.
\newblock Cambridge University Press.
\newblock Edited by G. Larry Bretthorst.

\bibitem[\protect\citeauthoryear{Jeffreys}{Jeffreys}{1946}]{Jeffreys46}
Jeffreys, H. (1946).
\newblock An invariant form for the prior probability in estimation problems.
\newblock {\em Proceedings of the Royal Statistical Society ({L}ondon) Series
  A\/}~{\em 186}, 453--461.

\bibitem[\protect\citeauthoryear{Jeffreys}{Jeffreys}{1961}]{Jeffreys61}
Jeffreys, H. (1961).
\newblock {\em Theory of Probability\/} (Third ed.).
\newblock London: Oxford University Press.

\bibitem[\protect\citeauthoryear{Kass and Raftery}{Kass and
  Raftery}{1995}]{KassR95}
Kass, R. and A.~E. Raftery (1995).
\newblock {Bayes} factors.
\newblock {\em Journal of the American Statistical Association\/}~{\em
  90\/}(430), 773--795.

\bibitem[\protect\citeauthoryear{Kass and Voss}{Kass and Voss}{1997}]{KassV97}
Kass, R. and P.~Voss (1997).
\newblock {\em Geometrical Foundations of Asymptotic Inference}.
\newblock Wiley Interscience.

\bibitem[\protect\citeauthoryear{Kass and Wasserman}{Kass and
  Wasserman}{1996}]{KassW96}
Kass, R. and L.~Wasserman (1996).
\newblock The selection of prior distributions by formal rules.
\newblock {\em Journal of the American Statistical Association\/}~{\em 91},
  1343--1370.

\bibitem[\protect\citeauthoryear{Kearns, Mansour, Ng, and Ron}{Kearns
  et~al.}{1997}]{KearnsMNR97}
Kearns, M., Y.~Mansour, A.~Ng, and D.~Ron (1997).
\newblock An experimental and theoretical comparison of model selection
  methods.
\newblock {\em Machine Learning\/}~{\em 27}, 7--50.

\bibitem[\protect\citeauthoryear{Kolmogorov}{Kolmogorov}{1965}]{Kolmogorov65}
Kolmogorov, A. (1965).
\newblock Three approaches to the quantitative definition of information.
\newblock {\em Problems Inform. Transmission\/}~{\em 1\/}(1), 1--7.

\bibitem[\protect\citeauthoryear{Kontkanen, Myllym\"aki, Buntine, Rissanen, and
  Tirri}{Kontkanen et~al.}{2004}]{KontkanenMBRT04}
Kontkanen, P., P.~Myllym\"aki, W.~Buntine, J.~Rissanen, and H.~Tirri (2004).
\newblock An {MDL} framework for data clustering.
\newblock In P.~D. Gr\"unwald, I.~J. Myung, and M.~A. Pitt (Eds.), {\em
  Advances in Minimum Description Length: Theory and Applications}. MIT Press.

\bibitem[\protect\citeauthoryear{Kontkanen, Myllym{\"a}ki, Silander, and
  Tirri}{Kontkanen et~al.}{1999}]{KontkanenMST99}
Kontkanen, P., P.~Myllym{\"a}ki, T.~Silander, and H.~Tirri (1999).
\newblock On supervised selection of {B}ayesian networks.
\newblock In K.~Laskey and H.~Prade (Eds.), {\em Proceedings of the 15th
  International Conference on Uncertainty in Artificial Intelligence (UAI'99)}.
  Morgan Kaufmann Publishers.

\bibitem[\protect\citeauthoryear{Lanterman}{Lanterman}{2004}]{Lanterman04}
Lanterman, A.~D. (2004).
\newblock Hypothesis testing for {P}oisson versus geometric distributions using
  stochastic complexity.
\newblock In P.~D. Gr\"unwald, I.~J. Myung, and M.~A. Pitt (Eds.), {\em
  Advances in Minimum Description Length: Theory and Applications}. MIT Press.

\bibitem[\protect\citeauthoryear{Lee}{Lee}{1997}]{Lee97}
Lee, P. (1997).
\newblock {\em Bayesian Statistics -- an introduction}.
\newblock Arnold \& Oxford University Press.

\bibitem[\protect\citeauthoryear{Li, Chen, Li, Ma, and Vit\'anyi}{Li
  et~al.}{2003}]{LiCLV03}
Li, M., X.~Chen, X.~Li, B.~Ma, and P.~Vit\'anyi (2003).
\newblock The similarity metric.
\newblock In {\em Proc. 14th ACM-SIAM Symp. Discrete Algorithms (SODA)}.

\bibitem[\protect\citeauthoryear{Li and Vit\'anyi}{Li and
  Vit\'anyi}{1997}]{LiV97}
Li, M. and P.~Vit\'anyi (1997).
\newblock {\em An Introduction to {K}olmogorov Complexity and Its
  Applications\/} (revised and expanded second ed.).
\newblock New York: Springer-Verlag.

\bibitem[\protect\citeauthoryear{Liang and Barron}{Liang and
  Barron}{2004a}]{LiangB04}
Liang, F. and A.~Barron (2004a).
\newblock Exact minimax predictive density estimation and {MDL}.
\newblock In P.~D. Gr\"unwald, I.~J. Myung, and M.~A. Pitt (Eds.), {\em
  Advances in Minimum Description Length: Theory and Applications}. MIT Press.

\bibitem[\protect\citeauthoryear{Liang and Barron}{Liang and
  Barron}{2004b}]{LiangB03}
Liang, F. and A.~Barron (2004b).
\newblock Exact minimax strategies for predictive density estimation.
\newblock To appear in {\em IEEE Transactions on Information Theory}.

\bibitem[\protect\citeauthoryear{Modha and Masry}{Modha and
  Masry}{1998}]{ModhaM98}
Modha, D.~S. and E.~Masry (1998).
\newblock Prequential and cross-validated regression estimation.
\newblock {\em Machine Learning\/}~{\em 33\/}(1), 5--39.

\bibitem[\protect\citeauthoryear{Myung, Balasubramanian, and Pitt}{Myung
  et~al.}{2000}]{MyungBP00}
Myung, I.~J., V.~Balasubramanian, and M.~A. Pitt (2000).
\newblock Counting probability distributions: Differential geometry and model
  selection.
\newblock {\em Proceedings of the National Academy of Sciences USA\/}~{\em 97},
  11170--11175.

\bibitem[\protect\citeauthoryear{Navarro}{Navarro}{2004}]{Navarro04}
Navarro, D. (2004).
\newblock Misbehaviour of the {F}isher information approximation to {M}inimum
  {D}escription {L}ength.
\newblock Under submission.

\bibitem[\protect\citeauthoryear{Pednault}{Pednault}{2003}]{Pednault03}
Pednault, E. (2003).
\newblock Personal communication.

\bibitem[\protect\citeauthoryear{Quinlan and Rivest}{Quinlan and
  Rivest}{1989}]{QuinlanR89}
Quinlan, J. and R.~Rivest (1989).
\newblock Inferring decision trees using the minimum description length
  principle.
\newblock {\em Information and Computation\/}~{\em 80}, 227--248.

\bibitem[\protect\citeauthoryear{Ripley}{Ripley}{1996}]{Ripley96}
Ripley, B. (1996).
\newblock {\em Pattern Recognition and Neural Networks}.
\newblock Cambridge University Press.

\bibitem[\protect\citeauthoryear{Rissanen}{Rissanen}{1978}]{Rissanen78}
Rissanen, J. (1978).
\newblock Modeling by the shortest data description.
\newblock {\em Automatica\/}~{\em 14}, 465--471.

\bibitem[\protect\citeauthoryear{Rissanen}{Rissanen}{1983}]{Rissanen83}
Rissanen, J. (1983).
\newblock A universal prior for integers and estimation by minimum description
  length.
\newblock {\em The Annals of Statistics\/}~{\em 11}, 416--431.

\bibitem[\protect\citeauthoryear{Rissanen}{Rissanen}{1984}]{Rissanen84}
Rissanen, J. (1984).
\newblock Universal coding, information, prediction and estimation.
\newblock {\em IEEE Transactions on Information Theory\/}~{\em 30}, 629--636.

\bibitem[\protect\citeauthoryear{Rissanen}{Rissanen}{1986}]{Rissanen86a}
Rissanen, J. (1986).
\newblock Stochastic complexity and modeling.
\newblock {\em The Annals of Statistics\/}~{\em 14}, 1080--1100.

\bibitem[\protect\citeauthoryear{Rissanen}{Rissanen}{1987}]{Rissanen87}
Rissanen, J. (1987).
\newblock Stochastic complexity.
\newblock {\em Journal of the Royal Statistical Society, series B\/}~{\em 49},
  223--239.
\newblock Discussion: pages 252--265.

\bibitem[\protect\citeauthoryear{Rissanen}{Rissanen}{1989}]{Rissanen89}
Rissanen, J. (1989).
\newblock {\em Stochastic Complexity in Statistical Inquiry}.
\newblock World Scientific Publishing Company.

\bibitem[\protect\citeauthoryear{Rissanen}{Rissanen}{1996}]{Rissanen96}
Rissanen, J. (1996).
\newblock Fisher information and stochastic complexity.
\newblock {\em IEEE Transactions on Information Theory\/}~{\em 42\/}(1),
  40--47.

\bibitem[\protect\citeauthoryear{Rissanen}{Rissanen}{2000}]{Rissanen00}
Rissanen, J. (2000).
\newblock {MDL} denoising.
\newblock {\em IEEE Transactions on Information Theory\/}~{\em 46\/}(7),
  2537--2543.

\bibitem[\protect\citeauthoryear{Rissanen}{Rissanen}{2001}]{Rissanen01}
Rissanen, J. (2001).
\newblock Strong optimality of the normalized {ML} models as universal codes
  and information in data.
\newblock {\em IEEE Transactions on Information Theory\/}~{\em 47\/}(5),
  1712--1717.

\bibitem[\protect\citeauthoryear{Rissanen, Speed, and Yu}{Rissanen
  et~al.}{1992}]{RissanenSY92}
Rissanen, J., T.~Speed, and B.~Yu (1992).
\newblock Density estimation by stochastic complexity.
\newblock {\em IEEE Transactions on Information Theory\/}~{\em 38\/}(2),
  315--323.

\bibitem[\protect\citeauthoryear{Rissanen and Tabus}{Rissanen and
  Tabus}{2004}]{RissanenT04}
Rissanen, J. and I.~Tabus (2004).
\newblock {K}olmogorov's structure function in {MDL} theory and lossy data
  compression.
\newblock In P.~D. Gr\"unwald, I.~J. Myung, and M.~A. Pitt (Eds.), {\em
  Advances in Minimum Description Length: Theory and Applications}. MIT Press.

\bibitem[\protect\citeauthoryear{Schwarz}{Schwarz}{1978}]{Schwarz78}
Schwarz, G. (1978).
\newblock Estimating the dimension of a model.
\newblock {\em The Annals of Statistics\/}~{\em 6\/}(2), 461--464.

\bibitem[\protect\citeauthoryear{Shafer and Vovk}{Shafer and
  Vovk}{2001}]{ShaferV01}
Shafer, G. and V.~Vovk (2001).
\newblock {\em Probability and Finance -- It's only a game!}
\newblock Wiley.

\bibitem[\protect\citeauthoryear{Shtarkov}{Shtarkov}{1987}]{Shtarkov87}
Shtarkov, Y.~M. (1987).
\newblock Universal sequential coding of single messages.
\newblock {\em (translated from) Problems of Information Transmission\/}~{\em
  23\/}(3), 3--17.

\bibitem[\protect\citeauthoryear{Solomonoff}{Solomonoff}{1964}]{Solomonoff64}
Solomonoff, R. (1964).
\newblock A formal theory of inductive inference, part 1 and part 2.
\newblock {\em Information and Control\/}~{\em 7}, 1--22, 224--254.

\bibitem[\protect\citeauthoryear{Solomonoff}{Solomonoff}{1978}]{Solomonoff78}
Solomonoff, R. (1978).
\newblock Complexity-based induction systems: comparisons and convergence
  theorems.
\newblock {\em IEEE Transactions on Information Theory\/}~{\em 24}, 422--432.

\bibitem[\protect\citeauthoryear{Speed and Yu}{Speed and Yu}{1993}]{SpeedY93}
Speed, T. and B.~Yu (1993).
\newblock Model selection and prediction: normal regression.
\newblock {\em Ann. Inst. Statist. Math.\/}~{\em 45\/}(1), 35--54.

\bibitem[\protect\citeauthoryear{Takeuchi}{Takeuchi}{2000}]{Takeuchi00}
Takeuchi, J. (2000).
\newblock On minimax regret with respect to families of stationary stochastic
  processes (in {J}apanese).
\newblock In {\em Proceedings IBIS 2000}, pp.\  63--68.

\bibitem[\protect\citeauthoryear{Takeuchi and Barron}{Takeuchi and
  Barron}{1997}]{TakeuchiB97}
Takeuchi, J. and A.~Barron (1997).
\newblock Asymptotically minimax regret for exponential families.
\newblock In {\em Proceedings SITA '97}, pp.\  665--668.

\bibitem[\protect\citeauthoryear{Takeuchi and Barron}{Takeuchi and
  Barron}{1998}]{TakeuchiB98}
Takeuchi, J. and A.~Barron (1998).
\newblock Asymptotically minimax regret by {B}ayes mixtures.
\newblock In {\em Proceedings of the 1998 International Symposium on
  Information Theory (ISIT 98)}.

\bibitem[\protect\citeauthoryear{Townsend}{Townsend}{1975}]{Townsend75}
Townsend, P. (1975).
\newblock The mind-body equation revisited.
\newblock In C.-Y. Cheng (Ed.), {\em Psychological Problems in Philosophy},
  pp.\  200--218. Honolulu: University of Hawaii Press.

\bibitem[\protect\citeauthoryear{Vapnik}{Vapnik}{1998}]{Vapnik98}
Vapnik, V. (1998).
\newblock {\em Statistical Learning Theory.}
\newblock John Wiley.

\bibitem[\protect\citeauthoryear{Vereshchagin and Vit\'anyi}{Vereshchagin and
  Vit\'anyi}{2002}]{VereshchaginV02}
Vereshchagin, N. and P.~M.~B. Vit\'anyi (2002).
\newblock Kolmogorov's structure functions with an application to the
  foundations of model selection.
\newblock In {\em Proc. 47th IEEE Symp. Found. Comput. Sci. (FOCS'02)}.

\bibitem[\protect\citeauthoryear{Viswanathan., Wallace, Dowe, and
  Korb}{Viswanathan. et~al.}{1999}]{ViswanathanWDK99}
Viswanathan., M., C.~Wallace, D.~Dowe, and K.~Korb (1999).
\newblock Finding cutpoints in noisy binary sequences - a revised empirical
  evaluation.
\newblock In {\em Proc. 12th Australian Joint Conf. on Artif. Intelligence},
  Volume 1747 of {\em Lecture Notes in Artificial Intelligence (LNAI)}, Sidney,
  Australia, pp.\  405--416.

\bibitem[\protect\citeauthoryear{Vit\'anyi}{Vit\'anyi}{2004}]{Vitanyi04}
Vit\'anyi, P.~M. (2004).
\newblock Algorithmic statistics and {K}olmogorov's structure function.
\newblock In P.~D. Gr\"unwald, I.~J. Myung, and M.~A. Pitt (Eds.), {\em
  Advances in Minimum Description Length: Theory and Applications}. MIT Press.

\bibitem[\protect\citeauthoryear{Wallace and Boulton}{Wallace and
  Boulton}{1968}]{WallaceB68}
Wallace, C. and D.~Boulton (1968).
\newblock An information measure for classification.
\newblock {\em Computing Journal\/}~{\em 11}, 185--195.

\bibitem[\protect\citeauthoryear{Wallace and Boulton}{Wallace and
  Boulton}{1975}]{WallaceB75}
Wallace, C. and D.~Boulton (1975).
\newblock An invariant {B}ayes method for point estimation.
\newblock {\em Classification Society Bulletin\/}~{\em 3\/}(3), 11--34.

\bibitem[\protect\citeauthoryear{Wallace and Freeman}{Wallace and
  Freeman}{1987}]{WallaceF87}
Wallace, C. and P.~Freeman (1987).
\newblock Estimation and inference by compact coding.
\newblock {\em Journal of the Royal Statistical Society, Series {B}\/}~{\em
  49}, 240--251.
\newblock Discussion: pages 252--265.

\bibitem[\protect\citeauthoryear{Webb}{Webb}{1996}]{Webb96}
Webb, G. (1996).
\newblock Further experimental evidence against the utility of {O}ccam's razor.
\newblock {\em Journal of Artificial Intelligence Research\/}~{\em 4},
  397--417.

\bibitem[\protect\citeauthoryear{Yamanishi}{Yamanishi}{1998}]{Yamanishi98}
Yamanishi, K. (1998).
\newblock A decision-theoretic extension of stochastic complexity and its
  applications to learning.
\newblock {\em IEEE Transactions on Information Theory\/}~{\em 44\/}(4),
  1424--1439.

\bibitem[\protect\citeauthoryear{Zhang}{Zhang}{2004}]{Zhang04}
Zhang, T. (2004).
\newblock On the convergence of {MDL} density estimation.
\newblock In Y.~Singer and J.~Shawe-Taylor (Eds.), {\em Proceedings of the
  Seventeenth Annual Conference on Computational Learning Theory (COLT' 04)},
  Lecture Notes in Computer Science. Springer-Verlag.

\end{thebibliography}

\cleardoublepage
%\addcontentsline{toc}{chapter}{\protect\numberline{}\protect\hbox{Index}}
%\printindex

\end{document}